\documentclass[12pt]{amsart}
\usepackage{fullpage}
\usepackage{hyperref,mathrsfs,stmaryrd,amssymb}
\usepackage[all]{xy}

\newtheorem{theorem}{Theorem}[subsection]
\newtheorem{lemma}[theorem]{Lemma}
\newtheorem{lemma-def}[theorem]{Lemma-Definition}
\newtheorem{prop}[theorem]{Proposition}
\newtheorem{cor}[theorem]{Corollary}

\theoremstyle{definition}
\newtheorem{hypothesis}[theorem]{Hypothesis}

\newtheorem{notation}[theorem]{Notation}
\newtheorem{remark}[theorem]{Remark}
\newtheorem{defn}[theorem]{Definition}
\newtheorem{example}[theorem]{Example}

\numberwithin{equation}{theorem}

\newcommand{\AAA}{\mathbf{A}}
\newcommand{\CC}{\mathbb{C}}
\newcommand{\FF}{\mathbb{F}}

\newcommand{\PP}{\mathbf{P}}
\newcommand{\QQ}{\mathbb{Q}}
\newcommand{\RR}{\mathbb{R}}
\newcommand{\ZZ}{\mathbb{Z}}
\newcommand{\calC}{\mathcal{C}}

\newcommand{\calE}{\mathcal{E}}
\newcommand{\calF}{\mathcal{F}}
\newcommand{\calG}{\mathcal{G}}
\newcommand{\calH}{\mathcal{H}}

\newcommand{\calL}{\mathcal{L}}
\newcommand{\calM}{\mathcal{M}}
\newcommand{\calO}{\mathcal{O}}

\newcommand{\calX}{\mathcal{X}}
\newcommand{\calY}{\mathcal{Y}}

\newcommand{\bv}{\mathbf{v}}

\newcommand{\frako}{\mathfrak{o}}
\newcommand{\frakC}{\mathfrak{C}}
\newcommand{\frakD}{\mathfrak{D}}
\newcommand{\frakE}{\mathfrak{E}}
\newcommand{\frakf}{\mathfrak{f}}
\newcommand{\frakF}{\mathfrak{F}}

\newcommand{\frakH}{\mathfrak{H}}

\newcommand{\frakS}{\mathfrak{S}}

\newcommand{\frakX}{\mathfrak{X}}

\newcommand{\frakZ}{\mathfrak{Z}}
\newcommand{\dual}{\vee}

\newcommand{\stacktag}[1]{\cite[\href{https://stacks.math.columbia.edu/tag/#1}{Tag #1}]{stacks-project}}

\DeclareMathOperator{\ab}{ab}

\DeclareMathOperator{\Coh}{\mathbf{Coh}}
\DeclareMathOperator{\CohHom}{\mathbf{CohHom}}

\DeclareMathOperator{\Conn}{\mathbf{Conn}}
\DeclareMathOperator{\ConnHom}{\mathbf{ConnHom}}
\DeclareMathOperator{\coker}{coker}

\DeclareMathOperator{\dR}{dR}
\DeclareMathOperator{\End}{End}
\DeclareMathOperator{\et}{et}
\DeclareMathOperator{\Ext}{Ext}
\DeclareMathOperator{\Frac}{Frac}
\DeclareMathOperator{\FConn}{\mathbf{F-Conn}}

\DeclareMathOperator{\FIsoc}{\mathbf{F-Isoc}}
\DeclareMathOperator{\Frob}{Frob}

\DeclareMathOperator{\GL}{GL}

\DeclareMathOperator{\Hom}{Hom}

\DeclareMathOperator{\Isoc}{\mathbf{Isoc}}

\DeclareMathOperator{\perf}{perf}
\DeclareMathOperator{\Pic}{Pic}
\DeclareMathOperator{\PGL}{PGL}
\DeclareMathOperator{\Quot}{\mathbf{Quot}}

\DeclareMathOperator{\rank}{rank}
\DeclareMathOperator{\red}{red}
\DeclareMathOperator{\Rep}{\mathbf{Rep}}

\DeclareMathOperator{\Sch}{\mathbf{Sch}}

\DeclareMathOperator{\Spec}{Spec}
\DeclareMathOperator{\Spf}{Spf}

\DeclareMathOperator{\tame}{tame}

\DeclareMathOperator{\Vect}{\mathbf{Vec}}
\DeclareMathOperator{\Weil}{\mathbf{Weil}}
\DeclareMathOperator{\width}{width}

\begin{document}
\title{\'Etale and crystalline companions, II}
\author{Kiran S. Kedlaya}
\date{March 24, 2026}
\thanks{Thanks to Tomoyuki Abe, Anna Cadoret, Marco D'Addezio, Johan de Jong, Pierre Deligne,
Valentina Di Proietto, H\'el\`ene Esnault, Thomas Grubb,
Teruhisa Koshikawa,
Raju Krishna\-moorthy,
Daniel Litt, Gyujin Oh,
Ambrus P\'al, and
Remy van Dobben de Bruyn for helpful discussions;
thanks also to Cadoret for providing an advance copy of \cite{cadoret}, and to an anonymous referee for identifying some significant errors in an earlier version.
The author was supported by NSF (grants DMS-1501214, DMS-1802161, DMS-2053473), UCSD (Warschawski Professorship),
and the IAS School of Mathematics (Visiting Professorship 2018--2019).
The author also benefited from the hospitality of IM PAN (Simons Foundation grant 346300, Polish MNiSW fund 2015--2019) during September 2018 and MSRI (NSF grants DMS-1440140, DMS-1928930) during May 2019 and January 2023.}

\begin{abstract}
Let $X$ be a smooth scheme over a finite field of characteristic $p$.
In answer to a conjecture of Deligne, we establish that for any prime $\ell \neq p$,
an $\ell$-adic Weil sheaf on $X$ which is algebraic (or irreducible with finite determinant)
admits a \emph{crystalline companion} in the category of overconvergent $F$-isocrystals, for which the Frobenius characteristic polynomials agree at all closed points (with respect to some fixed identification of the algebraic closures of $\QQ$ within fixed algebraic closures of $\QQ_\ell$ and $\QQ_p$).
The argument depends heavily on the free passage between $\ell$-adic and $p$-adic coefficients for curves provided by the Langlands correspondence for $\mathrm{GL}_n$ over global function fields (work of L. Lafforgue and T. Abe), and
on the construction of Drinfeld (plus adaptations by Abe--Esnault and Kedlaya) giving rise to \'etale companions
of overconvergent $F$-isocrystals.
As corollaries, we transfer a number of statements from crystalline to \'etale coefficient objects, including
properties of the Newton polygon stratification (results of Grothendieck--Katz and de Jong--Oort--Yang)
and Wan's theorem  (previously Dwork's conjecture) on $p$-adic meromorphicity of unit-root $L$-functions.
\end{abstract}

\maketitle

\section*{Introduction}

\subsection{Overview}
Let $k$ be a finite field of characteristic $p$ and let $X$ be a smooth scheme over $k$.
In a previous paper \cite{kedlaya-companions},
we studied the relationship between coefficient objects (of locally constant rank) in Weil cohomology with 
$\ell$-adic coefficients for various primes $\ell$. For $\ell \neq p$, such objects are \emph{lisse Weil $\overline{\QQ}_\ell$-sheaves}, while for $\ell = p$ they are \emph{overconvergent $F$-isocrystals}.

The purpose of this paper is to complete the proof of a conjecture of Deligne \cite[Conjecture~1.2.10]{deligne-weil2} which asserts that all coefficient objects ``look motivic'', that is, they have various features that would hold if they were to arise in the cohomology of some family of smooth proper varieties over $X$. The deepest aspect of this conjecture is the fact that coefficient objects do not occur in isolation; in a certain sense, coefficient objects in one category have ``companions'' in the other categories.

To make this more precise, we say that an $\ell$-adic coefficient $\calE$ on $X$ is \emph{algebraic} if
for each closed point $x \in X$, the characteristic polynomial of (geometric) Frobenius $F_x$ acting on the fiber $\calE_x$ has coefficients which are algebraic over $\QQ$. 

To define the companion relation, consider two primes $\ell, \ell'$ and fix
 an identification of the algebraic closures of $\QQ$ within $\overline{\QQ}_\ell$
and $\overline{\QQ}_{\ell'}$.
We say that an algebraic $\ell$-adic coefficient $\calE$ on $X$ and an algebraic $\ell'$-adic coefficient $\calE'$ on $X$ are \emph{companions} if for each closed point $x \in X$, the characteristic polynomials of Frobenius on $\calE_x, \calE'_x$ coincide; note that $\calE$ then determines $\calE'$ up to semisimplification
\cite[Theorem~3.3.1]{kedlaya-companions}.

With this definition, we can state a theorem answering \cite[Conjecture~1.2.10]{deligne-weil2}
(and more precisely \cite[Conjecture~0.1.1]{kedlaya-companions}
or \cite[Conjecture~1.1]{cadoret}),
 incorporating Crew's proposal \cite[Conjecture~4.13]{crew-mono}
to interpret (vi) by reading the phrase
``petit camarade cristalline'' to mean ``companion in the category of overconvergent $F$-isocrystals.'' (The definition of an overconvergent $F$-isocrystal was unavailable at the time of \cite{deligne-weil2}; it was subsequently introduced by Berthelot \cite{berthelot-mem}.)
See Corollary~\ref{C:deligne} for the proof.

\begin{theorem} \label{T:deligne}
Let $\calE$ be an $\ell$-adic coefficient which is irreducible with determinant of finite order.
(Recall that we allow $\ell=p$ here.)
In the following statements, $x$ is always quantified over all closed points of $X$, and $\kappa(x)$ denotes the residue field of $x$.
\begin{enumerate}
\item[(i)]
$\calE$ is \emph{pure of weight $0$}: for every algebraic embedding of $\overline{\QQ}_\ell$ into $\CC$ and all $x$, the images of the eigenvalues of $F_x$ all have complex absolute value $1$.
\item[(ii)]
For some number field $E$, 
$\calE$ is \emph{$E$-algebraic}: for all $x$, the characteristic polynomial of $F_x$ has coefficients in $E$.
(Beware that the roots of this polynomial need not belong to a single number field as $x$ varies.)
\item[(iii)]
$\calE$ is \emph{$p$-plain}: for all $x$, the eigenvalues of $F_x$ have trivial $\lambda$-adic valuation at all finite places $\lambda$ of $E$ not lying above $p$.
\item[(iv)]
For every place $\lambda$ of $E$ above $p$ and all $x$, every eigenvalue of $F_x$ has $\lambda$-adic valuation at most
$\frac{1}{2} \rank(\calE)$ times the valuation of $\#\kappa(x)$.
\item[(v)]
For every prime $\ell' \neq p$ and every place $\lambda$ of $E$ above $\ell'$, there exists an $\ell'$-adic coefficient $\calE'$ which is irreducible with determinant of finite order and is a companion of $\calE$ with respect to $\lambda$.
\item[(vi)]
As in (v), but with $\ell' = p$.
\end{enumerate}
\end{theorem}

Of the various aspects of Theorem~\ref{T:deligne}, all were previously known for $X$ of dimension 1, and
all but (vi) for general $X$ (see below); consequently, the new content can also be expressed as follows
(answering \cite[Conjecture~0.5.1]{kedlaya-companions}).
See Corollary~\ref{C:companion p} for the proof.

\begin{theorem}  \label{T:companion}
Any algebraic $\ell$-adic coefficient on $X$ admits an $\ell'$-adic companion.
(By \cite[Theorem~3.5.2]{kedlaya-companions}, this is only new for $\ell'=p$.)
\end{theorem}

In the remainder of this introduction, we summarize the preceding work in the direction of Theorem~\ref{T:deligne},
including the results of \cite{kedlaya-companions};
we then describe the new ingredients in this paper that lead to a complete proof.
See also the survey article \cite{kedlaya-isocrystals} for background on $p$-adic coefficients.

\subsection{Prior results: dimension 1}

One conceivable approach to proving Theorem~\ref{T:companion} would be to show that every $\ell$-adic coefficient arises as a realization of some motive, to which one could then apply the $\ell'$-adic realization functor to obtain the $\ell'$-adic companion. 
As part of the original formulation of \cite[Conjecture~1.2.10]{deligne-weil2}, Deligne pointed out that 
for $X$ of dimension 1, one could hope to execute this strategy for $\ell, \ell' \neq p$ by establishing the Langlands correspondence for $\mathrm{GL}_n$ over the function field of $X$ for all positive integers $n$;
at the time, this was done only for $n=1$ by class field theory, and for $n=2$ by the work of Drinfeld \cite{drinfeld}. 
An extension of Drinfeld's work to general $n$ was subsequently achieved by L. Lafforgue \cite{lafforgue}, which yields parts (i)--(v) of Theorem~\ref{T:deligne} when $\dim(X) = 1$
except with a slightly weaker inequality in part (iv); this was subsequently improved by V. Lafforgue \cite{lafforgue-hecke} to obtain (iv) as written (and a bit more).

This work necessarily omitted cases where $\ell=p$ or $\ell'=p$ due to the limited development of $p$-adic Weil cohomology (rigid cohomology) at the time. Building on recent advances in this direction, Abe \cite{abe-companion} has replicated Lafforgue's argument in $p$-adic cohomology; this completes  Theorem~\ref{T:deligne} when $\dim(X) = 1$ by adding the cases where $\ell=p$ or $\ell'=p$.

\subsection{Prior results: higher dimension}

In higher dimensions, no general method for associating motives to coefficient objects seems to be known.
The proofs of the various aspects of Theorem~\ref{T:deligne} for general $X$ thus proceed by using the case of curves as a black box. 

To begin with, suppose that $\ell \neq p$.
To make headway, one first shows that irreducibility is preserved by restriction to suitable curves;
this was shown by Deligne \cite[\S 1.7]{deligne-finite} by correcting an argument of L. Lafforgue \cite[\S VII]{lafforgue}. Since parts (i), (iii), (iv) of Theorem~\ref{T:deligne} are statements about individual closed points,
they follow almost immediately.

As for part (ii) of Theorem~\ref{T:deligne}, by restricting to curves one sees that the coefficients in question are all algebraic, but one needs a uniformity argument over these curves to show that the extension of $\QQ$ generated by all of the coefficients is finite. Such an argument was provided by Deligne \cite{deligne-finite}.
Building on this, Drinfeld \cite{drinfeld-deligne} was then able to establish part (v) of Theorem~\ref{T:deligne} using an idea of Wiesend \cite{wiesend1} to patch together tame Galois representations based on their restrictions to curves. (Esnault--Kerz \cite{esnault-kerz} refer to this technique as the method of \emph{skeleton sheaves}.)

It is not entirely automatic to extend these results to the case $\ell = p$, as several key arguments (notably preservation of irreducibility) are made in terms of properties of $\ell$-adic sheaves with no direct $p$-adic analogues. Generally, arguments that refer to residual representations do not transfer (although there are some crucial exceptions), whereas arguments that refer
only to monodromy groups or cohomology do transfer.
With some effort, one can replace the offending arguments with alternates that can be ported to the $p$-adic setting, and thus recover the previously mentioned results with $\ell=p$; this was carried out (in slightly different ways) by Abe--Esnault \cite{abe-esnault} and Kedlaya 
\cite{kedlaya-companions}.
One key example of this replacement is the reduction to the case of everywhere tame local ramification (see next paragraph): when $\ell \neq p$ this is an elementary observation, but when $\ell=p$ it requires the semistable reduction theorem for overconvergent $F$-isocrystals \cite{kedlaya-semi1, kedlaya-semi2, kedlaya-semi3, kedlaya-semi4}.

Crucially, this possibility of replacement only applies in cases where one starts with a $p$-adic coefficient;
it therefore does not apply to part (vi) of Theorem~\ref{T:deligne}, where the existence of a $p$-adic coefficient is itself at issue and cannot be established using a direct analogue of the $\ell$-adic construction.
However, from the above discussion, one can at least deduce  some reductions for the problem of constructing crystalline companions; notably, in any given case the existence of a crystalline companion can be checked after pullback along an open immersion with dense image, or an alteration in the sense of de Jong \cite{dejong-alterations}.
Also, we may ignore the case $\ell = p$, as we may move from $\ell = p$ to $\ell' = p$ via an intermediate prime different from $p$.
This means that in most of this paper, we can focus on the case of an \'etale $\ell$-adic coefficient which is \emph{tame} (that is, whose local monodromy representations are tamely ramified and quasi-unipotent) or even \emph{docile} (replacing ``quasi-unipotent'' by ``unipotent''). We can also focus most attention on the situation where $X$ embeds into a smooth curve fibration,
which provides some technical simplifications (for example, the existence of smooth lifts \'etale-locally on the base of the fibration).

\subsection{Uniformity, fibrations, and crystalline companions}

We now arrive at the methods of the present paper. Before proceeding, we point out that \S 1--3 consist entirely of background material (as does \S 4 for the most part, with the notable exception of \S\ref{subsec:restriction to divisors}); we thus jump ahead to \S\ref{sec:uniformities} for this part of the discussion.

At a superficial level, the basic strategy is the same as in Drinfeld's work. Given an \'etale coefficient object on $X$, the Langlands correspondence implies the existence of a crystalline companion for the restriction to any curve contained in $X$. In order to construct a crystalline companion on $X$ itself, we will construct a coherent sequence of mod-$p^n$ truncations whose inverse limit gives rise to the companion. This will require building and analyzing a moduli stack parametrizing all possible truncated coefficient objects; the presence of many points arising from companions on curves forces the moduli space to be large enough to give rise to an object over all of $X$. 

In order to further this analogy, we need some uniformity properties on companions on curves in order to obtain a finiteness property for the moduli spaces of truncations.
In the \'etale case,
the necessary uniformity assertion is a consequence of properties of tame \'etale fundamental groups. In the crystalline case, we instead make a careful analysis in \S\ref{sec:uniformities}
of the Harder--Narasimhan polygons of the underlying vector bundles of tame overconvergent $F$-isocrystals on curves.
 
\subsection{Companions on a curve fibration: a special case}
\label{subsec:intro companions on fibration}

We next give a more detailed account of how the aforementioned ingredients come together in \S\ref{sec:companions in fibration} and \S\ref{sec:companion points} in 
the key technical step: 
the construction of a crystalline companion of an algebraic \'etale coefficient on the total space of a smooth curve fibration over some smooth base space $S$.
For reasons explained below, we must also assume the least slope of the generic Newton polygon occurs with multiplicity 1; see the discussion in \S\ref{subsec:intro companions on fibration2} to find out how we work around this to prove our final result.

Since we will apply this construction in the context of an induction on dimension (Theorem~\ref{T:algebraic companions}), we are free to assume the existence of crystalline companions not just on curves, but on all varieties of dimension at most that of $S$ (Hypothesis~\ref{H:companion points new}). In particular, divisors on $X$ that dominate $S$ play a key role.

In \S\ref{sec:companions in fibration}, we apply the uniformity estimates from 
\S\ref{sec:uniformities} to construct some moduli stacks of mod-$p^n$ relative connections
which are of finite type over $\ZZ$ and hence noetherian (Lemma~\ref{L:uniformity for crystalline lattices}). Applying the existence of crystalline companions on curves in the fibration yields a collection of ``companion points'' on these moduli stacks, whose Zariski closures give rise to a profinite collection of candidates for the generic fiber of the desired crystalline companion (Proposition~\ref{P:carryover data}). However, while we do get some relationship between these objects and the crystalline companions on fibers, in terms of the $p$-adic local systems associated to successive quotients of the slope filtration, this control is diffused across the profinite set of candidates; we are thus stymied by the fact that this construction cannot by itself produce a \emph{single} candidate for the generic fiber.

To work around this, in \S \ref{subsec:rank 1 objects} we introduce a variant of Wiesend's construction to identify the local system corresponding to the unit root of the desired crystalline companion. This requires some care as this representation is in general highly wildly ramified; to overcome this we are forced to limit attention to the abelian case, whence our condition that the least slope of the generic Newton polygon occurs with multiplicity 1.
Using the stack-theoretic input, we are able to verify the hypotheses of this construction to obtain a candidate for the unit root of our putative crystaline companion (Proposition~\ref{P:putative unit root}).
As a corollary, we gain some more control on Newton polygons in the \'etale case (Corollary~\ref{C:open constant NP}). 

We then take advantage of the fact that 
irreducible overconvergent $F$-isocrystals are uniquely identified by the
first steps of their (convergent) slope filtrations:
this is the \emph{minimal slope theorem}
of Tsuzuki \cite{tsuzuki-minimal} and D'Addezio \cite{daddezio}. In our context,
the minimal slope theorem collapses our profinite set of candidates for the generic fiber of the crystalline companion down to a single object (Proposition~\ref{P:generic companion}).

Next, using particular divisors on $X$ which dominate $S$,
we extend the local systems corresponding to the successive quotients of the slope filtration of the generic fiber 
across the entire fibration (Corollary~\ref{C:extend slope steps}).
Finally, we again use a suitable divisor to glue together these extensions with our putative generic fiber to obtain a full crystalline companion (Proposition~\ref{P:full companion on fibration}).

\subsection{Companions on a curve fibration: general case}
\label{subsec:intro companions on fibration2}

We next explain the content of \S\ref{subsec:main results}, where we apply the argument described in \S\ref{subsec:intro companions on fibration} to complete the construction of crystalline companions.
 
The main remaining issue is to overcome the restriction that the least slope of the generic Newton polygon occurs with multiplicity 1, which was needed to adapt Wiesend's construction to handle wild ramification.
When the generic Newton polygon has more than one distinct slope,
we can take some intermediate exterior power to arrive at a case where the least slope occurs with multiplicity 1; the cost of this is that taking the exterior power may reduce the Tannakian automorphism group by a finite quotient. To recover the desired companion, we again use a suitable divisor to reduce to the trivialization of a certain  Brauer class (Lemma~\ref{L:deduce from exterior power}). 

This still leaves the case where the generic Newton polygon has only one slope, which is to say it is \emph{isoclinic}. (Note that even for a fixed \'etale coefficient, whether or not this holds may depend on the choice of a $p$-adic place; see \S\ref{subsec:Newton}.)
In this case, the desired crystalline companion corresponds to an \emph{everywhere unramified} 
$\overline{\QQ}_p$-local system; we can construct the latter using a slight modification of Drinfeld's construction in the \'etale case (Lemma~\ref{L:isoclinic companions}).

\subsection{Applications}

We conclude the paper by giving some applications of the construction of crystalline companions (for more on which see the remainder of \S\ref{sec:corollaries}).
These include properties of Newton polygons
(Theorem~\ref{T:stratification1}, Theorem~\ref{T:purity}) and 
Wan's theorem on the $p$-adic meromorphy of unit-root $L$-functions (Theorem~\ref{T:Wan}).

As remarked upon in \cite[Remark~2.1.6]{kedlaya-companions}, the existence of companions on curves suggests a new method for counting lisse Weil sheaves on curves, by directly relating these counts to the zeta functions of moduli spaces of vector bundles (see \emph{loc. cit.} for some references on this question).
Theorem~\ref{T:companion} in turn provides an opportunity (albeit one not acted upon here)
to make similar arguments on higher-dimensional varieties
where techniques based on the Langlands correspondence do not apply, although one can at least use them to establish a finiteness result
\cite[Theorem~1.1]{esnault-kerz} (see also \cite[Corollary~3.7.6]{kedlaya-companions}).

Another application of crystalline companions is to improved ``cut-by-curves'' criteria for
detecting whether a convergent $F$-isocrystal is overconvergent, building on the work of Shiho
\cite{shiho-cut-over} (see also \cite[Theorem~5.16]{kedlaya-isocrystals}). See \cite{gku} for details.

Yet another application is to extend Drinfeld's theorem on the independence of $\ell$ of the $\ell$-adic pro-semisimple completion of the fundamental group \cite[Theorem~1.4.1]{drinfeld-pro-simple} to include a corresponding statement for $\ell=p$, modeled on the case $\dim(X) = 1$ covered by \cite[Theorem~5.1.1]{drinfeld-pro-simple}. This can also be done in the context of ``Drinfeld's lemma'', where one considers products of varieties equipped with partial Frobenius maps; see  \cite[Proposition~4.2.4]{kedlaya-xu}.

We expect many additional applications to arise in due course, some of which do not explicitly refer to any
$p$-adic behavior. For example, the existence of companions strengthens the work
of Krishnamoorthy--P\'al on the existence of abelian varieties associated to $\ell$-adic representations
\cite{krishnamoorthy-pal, krishnamoorthy-pal2}. It is an intriguing open question whether the existence of companions
can be used to make even further progress on the existence of motives associated to \'etale or crystalline coefficients.

We note in passing that in place of coefficient objects of rank $n$, which are implicitly $\GL_n$-torsors, one may consider similar objects with $\GL_n$ replaced by a more general reductive group over $\QQ$. We will not treat the resulting conjecture in any great detail, but see \S\ref{sec:reduction} for a limited discussion.

\section{Background on algebraic stacks}

We begin with some background material, mostly related to algebraic stacks, that will play a crucial role in our construction of crystalline companions. We follow the conventions of the Stacks Project \stacktag{026M}; a more informal overview can be found in \stacktag{072I}.

 Before proceeding, we fix some geometric conventions that run throughout the paper.

\begin{notation} \label{running notation}
Throughout this paper, let $k$ be a perfect field of characteristic $p$;
our main results require $k$ to be finite, but we will impose this hypothesis explicitly when needed.
Let $K$ denote the fraction field of the ring $W(k)$ of $p$-typical Witt vectors with coefficients in $k$.
Let $X$ denote a smooth (but not necessarily geometrically irreducible) separated scheme of finite type over $k$. 
Let $X^\circ$ denote the set of closed points of $X$ and let $|X|$ denote the entire underlying topological space (compare Definition~\ref{D:schematic points}).
We interpret the (profinite) \'etale fundamental group $\pi_1(X)$ as a groupoid rather than a group, and thus do not specify a basepoint; note that its abelianization $\pi_1^{\ab}(X)$ is unambiguously a group.
\end{notation}

\begin{defn}
By a \emph{curve} over $k$, we will always mean a scheme which is smooth of dimension 1 and geometrically irreducible over $k$, but not necessarily proper over $k$ (this condition will be specified separately as needed).
A \emph{curve in $X$} is a locally closed subscheme of $X$ which is a curve over $k$ in the above sense.
\end{defn}

\begin{defn}
A \emph{smooth pair} over a base scheme $S$ is a pair $(Y,Z)$ in which $Y$ is a smooth $S$-scheme and $Z$ is a relative strict normal crossings 
divisor on $Y$; we refer to $Z$ as the \emph{boundary} of the pair. (Note that $Z = \emptyset$ is allowed.)
A \emph{good compactification} of $X$ is a smooth pair $(\overline{X}, Z)$ over $k$ with $\overline{X}$ projective (not just proper) over $k$, together with an isomorphism $X \cong \overline{X} \setminus Z$; we will generally treat the latter as an identification.
\end{defn}

\subsection{Algebraic stacks}

We start with very brief definitions, together with copious pointers to \cite{stacks-project} which are crucial for making any sense of the definitions.
As a critical link between schemes and stacks, we introduce the category of algebraic spaces as per \stacktag{025X}.
\begin{defn}
Let $\Sch$ denote the category of schemes. For $S \in \Sch$, let $\Sch_S$ denote the category of schemes over $S$,
equipped with the fppf topology.

An \emph{algebraic space} over $S$ is a sheaf $F$ on $\Sch_S$ valued in sets such that the diagonal $F \to F \times F$ is representable, and there exists a surjective \'etale morphism $h_U \to F$ for some $U \in \Sch$ (writing
$h_U$ for the functor represented by $U$). These form a category in which morphisms are natural transformations of functors, containing $\Sch$ as a full subcategory via the Yoneda embedding $U \mapsto h_U$.
\end{defn}

With this definition in hand, we can define algebraic stacks.
\begin{defn}
As per \stacktag{026N}, by an \emph{algebraic stack} over $S$, we will mean
a stack $\calX$ in groupoids over $\Sch_S$ for the fppf topology whose diagonal $\calX \to \calX \times \calX$ 
is representable by algebraic spaces,
and for which there exists a surjective smooth morphism $\Sch_U \to \calX$ for some scheme $U$.
These form a 2-category as per \stacktag{03YP}, which contains $\Sch_S$ as a full subcategory
via the operation $U \mapsto \Sch_U$.

A \emph{Deligne--Mumford (DM) stack} over $S$ is an algebraic stack $\calX$ for which the surjective smooth morphism
$\Sch_U \to \calX$ can be taken to be \'etale.

The category of algebraic stacks admits fiber products, or more precisely 2-fiber products; see \stacktag{04T2}.
\end{defn}

\begin{remark}
Any property of schemes which obeys sufficiently strong locality properties admits a natural generalization for algebraic stacks, which obeys corresponding locality properties and moreover specializes back to the original property when applied to the stacks corresponding to ordinary schemes.
For example, such a generalization exists for the properties \emph{reduced} \stacktag{04YJ} and \emph{locally noetherian} \stacktag{04YE}; there is also a construction of the
\emph{reduced closed substack} of a given stack \stacktag{0509}.

Similarly, any property of morphisms of schemes which obeys sufficiently strong locality and descent properties admits a natural generalization for algebraic stacks, which obeys corresponding locality and descent properties and moreover specializes back to the original property when applied to the stacks corresponding to ordinary schemes. For example, such generalizations exist for the properties 
\emph{quasicompact} \stacktag{050S},
\emph{quasiseparated} and
\emph{separated} \stacktag{04YV},
\emph{finite type} \stacktag{06FR}, \emph{smooth} \stacktag{075T}, \emph{universally closed} \stacktag{0511},
 and \emph{proper} \stacktag{0CL4}.
For the properties of being an \emph{open immersion}, \emph{closed immersion}, or \emph{immersion}, we make a similar construction but also require that the morphism be
representable (in schemes) \stacktag{04YK}.
\end{remark}

\begin{defn} \label{D:schematic points}
Let $\calX$ be an algebraic stack. A \emph{schematic point} of $\calX$ is a morphism of the form $\Spec L \to \calX$
where $L$ is an arbitrary field. A \emph{point} of $\calX$ is an equivalence class of schematic points under the relation that two schematic points
$\Spec L_1 \to \calX$, $\Spec L_2 \to X$ are equivalent if there exists 
a 2-commutative diagram
\[
\xymatrix{
\Spec L_3 \ar[r] \ar[d] & \Spec L_1 \ar[d] \\
\Spec L_2 \ar[r] & \calX
}
\]
with $L_3$ being a third field. The set of equivalence classes of points is denoted $\left| \calX \right|$;
this recovers the usual underlying set of a scheme. There is a natural quotient topology on $\left| \calX \right|$ which recovers the Zariski topology of a scheme
\stacktag{04Y8}.

There is no natural notion of ``closed points'' on a general algebraic stack $\calX$. Instead, it is more natural to speak of 
\emph{points of finite type} of $\calX$, meaning schematic points $\Spec L \to \calX$ for which the structure morphism is locally of finite type \stacktag{06FW};
these form a dense subset of $\left|\calX \right|$ \stacktag{06G2}.
\end{defn}

\begin{remark} \label{R:points of finite type}
Let $\calX$ be an algebraic stack of finite type over $k$.
If $\Spec L \to \calX$ is a point of finite type, then by the Nullstellensatz \stacktag{00FY} the field $L$ must be a finite extension of $k$;
in particular, if $\calX$ is a scheme of finite type over $k$, then points of finite type are the same as closed points up to a finite base extension. 

By contrast, suppose that $\calX = \Spec R$ where $R$ is a discrete valuation ring.
Then the generic point of $\calX$ is not a closed point, but it is a point of finite type because $\Frac R$ is generated as an $R$-algebra by the inverse of any single uniformizer.
\end{remark}

\begin{defn} \label{D:scheme-theoretic image}
For $f\colon \calY \to \calX$ a morphism of algebraic stacks, the \emph{scheme-theoretic image}
of $f$ is the smallest closed substack of $\calX$ through which $f$ factors. See \stacktag{0CPU}
for proof that such a substack always exists.
\end{defn}

\subsection{Moduli stacks of smooth curves}

As a reminder of the practical meaning of some of the previous definitions, we recall the basic properties of moduli stacks of smooth curves.

\begin{defn} \label{D:smooth curve fibration}
For $g,n \geq 0$, a \emph{smooth $n$-pointed genus-$g$ curve fibration} (or for short a \emph{smooth curve fibration}) consists of a morphism $f\colon Y \to S$ and $n$ morphisms $s_1,\dots,s_n\colon S \to Y$ (all in $\Sch$) satisfying the following conditions.
\begin{itemize}
\item
The morphism $f$ is smooth and proper of relative dimension 1.
\item
Each geometric fiber of $f$ has genus $g$.
\item
The morphisms $s_1,\dots,s_n$ are sections of $f$ whose images are pairwise disjoint.
\end{itemize}
We refer to the union of the images of $s_1,\dots,s_n$ in $Y$ as the \emph{pointed locus} and the complement as the \emph{unpointed locus}.
Let $M_{g,n}$ be the category over $\Sch$ whose fiber over $S \in \Sch$ consists of smooth $n$-pointed genus-$g$ curve fibrations over $S$.
\end{defn}

\begin{remark} \label{R:stable curves are schemes}
Note that the definition of the full moduli stack of curves in \stacktag{0DMJ} requires consideration of families of curves in which the total space is an algebraic space rather than a scheme; this does not change anything over the spectrum of an artinian local ring or a noetherian complete local ring \stacktag{0AE7}, but does make a difference over a more general base \stacktag{0D5D}.

However, this discrepancy does not arise for smooth curve fibrations: if $S$ is a scheme and $f\colon Y \to S$ is a smooth $n$-pointed genus-$g$ fibration in the category of algebraic spaces, then $Y$ is a scheme. That is because the hypotheses ensure that the relative
canonical bundle (for the logarithmic structure defined by $s_1,\dots,s_n$) is ample, so we may 
realize $Y$ as a closed subscheme of a particular projective bundle over $S$ (compare \stacktag{0E6F}).
\end{remark}

\begin{prop} \label{P:stable reduction}
The category $M_{g,n}$ is a smooth separated DM stack over $\ZZ$.
\end{prop}
\begin{proof}
In the case $n=0$ this is included in \stacktag{0E9C}. The general case follows by induction on $n$ because the morphism $M_{g,n+1} \to M_{g,n}$ that forgets the last marked point is representable and smooth (it is the universal curve over $M_{g,n}$ with the marked sections removed).
\end{proof}

\begin{remark} \label{R:distinguishable points}
At no point will we use the fact that the $n$ marked points in a smooth curve fibration are \emph{distinguishable}. That is,
it would be sufficient for our purposes to replace the $n$ sections $S \to Y$ with a single closed immersion $Z \to Y$ such that $Z \to Y \to S$ is finite \'etale of constant degree $n$.
\end{remark}

\subsection{Alterations and elementary fibrations}

We record two crucial results that allow us to ``blow up into controlled situations''.

\begin{defn}
An \emph{alteration} of a scheme $Y$ is a morphism $f\colon Y' \to Y$ which is proper, surjective, and generically finite \'etale. This corresponds to a \emph{separable alteration} in the sense of de Jong \cite{dejong-alterations}.
\end{defn}

\begin{prop}[de Jong] \label{P:alterations}
For $X$ smooth of finite type over $k$ (as per our running convention),
there exists an alteration $f\colon X' \to X$ such that $X'$ is smooth (but not necessarily geometrically irreducible
over $k$) and admits a good compactification.
\end{prop}
\begin{proof}
Keeping in mind that $k$ is perfect, see \cite[Theorem~4.1]{dejong-alterations}.
\end{proof}

The following is a variant of \cite[Lemma~3.1.8]{kedlaya-companions}.
\begin{cor} \label{C:stable curve fibration}
For $X$ smooth of finite type over $k$ (as per our running convention),
there exist a finite extension $k'$ of $k$,
an alteration $X' \to X \times_k k'$,
an open dense subscheme $U$ of $X'$, and a diagram
\[
\xymatrix{
U \ar[rd] \ar[r] & \overline{X}' \ar[d] \\
& S
}
\]
in which $S$ is smooth over $k'$
and $\overline{X}' \to S$ is a smooth curve fibration
with unpointed locus equal to $U$.
For any given finite subset $T$ of $X^\circ$,
we may further ensure that every point of $T(k')$ lifts to a point of $X'$ contained in $U$.
\end{cor}
\begin{proof}
At any stage in the proof, we are free to replace $X$ with either an alteration or an open dense subscheme, or to replace $k$ with a finite extension.
We may thus assume at once that $X$ is quasiprojective; we then proceed as in \cite[Proposition~3.3]{artin-sga4}
(compare \cite[Lemma~3.1.8]{kedlaya-companions}).

 Choose a very ample line bundle $\calL$ on $X$, put $n := \dim(X)$, and choose a point $z \in X^\circ$. After possibly enlarging $k$, if we make a generic choice of $n$ sections
$H_1,\dots,H_n$ of $\calL$ containing $z$, then the intersection $C$ will be zero-dimensional (but nonempty; see Remark~\ref{R:contracted section} below).
Let $\tilde{X}$ be the blowup of $X$ at $C$ and let $f\colon \tilde{X} \to \PP^{n-1}_k$ be the morphism defined by $H_1,\dots,H_n$. 

For $S$ an open dense subscheme of $\PP^{n-1}_k$, write $\tilde{X}_S := \tilde{X} \times_{\PP^{n-1}_k} S$.
For suitable $S$, the map $f$ is an \emph{elementary fibration} in the sense of \cite[Definition~3.1]{artin-sga4}; that is, it fits into a commutative diagram of the form
\[
\xymatrix{
\tilde{X}_S \ar^{j}[r] \ar_{f}[rd] & \overline{X} \ar_{\overline{f}}[d] & Z \ar_{i}[l] \ar^{g}[ld] \\
& S &
}
\]
in which:
\begin{itemize}
\item
$j$ is an open immersion with dense image in each fiber, and $i$ is a closed immersion such that
$\tilde{X}_S = \overline{X} \setminus Z$;
\item
$\overline{f}$ is smooth projective with fibers which are geometrically irreducible of dimension $1$;
\item
$g$ is finite \'etale and  surjective.
\end{itemize}
To make $\overline{f}$ into a smooth curve fibration, we must force $g$ to become a disjoint union of sections; this can be achieved by replacing $S$ with a finite \'etale cover. More precisely, take any component of $Z$ which does not map isomorphically to $S$; this component is itself a finite \'etale cover of $S$, and pulling back along it produces a fibration in which the inverse image of $Z$ splits off a component which maps isomorphically to $S$. We may repeat the construction to achieve the desired result.
\end{proof}

\begin{remark} \label{R:contracted section}
The construction in the proof of  Corollary~\ref{C:stable curve fibration} has the following useful side effect (see  Lemma~\ref{L:contracted section2}; compare also
\cite[Remark~3.1.9]{kedlaya-companions}):
the map $f$ admits a section whose image contracts to a point in $X$. In fact we can even specify this point in advance.
\end{remark}

\subsection{Moduli of coherent sheaves}
\label{subsec:moduli of coherent}

We now introduce a different moduli stack that will be more closely related to crystals.

\begin{hypothesis} \label{H:moduli of coherent}
Throughout \S\ref{subsec:moduli of coherent}, let $f\colon Y \to B$ be a separated morphism of finite presentation of algebraic spaces over some base scheme $S$.
\end{hypothesis}

\begin{defn} \label{D:stack of coherent}
Let $\Coh_{Y/B}$ denote the category in which:
\begin{itemize}
\item
the objects are triples $(T, g, \calF)$ in which $T$ is a 
scheme over $S$, $g\colon T\to B$ is a morphism over $S$,
and (writing $Y_T := Y \times_{B,g} T$)
$\calF$ is a quasicoherent $\calO_{Y_T}$-module of finite presentation
which is flat over $T$ and has support which is proper over $T$;
\item
the morphisms $(T',g',\calF') \to (T,g,\calF)$ consist of pairs $(h,\psi)$ in which
$h\colon T' \to T$ is a morphism of schemes over $B$ and 
(writing $h'\colon Y_{T'} \to Y_{T}$ for the base extension of $h$ along $f$)
$\psi\colon (h')^* \calF \to \calF'$ is an isomorphism
of $\calO_{Y_{T'}}$-modules.
\end{itemize}
These form a stack over $B$ via the functor $(T,g,\calF) \mapsto (T,g)$.
\end{defn}

\begin{prop}\label{P:Coh stack}
The category $\Coh_{Y/B}$ is an algebraic stack over $S$.
\end{prop}
\begin{proof}
See \stacktag{09DS}.
\end{proof}

\begin{prop} \label{P:Coh stack2}
The morphism $\Coh_{Y/B} \to B$ is quasiseparated and locally of finite presentation.
\end{prop}
\begin{proof}
See \stacktag{0DLZ}. Additional references, which impose more restrictive hypotheses but would still suffice for our purposes,
are \cite[Th\'eor\`eme~4.6.2.1]{laumon-moret-bailly} and \cite[Theorem~2.1]{lieblich}.
\end{proof}

Unless $f$ is finite, we cannot hope for $\Coh_{Y/B} \to B$ to be quasicompact. However, when $f$ is projective, we can cover $\Coh_{Y/B}$ with open substacks which are themselves quasicompact over $B$.

\begin{defn} \label{D:truncated stack}
Assume that $f$ is projective and $B$ is quasicompact,
and let $\calL$ be a line bundle on $Y$ which is very ample relative to $f$. 
For any object $(T, g, \calF) \in\Coh_{Y/B}$, we may define the associated \emph{Hilbert function}
\[
P\colon T \mapsto \QQ[t], \qquad P(x)(t) = \chi(Y_x, \mathbf{L} \iota_x^* (\calF \otimes \calL^{\otimes t}))
\]
where $\iota_x\colon x \to T$ denotes the canonical inclusion.
This function is locally constant \stacktag{0D1Z}.

For $P \in \QQ[t]$, let $\Coh_{Y/B}^{P,\calL}$
be the substack of $\Coh_{Y/B}$ consisting of those triples $(T,g,\calF)$ for which
the Hilbert function of $\calF$ (with respect to $\calL$) is identically equal to $P$.
As per \stacktag{0DNF}, $\Coh_{Y/B}^{P,\calL}$ is a closed-open substack of $\Coh_{Y/B}$ and
$\Coh_{Y/B}$ is equal to the disjoint union of the $\Coh_{Y/B}^{P,\calL}$ over all $P$ (because of the flatness condition in Definition~\ref{D:stack of coherent}).

For $m$ a positive integer, let $\Coh_{Y/B}^{P,\calL,m}$
be the locally closed substack of $\Coh_{Y/B}^{P,\calL}$ consisting of those triples $(T,g,\calF)$ for which
$f^{*} f_{*} (\calF \otimes g^*\calL^{\otimes m}) \to \calF \otimes \calL^{\otimes m}$ is surjective
and $R^i f_* (\calF \otimes g^*\calL^{\otimes m}) = 0$ for all $i>0$.
Note that $\Coh_{Y/B}^{P,\calL}$ is the union of the $\Coh_{Y/B}^{P,\calL,m}$ over all $m$.

Let $\Vect_{Y/B}$ be the substack of $\Coh_{Y/B}$ consisting of objects $(T,g,\calF)$ in which $\calF$ is a locally free $\calO_{Y_T}$-module. By the previous analysis, this is a closed-open substack of $\Coh_{Y/B}$; set $\Vect_{Y/B}^{P,\calL}
:= \Vect_{Y/B} \times_{\Coh_{Y/B}} \Coh_{Y/B}^{P,\calL}$
and $\Vect_{Y/B}^{P,\calL,m}
:= \Vect_{Y/B} \times_{\Coh_{Y/B}} \Coh_{Y/B}^{P,\calL,m}$.
\end{defn}

\begin{prop} \label{P:truncated stack}
Assume that $f$ is projective and $B$ is quasicompact.
Then for any $P,\calL,m$ as in Definition~\ref{D:truncated stack},
$\Coh_{Y/B}^{P,\calL,m}$ is quasicompact (as then is $\Vect_{Y/B}^{P,\calL,m}$).
\end{prop}
\begin{proof}
It suffices to produce a quasicompact algebraic space $W$ which surjects onto $\Coh_{Y/B}^{P,\calL,m}$. 
Let $Y \to \PP^n_B$ be the projective embedding defined by $\calL^{\otimes m}$.
Let $P_m$ be the polynomial with $P_m(t) = P(m+t)$ and put $r = P_m(0)$;
then each fiber of $\calF \otimes \calL^{\otimes m}$ is globally generated by its $r$-dimensional space of global sections. We may thus take $W$ to be the Quot space $\Quot^{P_m}_{\calO^{\oplus r}_{\PP^n_B}/\PP^n_B/B}$,
which is proper over $B$ by \stacktag{0DPA}. 
\end{proof}

\begin{remark} \label{R:Coh stack}
If $f$ is projective, then by \stacktag{0DM0} and \stacktag{0CLW}, the morphism  $\Coh^{P,\calL,m}_{Y/B} \to B$ is universally closed. However, in general it is  not separated and hence not proper.
\end{remark}

\begin{remark} \label{R:dual vb stack}
Assume that $f$ is projective and $B$ is quasicompact.
Choose $e_1,\dots,e_n \in \ZZ$ and let $\calG$ be a fixed vector bundle on $Y$.
It follows from Proposition~\ref{P:truncated stack}
that  the functor from the $n$-fold fiber product of $\Vect_{Y/B}^{P,\calL,m}$ over $B$
to $\Vect_{Y/B}$ given by
\[
((T,g,\calF_1), \dots, (T,g,\calF_n)) \mapsto (T,g,\calF_1^{\otimes e_1} \otimes \cdots \otimes \calF_n^{\otimes e_n} \otimes g^* \calG)
\]
factors through 
$\bigcup_{P',m} \Vect_{Y/B}^{P',\calL,m'}$ for some \emph{finite} union over pairs $(P',m')$.
\end{remark}

We can also consider moduli stacks of morphisms between coherent sheaves.
\begin{defn}
Let $\CohHom_{Y/B}$ denote the category of tuples $(T,g,\calF_1, \calF_2, h)$ in
which $(T,g,\calF_1),(T,g,\calF_2) \in \Coh_{Y/B}$ and $h \in \Hom_{\calO_{Y_T}}(\calF_1, \calF_2)$.
\end{defn}

\begin{prop} \label{P:vect hom}
Assume that $f$ is projective and $B$ is noetherian. Then 
the two morphisms $\CohHom_{Y/B} \to \Coh_{Y/B}$ taking
$(T,g,\calF_1, \calF_2, h)$ to $(T,g,\calF_1)$ and $(T,g,\calF_2)$
induce a morphism $\CohHom_{Y/B} \to \Coh_{Y/B} \times_B \Coh_{Y/B}$
which is representable (in algebraic spaces), affine, and of finite presentation.
\end{prop}
\begin{proof}
This reduces to \stacktag{08JX} as in the proof of \stacktag{08JY} (whose condition (4) holds because $f$ is proper).
\end{proof}

\begin{cor} \label{C:vect hom}
The category $\CohHom_{Y/B}$ is an algebraic stack over $S$.
The morphism $\CohHom_{Y/B} \to B$ is quasiseparated and locally of finite presentation.
\end{cor}
\begin{proof}
This follows by combining Proposition~\ref{P:Coh stack}, Proposition~\ref{P:Coh stack2},
and Proposition~\ref{P:vect hom}.
\end{proof}

\subsection{Rigidification}

As a counterpart to Remark~\ref{R:Coh stack}, we introduce the following ``rigidification'' statement.
\begin{lemma} \label{L:faithful restriction Vec}
Let $B$ be a noetherian (hence quasicompact) and connected scheme over some base scheme $S$.
Let $f \colon Y \to B$ be a smooth projective morphism of relative dimension $1$.
Let $\calL$ be a line bundle on $Y$ which is very ample relative to $f$.
Let $D$ be the zero locus of a nonzero section of $\calL^{\otimes d}$ for some positive integer $d$. Then for $d$ sufficiently large (in a sense that is invariant under base change on $B$), the following statements hold.
\begin{enumerate}
\item[(a)]
The restriction functor $\Vect^{P,\calL,m}_{Y/B} \to \Vect_{D/B}$ is faithful.
Moreover, for any two triples $(T,g,\calF_1), (T,g,\calF_2) \in \Vect^{P,\calL,m}_{Y/B}$,
whether or not a given morphism between the corresponding objects in $\Vect_{D/B}$ arises from a morphism in $\Vect^{P,\calL,m}_{Y/B}$ can be checked after replacing $T$ with an fpqc covering of $T$.
\item[(b)]
The morphism  $\Vect^{P,\calL,m}_{Y/B} \to \Vect_{D/B}$ is separated and universally closed.
\end{enumerate}
\end{lemma}
\begin{proof}
We may assume that $B$ is an affine scheme. 
Per Remark~\ref{R:dual vb stack}, we may fix a finite set of pairs $(P',m')$ for which the functor
$\Vect_{Y/B}^{P,\calL,m} \times_B \Vect_{Y/B}^{P,\calL,m} \to \Vect_{Y/B}$ 
taking a pair $(T,g,\calF_1), (T,g,\calF_2)$ to the Serre dual of $\calF_1^\dual \otimes \calF_2$
factors through $\bigcup_{P',m'} \Vect_{Y/B}^{P',\calL,m'}$. We prove both claims with $d := m + 1 + \max\{m'\}$.

Let $g \colon T \to B$ be a morphism of schemes over $S$ with $T$ affine.
Let $(T,g,\calF_1), (T,g,\calF_2) \in \Vect^{P,\calL,m}_{Y/B}$ be arbitrary for the given $T,g$.
For $\calF := \calF_1^\dual \otimes \calF_2$, we have an exact sequence
\begin{equation} \label{eq:bundle from section}
H^0(Y_T, \calF \otimes g^* \calL^{\otimes m-d}) \to H^0(Y_T, \calF \otimes g^* \calL^{\otimes m}) \to H^0(D_T, \calF \otimes g^* \calL^{\otimes m}) \to H^1(Y_T, \calF \otimes g^* \calL^{\otimes m-d}).
\end{equation}
By Serre duality plus the choice of $d$, the first term vanishes and the last term is a finite projective $\calO(B_T)$-module whose formation commutes with base extension along $T$.

From the vanishing of the first term of \eqref{eq:bundle from section} in the case $\calF = \calF_1^\dual \otimes \calF_2$, we deduce that 
$\Vect^{P,\calL,m}_{Y/B} \to \Vect_{D/B}$ is faithful. By looking at the last term of \eqref{eq:bundle from section}, we see that whether or not a particular morphism in 
$\Hom_{\Vect_{D/B}}((T,g,\calF_1), (T,g,\calF_2))$ lifts to 
$\Hom_{\Vect^{P,\calL,m}_{Y/B}}((T,g,\calF_1), (T,g,\calF_2))$  
can be checked after any dominant base extension on $T$. 
Since this includes an fpqc covering, this yields (a).

To check (b), we use the valuative criterion \stacktag{0CLT},  \stacktag{0CLW};
since $B$ is noetherian, we may restrict to discrete valuation rings \stacktag{0CM1}.
That is, we may retain the setup of (a) but with the further condition that $T$ is the spectrum of a discrete valuation ring; let $U$ be the spectrum of the fraction field.
Given an object $(U, g|_U, \calF) \in \Vect^{P,\calL,m}_{Y/B}$, the graded $\calO(U)$-module $\bigoplus_{i=0}^\infty H^0(Y_U, \calF \otimes g^* \calL^{\otimes im})$  gives rise to $\calF$ via the Proj construction. Given an extension of $\calF|_{D_U}$ across $D_T$, take the graded $\calO(T)$-submodule consisting of sections whose restrictions to $D_U$ extend across $D_T$; its Proj defines a reflexive coherent sheaf on $Y_T$ which by (a) extends $\calF$. Since $Y_T$ is regular of dimension 2, any such sheaf is locally free \stacktag{0B3N}. (This is the only step where we use that we started with a discrete valuation ring; this restriction can be eliminated using \cite[Theorem~2.7]{kedlaya-ainf}.)
\end{proof}

\begin{remark} \label{R:reflect fpqc descent data}
A typical application of Lemma~\ref{L:faithful restriction Vec} will be to ``reflect fpqc descent data'' in the following sense. 
Fix morphisms $g \colon T \to B$ and $h \colon T' \to T$ of schemes over $S$ with $h$ an fpqc covering.
Let $(T',g',\calF') \in \Vect^{P,\calL,m}_{Y/B}$ be an object whose image in $\Vect_{D/B}$
is the pullback of an object $(T,g,\calF|_D)$. 
The latter defines an descent datum on $\calF'|_D$ for the fpqc topology, which by Lemma~\ref{L:faithful restriction Vec}(a) lifts to a descent datum on $\calF$; the latter is effective by faithfully flat descent for quasicoherent sheaves \stacktag{023T}.
We conclude that there is a unique (up to unique isomorphism) object $(T,g,\calF) \in \Vect^{P,\calL,m}_{Y/B}$ which pulls back compatibly to both $\calF'$ and $\calF|_D$.
\end{remark}

\begin{cor} \label{C:faithful restriction split Azumaya}
With notation as in Lemma~\ref{L:faithful restriction Vec},
suppose that $S$ is a scheme over $\QQ$.
For $d$ sufficiently large (in a sense that is invariant under base change on $B$), the following statement holds.
Given $(T,g,\calF) \in  \Vect^{P,\calL,m}_{Y/B}$ where $\calF$ is equipped with the structure of an Azumaya algebra of degree $n$ on $Y_T$ (in the sense of \stacktag{0A2J}), together with a splitting $\calF|_{D_T} \cong \calG^\dual \otimes \calG$ for some vector bundle $\calG$ on $D_T$, there is a unique (up to unique isomorphism)
vector bundle $\calH$ on $Y_T$ equipped with isomorphisms $\calH|_{D_T} \cong \calG$,
$\calH^\dual \otimes \calH \cong \calF$ which form a commutative diagram:
\[
\xymatrix{
\calF \ar[r] \ar[d] & \calF|_{D_T} \ar[d] \\
\calH^\dual \otimes \calH \ar[r] & \calG^\dual \otimes \calG
}
\]
\end{cor}
\begin{proof}
By Lemma~\ref{L:faithful restriction Vec}, the fact that uniqueness must be up to \emph{unique} isomorphism can be checked locally for the \'etale topology (or even the fpqc topology) on $T$.
This in turn implies that uniqueness of $\calH$ can be checked locally for the \'etale topology on $T$,
which in turn implies that the existence of $\calH$ can be checked locally for the \'etale topology on $T$.

An Azumaya algebra of degree $n$ on $Y_T$ gives rise to a class in
$H^2_{\et}(Y_T, \mu_n)$. When $T$ is a geometric point, this means that the stack of local splittings is
represented by a torsor for the group scheme $\Pic(Y_T)[n]$.
For general $T$, we correspondingly get a torsor for the $n$-torsion of the relative Picard scheme of $Y_T/T$, which in particular is finite \'etale over $T$. Since we are free to pull back by such a cover thanks to the previous paragraph, we may reduce to the case where $\calF$ is split on each fiber;
this implies that $\calF$ splits locally for the \'etale topology on $T$, which by the previous
paragraph is sufficient.
\end{proof}

\begin{remark} \label{R:tannakian to azumaya}
When applying Remark~\ref{C:faithful restriction split Azumaya}, keep in mind that for a scheme $X$ over $\QQ$, specifying a $\QQ$-linear tensor functor from $\Rep_{\QQ}(\PGL_{n,\QQ})$ (the category of algebraic representations of $\PGL_{n,\QQ}$
on finite-dimensional $\QQ$-vector spaces) to the category of vector bundles on $X$
gives rise to an Azumaya algebra of degree $n$ on $X$
(and conversely, but we will not need this).
We will prefer the Tannakian point of view when introducing Azumaya algebra objects in other categories, notably in Definition~\ref{D:Azumaya algebra}.
\end{remark}

\subsection{Moduli stacks of connections}
\label{subsec:moduli conn}

The notion of a moduli stack of connections was first introduced in the setting of complex geometry by Simpson \cite{simpson1, simpson2}. In fact Simpson considered three different versions of this concept, of which the ``de Rham'' version $\calM_{\dR}$ adapts to algebraic geometry over a general base.

\begin{hypothesis}
Throughout \S\ref{subsec:moduli conn},
let $(Y,Z)$ be a smooth pair over a noetherian base scheme $B$ with $Y \to B$ proper.
Let $\calL$ be a line bundle on $Y$ which is very ample relative to $B$.
Fix $P \in \QQ[t]$ and a positive integer $m$.
\end{hypothesis}

\begin{defn} \label{D:mDR}
Let $\calM_{\dR, (Y,Z)/B}$ be the category fibered in groupoids over $\Sch_S$ in which the objects 
are tuples $(T,g,\calF,\nabla)$ where $g\colon T \to B$ is a morphism of schemes over $S$;
$\calF$ is a quasicoherent, finitely presented, $T$-flat $\calO_{Y_T}$-module; and
$\nabla \colon \calF \to \calO \otimes \Omega^1_{(Y_T,Z_T)/T}$ is an integrable logarithmic
(with respect to the canonical logarithmic structure on $Y_T$ defined by $Z_T$) connection.
By forgetting about the connection, we get a natural morphism $\calM_{\dR, (Y,Z)/B} \to \Coh_{Y/B}$.
\end{defn}

\begin{prop} \label{P:conn1}
The category $\calM_{\dR,(Y,Z)/B}$ is an algebraic stack locally of finite presentation over $B$ with affine diagonal.
\end{prop}
\begin{proof}
See \cite[Proposition~5.3]{oh-shimizu}.
\end{proof}

As the existence of an integrable logarithmic connection does not enforce local freeness on the underlying sheaf (except when $B$ is a $\QQ$-scheme and $Z = \emptyset$), we also state the following.
\begin{lemma} \label{L:conn2}
The category $\calM_{\dR,(Y,Z)/B} \times_{\Coh_{Y/B}} \Vect_{Y/B}$ is an open substack of 
$\calM_{\dR,(Y,Z)/B}$.
\end{lemma}
\begin{proof}
See \cite[Proposition~5.5]{oh-shimizu}.
\end{proof}

We now add a condition on nilpotence of residues, keeping in mind that there is some subtlety when doing this in the presence of nilpotents on the base (see Remark~\ref{R:nilpotent base}).

\begin{defn}
Let $\Conn_{(Y,Z)/B}$ be the closed substack of $\calM_{\dR,(Y,Z)/B} \times_{\Coh_{Y/B}} \Vect_{Y/B}$ whose residue along any component of $Z_T$ has the same characteristic polynomial as the zero endomorphism.

For $P \in \QQ[t]$ and $m$ a positive integer, set
\[
\Conn_{(Y,Z)/B}^{P,\calL,m} := \Conn_{(Y,Z)/B} \times_{\Vect_{Y/B}} \Vect^{P,\calL,m}_{Y/B}.
\]
\end{defn}

\begin{lemma} \label{L:quasicompact connection stack}
The category $\Conn_{(Y,Z)/B}$ is an algebraic stack locally of finite presentation over $B$ with affine diagonal. The stack $\Conn^{P,\calL,m}_{(Y,Z)/B}$
is qcqs and of finite type over $B$.
\end{lemma}
\begin{proof}
The first assertion follows from Proposition~\ref{P:conn1}, Lemma~\ref{L:conn2}, and the fact that 
$\Conn_{(Y,Z)/B}$ is cut out of $ \calM_{\dR,(Y,Z)/B} \times_{\Coh_{Y/B}} \Vect_{Y/B}$  by a closed condition. Since $B$ is
noetherian, we may replace ``finite presentation'' with ``finite type'' in that assertion;
by this plus Proposition~\ref{P:truncated stack}, the stack 
$\Conn^{P,\calL,m}_{(Y,Z)/B}$ is quasicompact.
\end{proof}

We next consider morphisms of connections.

\begin{defn}
Let $\ConnHom_{(Y,Z)/B}$ denote the category of tuples $(T,g,\calF_1,\nabla_1,\calF_2,\nabla_2,h)$
in which $(T,g,\calF_1,\nabla_1), (T,g,\calF_2,\nabla_2) \in \Conn_{(Y,Z)/B}$ and $h \in \Hom_{\calO_{Y_T}}(\calF_1,\calF_2)$ is a \emph{horizontal} morphism (i.e., it commutes with the connections).
\end{defn}

\begin{prop} \label{P:conn hom}
The two morphisms $\ConnHom_{(Y,Z)/B} \to \Conn_{(Y,Z)/B}$ taking
$(T,g,\calF_1,\nabla_1,\calF_2,\nabla_2,h)$ to $(T,g,\calF_1,\nabla_1)$ and $(T,g,\calF_2,\nabla_2)$
induce a morphism $\ConnHom_{(Y,Z)/B} \to \Conn_{(Y,Z)/B} \times_B \Conn_{(Y,Z)/B}$
which is representable (in algebraic spaces), affine, and of finite presentation.
\end{prop}
\begin{proof}
This follows from Proposition~\ref{P:vect hom} and the fact that horizontality is a closed condition on $h$.
\end{proof}

\begin{cor} \label{C:conn hom}
The category $\ConnHom_{(Y,Z)/B}$ is an algebraic stack over $S$.
The morphism $\ConnHom_{(Y,Z)/B} \to B$ is locally of finite presentation with affine diagonal.
\end{cor}
\begin{proof}
This follows by combining Proposition~\ref{P:conn1} with Proposition~\ref{P:conn hom}.
\end{proof}

\begin{remark}
In the original work of Simpson, one gets some finiteness properties for moduli stacks of connections
from the fact that the existence of an integrable logarithmic connection on a vector bundle implies a bound on the first Chern class of the bundle (the classical case being Atiyah's result that when there are no singularities at all, the first Chern class of the bundle vanishes). The failure of this point in positive characteristic will force us to detour through some additional discussion; see  \S\ref{sec:uniformities}.
\end{remark}

In positive characteristic, we have the following construction introduced by Katz \cite{katz}.
\begin{defn} \label{D:conn}
Let $p$ be a prime and suppose that $B$ is an $\FF_p$-scheme.
For any object $(T,g,\calF,\nabla) \in \Conn_{(Y,Z)/B}$ and any $\calO_T$-linear derivation
$D$ on $\calF$,
$D^p$ is again an $\calO_T$-linear derivation on $\calF$; this applies in particular when $\calF = \calO_T$. We thus obtain a $\calO_{Y_T}$-linear map 
\[
\psi\colon \calF \to \calF \otimes_{\calO_{Y_T}} \Omega^1_{(Y_T,Z_T)/B}
\]
whose contraction with a derivation $D$ is
\[
\langle \psi, D\rangle = \langle \nabla, D \rangle^p - \langle \nabla, D^p \rangle.
\]
We call $\psi$ the \emph{$p$-curvature} of $\calF$.

Let $\Conn^p_{(Y,Z)/B}$ be the closed substack of $\Conn_{(Y,Z)/B}$ consisting of objects $(T,g,\calF,\nabla)$
for which the $p$-curvature of $\calF$ vanishes. Set
\[
\Conn_{(Y,Z)/B}^{p,P,\calL,m} := \Conn^p_{(Y,Z)/B} \times_{\Coh_{Y/B}} \Vect^{P,\calL,m}_{Y/B}.
\]
\end{defn}

\begin{lemma} \label{L:quasicompact connection stack2}
The stack $\Conn_{(Y,Z)/B}^{p,P,\calL,m}$ is qcqs and of finite type over $B$.
\end{lemma}
\begin{proof}
This is immediate from Lemma~\ref{L:quasicompact connection stack}.
\end{proof}

\begin{remark} \label{R:nilpotent base}
While for our purposes we can get by considering only the distinction between zero and nonzero $p$-curvature, in practice it is more natural to distinguish between \emph{nilpotent} $p$-curvature and the contrary. However, from the stack-theoretic point of view this is a bit subtle because one cannot uniformly bound the index of nilpotency except over a reduced base. See \cite[\S 5]{oh-shimizu} for more discussion.
\end{remark}

\section{Coefficient objects}
\label{sec:companions}

We review some relevant properties of coefficient objects and companions.
Since we will be using terminology and notation from both \cite{kedlaya-isocrystals} and \cite{kedlaya-companions},
often with little comment, we recommend keeping those sources handy while reading.
Note also that the order of presentation here does not match the logical order in which concepts must be developed at the foundational level (in either the \'etale or crystalline settings).

\begin{hypothesis}
Throughout \S\ref{sec:companions}, assume that $k$ is finite.
\end{hypothesis}

\subsection{Coefficient objects and algebraicity}

\begin{defn}
By a \emph{coefficient object} on $X$, we will mean an object of one of the categories $\Weil(X) \otimes \overline{\QQ}_\ell$, the category of lisse Weil $\overline{\QQ}_\ell$-sheaves on $X$ for some prime $\ell \neq p$ (an \emph{\'etale coefficient});
or $\FIsoc^\dagger(X) \otimes \overline{\QQ}_p$, the category of overconvergent $F$-isocrystals on $X$ with coefficients in $\overline{\QQ}_p$, in the sense of \cite[Definition~9.2]{kedlaya-isocrystals} (a \emph{crystalline coefficient}).
By the \emph{full coefficient field} of a coefficient object (or its ambient category), we mean the algebraic closure of the field $\QQ_\ell$ in the former case and $\QQ_p$ in the latter case (without completion).

A coefficient object on $X$ is \emph{absolutely irreducible} if for any finite extension $k'$ of $k$, the pullback of $\calE$ to $X_{k'}$ is irreducible.
\end{defn}

\begin{lemma} \label{L:restriction is fully faithful}
Let $U$ be an open dense subscheme of $X$. For any category of coefficient objects, restriction from coefficient objects over $X$ to coefficient objects over $U$ is fully faithful.
\end{lemma}
\begin{proof}
See \cite[Lemma~1.2.2]{kedlaya-companions}.
\end{proof}

\begin{defn}
We say that a coefficient object on $X$ is \emph{constant} if it arises by pullback from $\Spec(k)$.
We refer to the operation of tensoring with a constant coefficient object of rank $1$, or the result of this operation, as a \emph{constant twist}.
\end{defn}

\begin{lemma} \label{L:decompose to finite order}
Suppose that $X$ is irreducible.
For every coefficient object $\calE$ on $X$,
there exists a constant twist of $\calE$ whose determinant is of finite order. 
\end{lemma}
\begin{proof}
See \cite[Lemma~1.1.3]{kedlaya-companions}.
\end{proof}

\begin{defn} \label{D:algebraic}
For $\calE$ a coefficient object on $X$ and $x \in X^\circ$, let $F_x$ denote the linearized Frobenius action on $\calE_x$, and let $P(\calE_x, T)$ denote the reverse characteristic polynomial (Fredholm determinant) of $F_x$ in the variable $T$. We say that $\calE$ is \emph{algebraic} if $P(\calE_x, T) \in \overline{\QQ}[T]$ for all $x \in X^\circ$; in this case, we have $P(\calE_x,T) \in E[T]$ for some number field $E$ independent of $x$
\cite[Theorem~3.4.2]{kedlaya-companions}. To identify such a number field, we may say that $\calE$ is \emph{$E$-algebraic}.
(Beware that even if $\calE$ is algebraic, the \emph{roots} of the polynomials $P(\calE_x, T)$ need not belong to a single number field independent of $x$.)

For $\calE$ algebraic, the \emph{$L$-function} of $\calE$ is the power series
\[
L(\calE, T) = \prod_{x \in X^\circ} P(\calE_x, T^{[\kappa(x):k]})^{-1} \in \overline{\QQ} \llbracket T \rrbracket,
\]
which represents a rational function by virtue of the Lefschetz trace formula (see \cite[Theorem~9.6]{kedlaya-isocrystals} in the crystalline case).
\end{defn}

\begin{lemma} \label{L:irreducible to algebraic}
Let $\calE$ be a coefficient object on $X$. If $\calE$ is irreducible and $\det(\calE)$ is of finite order, then
$\calE$ is algebraic.
\end{lemma}
\begin{proof}
See \cite[Corollary~3.4.3]{kedlaya-companions}.
\end{proof}

\begin{cor} \label{C:on curve constant twist}
Suppose that $X$ is irreducible.
Let $\calE$ be an irreducible coefficient object on $X$. Then some constant twist of $\calE$ is algebraic.
\end{cor}
\begin{proof}
Combine Lemma~\ref{L:decompose to finite order} and Lemma~\ref{L:irreducible to algebraic}.
\end{proof}

\begin{remark} \label{R:independence of base field}
Note that in a certain sense, the categories of coefficient objects on $X$ do not depend on the base field $k$;
see \cite[Definition~9.2]{kedlaya-isocrystals} for discussion of this point in the crystalline case. For this reason, it will often be harmless for us to assume that $X$ is absolutely irreducible.
\end{remark}

\begin{lemma} \label{L:bound coefficient extension}
Fix a category $\calC$ of coefficient objects and an embedding of $\overline{\QQ}$ into the full coefficient field $\overline{\QQ}_*$ of $\calC$.
Let $E$ be a number field within $\overline{\QQ}$ and let $L_0$ be the completion of $E$ in $\overline{\QQ}_*$. 
Let $\calE$ be an $E$-algebraic coefficient object on $X$ in $\calC$ of rank $r$.
Then there exists a finite extension $L$ of $L_0$, depending only on $L_0$ and $r$
(not on $X$ or $\calE$ or $E$), for which $\calE$ can be realized as an object of
$\Weil(X) \otimes L$ (in the \'etale case) or $\FIsoc^\dagger(X) \otimes L$ (in the crystalline case).
\end{lemma}
\begin{proof}
See \cite[Corollary~3.3.5]{kedlaya-companions}.
\end{proof}

\subsection{Monodromy groups and Tannakian categories}

\begin{defn}
By a \emph{geometric coefficient} object on $X$, we will mean one of the categories $\Weil(X_{\overline{k}}) \otimes \overline{\QQ}_\ell$ for $\ell \neq p$ (an \emph{\'etale geometric coefficient});
or $\Isoc^\dagger(X) \otimes \overline{\QQ}_p$, the category of overconvergent isocrystals on $X$ without Frobenius structure with coefficients in $\overline{\QQ}_p$ (a \emph{crystalline geometric coefficient}).

We will only ever consider geometric coefficient objects that appear in the Tannakian category generated by some coefficient object. Any such object which is irreducible can be promoted to a coefficient object after a finite base extension on $k$ \cite[Remark~1.3.9]{kedlaya-companions}. Consequently, a coefficient object $\calE$ is absolutely irreducible if and only if it is irreducible as a geometric coefficient object.
\end{defn}

\begin{defn} \label{D:monodromy}
Let $\calE$ be a coefficient object on $X$.
We write $G(\calE)$ and $\overline{G}(\calE)$ for the arithmetic and geometric monodromy groups of $\calE$, respectively,
defined in the sense of Crew in the crystalline case
\cite[Definition~1.3.4]{kedlaya-companions}.
In both cases these are Tannakian automorphism groups over the full coefficient field of $\calE$. In particular, $\calE$ is semisimple if and only if $G(\calE)^\circ$ is a semisimple algebraic group.
\end{defn}

\begin{lemma} \label{L:monodromy}
Let $\calE$ be a coefficient object on $X$.
\begin{enumerate}
\item[(a)] If $\calE$ is semisimple, then $\overline{G}(\calE)^\circ$ is a semisimple algebraic group.
\item[(b)] There exists a finite \'etale cover $f\colon X' \to X$ such that $\overline{G}(f^* \calE) = \overline{G}(\calE)^\circ$.
\end{enumerate}
\end{lemma}
\begin{proof}
For (a), see \cite[Proposition~1.3.12]{kedlaya-companions}.
For (b), see \cite[Proposition~1.3.11]{kedlaya-companions}
\end{proof}

\begin{remark} \label{R:same Tannakian category}
The previous discussion has the following key consequence. Let $\calE$ be a geometrically irreducible coefficient object of rank $r$. By Lemma~\ref{L:monodromy}, $\overline{G}(\calE)^\circ$ is a simple Lie group and so its associated Lie algebra is also simple.
For $i=1,\dots,r-1$, the natural surjective map $\overline{G}(\calE) \to \overline{G}(\wedge^i \calE)$ must therefore induce an isomorphism of Lie algebras, so its kernel and cokernel are finite. Moreover, if $\overline{G}(\calE)$ is connected, then $\overline{G}(\calE) \to \overline{G}(\wedge^i \calE)$ is surjective and its kernel is contained in the center of $\overline{G}(\calE)$.

By the same token, $\overline{G}(\calE) \to \overline{G}(\calE^\dual \otimes \calE)$ has finite kernel and cokernel.
Moreover, if $\overline{G}(\calE)$ is connected, then $\overline{G}(\calE) \to \overline{G}(\calE^\dual \otimes \calE)$ is surjective and its kernel is \emph{equal} to the center of $\overline{G}(\calE)$.
\end{remark}

We only need the following in the \'etale case, but it will also hold in the crystalline case; see Theorem~\ref{T:preservation of monodromy group}.
\begin{lemma} \label{L:pin monodromy restriction}
Let $\calE$ be an \'etale coefficient object on $X$.
Then there exists a finite subset $H$ of $X^\circ$ such that
for any connected locally closed subscheme $D$ of $X$ which is smooth over $k$, the inclusions  
\[
G(\calE|_D) \to G(\calE), \qquad \overline{G}(\calE|_D) \to \overline{G}(\calE)
\]
are isomorphisms.
\end{lemma}
\begin{proof}
In the \'etale case, $\calE$ can be interpreted as a continuous representation $\pi_1(X) \to \GL_n(\overline{\QQ}_\ell)$ and
$G(\calE)$ (resp. $\overline{G}(\calE)$) can be interpreted as the Zariski closure of the image of $\pi_1(X)$ (resp. $\pi_1(X_{\overline{k}})$) in $\GL_n(\overline{\QQ}_\ell)$. The claim is then a consequence of the noetherian property of the scheme $\GL_n$.
\end{proof}

\subsection{The companion relation for coefficients}

We next introduce the companion relation, together with a number of reduction steps for the construction of crystalline companions.

\begin{defn} \label{D:companions}
Fix coefficient objects $\calE$ and $\calF$ on $X$,
as well as an isomorphism $\iota$ between the algebraic closures of $\QQ$ in the full coefficient fields of $\calE$ and $\calF$. The map $\iota$ will often not be mentioned explicitly; in other case, it will be specified implicitly in terms of a place of the algebraic closure of $\QQ$ in one of the full coefficient fields.

We say that $\calE$ and $\calF$ are \emph{companions}
(with respect to $\iota$) if for each $x \in X^\circ$, the coefficients of $P(\calE_x, T)$ and $P(\calF_x, T)$ are identified via $\iota$; in particular, this can only occur if both $\calE$ and $\calF$ are algebraic.
Given $\calE$, the companion relation determines $\calF$ up to semisimplification 
(see Lemma~\ref{L:companion irreducible} below).
In case $\calF \in \FIsoc^\dagger(X) \otimes \overline{\QQ}_p$, we also say that $\calF$ is a \emph{crystalline companion} of $\calE$.
\end{defn}

\begin{lemma} \label{L:companion irreducible}
Let $\calE$ and $\calF$ be coefficient objects on $X$ which are companions.
\begin{enumerate}
\item[(a)]
If $\calE$ is irreducible, then so is $\calF$.
\item[(b)]
If $\calE$ is absolutely irreducible, then so is $\calF$.
\item[(c)]
If $\det(\calE)$ is of finite order, then so is $\det(\calF)$.
\item[(d)]
If $\calF'$ is another coefficient object in the same category as $\calF$ which is also a companion of $\calE$, then $\calF$ and $\calF'$ have isomorphic semisimplifications.
In particular, if one of $\calF$ and $\calF'$ is absolutely irreducible, then $\calF$ and $\calF'$ are isomorphic.
\end{enumerate}
\end{lemma}
\begin{proof}
Part (a) is \cite[Theorem~3.3.1(a)]{kedlaya-companions}. Part (b) is a formal consequence of (a).
Part (c) is \cite[Corollary~3.3.4]{kedlaya-companions}.
Part (d) is \cite[Theorem~3.3.2(b)]{kedlaya-companions} for the first assertion, plus (b) for the second assertion.
\end{proof}

\begin{lemma} \label{L:extension of companions}
Let $U$ be an open dense subscheme of $X$. Let $\calE$ be a semisimple coefficient object on $U$. Let $\calF$ be a coefficient object on $X$ whose restriction to $U$ is a companion of $\calE$. Then $\calE$ extends to a coefficient object on $X$, and any such extension is a companion of $\calF$.
\end{lemma}
\begin{proof}
By \cite[Corollary~3.3.3]{kedlaya-companions}, there exists an extension of $\calE$ which is a companion of $\calF$. By Lemma~\ref{L:restriction is fully faithful}, this is the unique extension of $\calE$ to $X$.
\end{proof}

\begin{remark} \label{R:radicial}
Let $f\colon X' \to X$ be a radicial morphism. Then for any fixed category of coefficient objects, pullback via $f$
defines an equivalence of categories between the coefficients over $X$ and the coefficients over $X'$.
\end{remark}

\begin{lemma} \label{L:alteration reduction}
Let $\calE$ be a coefficient object on $X$. 
Let $f\colon X' \to X$ be a dominant morphism (this includes the case of a base extension on $k$).
If $f^* \calE$ admits a crystalline companion, then so does $\calE$.
\end{lemma}
\begin{proof}
By restriction to a rational multisection of $f$, we may put ourselves in the position where $f$ is the composition
of an open immersion with dense image and an alteration. We may handle the two cases separately;
the former is treated by Lemma~\ref{L:extension of companions},
while for the latter see \cite[Corollary~3.6.3]{kedlaya-companions}.
\end{proof}

\begin{remark} \label{R:same Tannakian category2}
In the notation of Remark~\ref{R:same Tannakian category},
if $\overline{G}(\calE)$ is connected and there exist some $i \in \{1,\dots,r-1\}$ and some nonconstant irreducible subobject of $\wedge^i \calE$ that admits a crystalline companion, then $\calE^\dual \otimes \calE$ also admits a crystalline companion.
\end{remark}

\subsection{Tame and docile coefficients}

Using the fact that the existence of companions can be checked after an alteration 
(Lemma~\ref{L:alteration reduction}), we will be able to limit our attention to $p$-adic coefficient objects of a relatively simple sort. We repeat here \cite[Definition~1.4.1]{kedlaya-companions}.
\begin{defn} \label{D:tame docile}
Let $X \hookrightarrow \overline{X}$ be an open immersion with dense image. 
Let $D$ be an irreducible divisor of $\overline{X}$ with generic point $\eta$.
\begin{itemize}
\item
For $\ell$ a prime not equal to $p$, 
an object $\calE$ of $\Weil(X) \otimes \overline{\QQ}_\ell$ is \emph{tame} (resp.\ \emph{docile}) along $D$ if 
the action of the inertia group at $\eta$ on $\calE$ is tamely ramified
(resp. tamely ramified and unipotent). 
\item
An object $\calE$ of $\FIsoc^\dagger(X) \otimes \overline{\QQ}_p$ is \emph{tame} 
(resp.\ \emph{docile}) along $D$ if $\calE$ has
$\QQ$-unipotent monodromy in the sense of \cite[Definition~1.3]{shiho-log}
(resp. unipotent monodromy in the sense of \cite[Definition~4.4.2]{kedlaya-semi1}) along $D$.
\end{itemize}
\end{defn}

\begin{lemma} \label{L:alter to docile}
For any prime $\ell \neq p$ and any object $\calE$ of $\Weil(X) \otimes \overline{\QQ}_\ell$, there exists an alteration
$f\colon X' \to X$ such that $X'$ admits a good compactification with respect to which $f^* \calE$ is docile.
\end{lemma}
\begin{proof}
See \cite[Proposition~1.4.6]{kedlaya-companions}. This also covers the case $\ell = p$, but we will not need that here.
\end{proof}

\subsection{Companions on curves}

We next recall the critical results on the existence of companions on curves, as extracted from the work of L. Lafforgue \cite{lafforgue} and Abe \cite{abe-companion} on the global Langlands correspondence in positive characteristic.

\begin{theorem}[L. Lafforgue, Abe] \label{T:companion dimension 1}
For $X$ of dimension $1$, every coefficient object on $X$ which is irreducible with determinant of finite order is uniformly algebraic and admits companions in all categories of coefficient objects, which are again irreducible with determinant of finite order.
\end{theorem}
\begin{proof}
See \cite[Theorem~2.2.1]{kedlaya-companions} and references therein.
\end{proof}

Although we frame the following statement as a corollary of Theorem~\ref{T:companion dimension 1},
it also incorporates significant intermediate results of V. Lafforgue, Deligne, Drinfeld, Abe--Esnault, and Kedlaya;
see \cite{kedlaya-companions} for more detailed attributions.
\begin{cor} \label{C:deligne prior}
Parts (i)--(v) of Theorem~\ref{T:deligne} hold in general. Under the additional hypothesis $\dim(X) = 1$,
part (vi) of Theorem~\ref{T:deligne} also holds.
\end{cor}
\begin{proof}
The first assertion is \cite[Theorem~0.4.1]{kedlaya-companions}. The second assertion is 
\cite[Theorem~0.2.1]{kedlaya-companions}.
\end{proof}

We record some refinements of this statement which we will also need.
\begin{cor} \label{C:algebraic companion}
For $X$ of dimension $1$, any algebraic coefficient object admits
companions in all categories of coefficient objects.
\end{cor}
\begin{proof}
See \cite[Theorem~2.2.1, Corollary~2.2.3]{kedlaya-companions}.
\end{proof}

\begin{cor} \label{C:companions are tame}
For any tame (resp.\ docile) coefficient object on $X$, its companions are also tame (resp.\ docile).
\end{cor}
\begin{proof}
For the proof when $\dim(X) = 1$ (which is the only case we will use here), see \cite[Corollary~2.4.3]{kedlaya-companions}. The general case follows from this case using \cite[Lemma~1.4.9]{kedlaya-companions}.
\end{proof}

\subsection{Newton polygons}

The following argument corrects an error in the proof of \cite[Theorem~1.3.3(i,ii)]{drinfeld-kedlaya}: the proof of \cite[Lemma~5.3.1]{drinfeld-kedlaya} is incorrect, but is supplanted by part (e) of the following statement. We state these results because we need them as part of the construction of crystalline companions, but \emph{a posteriori} we will get much more precise results (see \S\ref{subsec:Newton}).

\begin{lemma-def}\label{L:generic Newton polygon}
Fix a category $\calC$ of coefficient objects and a normalized $p$-adic valuation $v_p$ of the algebraic closure of $\QQ$ in the full coefficient field of $\calC$.
For $\calE$ an algebraic object of $\calC$,
there exists a unique function $x \mapsto N_x(\calE)$ from $X$ to the set of polygons in $\RR^2$ with the following properties.
\begin{enumerate}
\item[(a)]
For $x \in X^\circ$, $N_x(\calE)$ equals the lower convex hull of the set of points
\[
\left\{ \left( i, \frac{1}{[\kappa(x)\colon\! k]} v_p(a_{i}) \right)\colon i = 0,\dots,d \right\} \subset \mathbb{R}^2
\]
where $P(\calE_x, T) = \sum_{i=0}^d a_i T^i \in \overline{\QQ}[T]$ with $a_0 = 1$. 
\item[(b)]
For $x \in X$ with Zariski closure $Z$, for all $z \in Z$, $N_z(\calE)$ lies on or above $N_x(\calE)$ with the same right endpoint. In particular, the right endpoint is constant on each connected component of $X$.
\item[(c)]
In (b), the set $U$ of $z \in Z^\circ$ for which equality holds is nonempty.
\item[(d)]
In (c), for each curve $C$ in $Z$, the inverse image of $U$ in $C$ is either empty or the complement of a finite set of closed points.
\item[(e)]
For $x \in X$, the vertices of $N_x(\calE)$ all belong to $\ZZ \times \frac{1}{N} \ZZ$ for some positive integer $N$ (which may depend on $X$ and $\calE$, but not on $x$).
\end{enumerate}
\end{lemma-def}
\begin{proof}
Using Lemma~\ref{L:decompose to finite order}, we may reduce to the case where $\calE$ is irreducible with determinant of finite order.
By Corollary~\ref{C:deligne prior}, we may apply part (iv) of Theorem~\ref{T:deligne} to obtain uniform upper and lower bounds on $N_x(\calE)$ for all $x \in X^\circ$.

For each curve $C$ in $X$, let $\calF_C$ denote a crystalline companion of $\calE|_C$ given by
Corollary~\ref{C:algebraic companion}; it has coefficients in some finite extension $L$ of $\QQ_p$,
which by Lemma~\ref{L:bound coefficient extension} can be chosen independently of $C$.
By comparison with the usual definition of Newton polygons for crystalline coefficients (see Definition~\ref{D:partially overconvergent isocrystals} and Lemma~\ref{L:convergent Newton polygon}),
we deduce that all of the claims hold with $X$ replaced by $C$, taking $N = [L:\QQ_p]$.
This implies (e) for $x \in X^\circ$ by choosing $C$ to pass through $x$.

Now let $x \in X$ be arbitrary and let $Z$ be the Zariski closure of $x$.
As $y$ varies over $Z^\circ$, the vertices of $N_y(\calE)$ are bounded above and below (by the first paragraph) and belong to $\ZZ \times \frac{1}{N} \ZZ$ (by the second paragraph),
and so $N_y(\calE)$ can assume only finitely many distinct values; we can thus find $y \in Z^\circ$ for which $N_y(\calE)$ is minimal
(note that \emph{a priori} there may be multiple minima). 
We then choose some such $y$ and set $N_x(\calE) := N_y(\calE)$;
this definition clearly satisfies (a), (c), and (e).

To check (d), we compare it to the statement of (d) for $\calF_C$ which is already known (see above); the only possible discrepancy is the generic point $\eta$ of $C$.
If $\eta \in U$, then by (c) and (d) for $\calF_C$ we deduce that $U$ is the complement of a finite set of closed points.
Otherwise, by (b) for $\calF_C$ we deduce that $U$ is empty.

To check (b), let $Z'$ be the Zariski closure of $z$. By (c), there exists $z' \in Z^{\prime \circ}$ for which $N_{z'}(\calE) = N_z(\calE)$. Choose a curve $C$ in $Z$ containing $z'$
and the point $y$ from the third paragraph. By (b) for $\calF_C$, we deduce the claim.

To establish uniqueness, we induct on $\dim(X)$. We may assume that $X$ is irreducible with generic point $\eta$.
For $\dim(X) = 1$, the function $N_x(\calE)$ is determined by (a) for $x \in X^\circ$ and by (c) and (d) for $x = \eta$.
For $\dim(X) > 1$, the induction hypothesis determines $N_x(\calE)$ for $x \neq \eta$. Meanwhile, (b) and (c) imply
that $N_\eta(\calE)$ equals a minimal value of $N_x(\calE)$ for $x \in X^\circ$, so it suffices to check that there is a unique such value. Choose $x,y \in X^\circ$ for which $N_x(\calE)$ and $N_y(\calE)$ are minimal. Choose a curve $C$ in $X$ containing both $x$ and $y$; by (d), the function $z \mapsto N_z(\calE)$
is constant on some nonempty open subset of $C$. By (b) and (c), the constant value lies below both $N_x(\calE)$ and $N_y(\calE)$, and so by minimality is equal to both of them.
\end{proof}

\begin{remark}
In the crystalline case of Lemma-Definition~\ref{L:generic Newton polygon}, if $L$ is a finite extension of $\QQ_p$ for which $\calE \in \FIsoc^\dagger(X) \otimes L$, then the $y$-coordinates of the vertices of $N_x(\calE)$ for all $x \in X^\circ$ belong to $[L:\QQ_p]^{-1} \ZZ$ (see Lemma~\ref{L:convergent Newton polygon}(e)). In the \'etale case, one can derive a similar bound on denominators using Lemma~\ref{L:bound coefficient extension}, depending on the rank of $\calE$ and the degree of the minimal number field $E$
for which $\calE$ is $E$-algebraic.
\end{remark}

\begin{remark}
By Lemma-Definition~\ref{L:generic Newton polygon}(a),
the definition of $N_x(\calE)$ remains unchanged upon replacing $\calE$ with a companion \emph{provided} that we maintain the choice of $v_p$. By contrast, if we change $v_p$ then the definition need not be preserved, even up to a constant twist; see 
Definition~\ref{D:absolutely unit-root}.
\end{remark}

\subsection{Weights}
\label{subsec:weights}

\begin{hypothesis}
Throughout \S\ref{subsec:weights}, fix a category $\calC$ of coefficient objects and an algebraic (but not topological!) embedding $\iota$ of the full coefficient field of $\calC$ into $\CC$.
\end{hypothesis}

\begin{defn}
For $\calE \in \calC$ and $x \in X^\circ$, we define the \emph{$\iota$-weights} of $\calE$ at $x$ as per \cite[Definition~3.1.2]{kedlaya-companions}. 
We say that $\calE$ is \emph{$\iota$-pure of weight $w$} if for all $x \in X^\circ$, the $\iota$-weights of $\calE$ at $x$ are all equal to $w$.
\end{defn}

\begin{lemma} \label{L:Weil II}
Suppose that $\calE \in \calC$ is $\iota$-pure of weight $w$. Then the eigenvalues of Frobenius on $H^1(X, \calE)$ all have $\iota$-absolute value at least $q^{(w+1)/2}$.
\end{lemma}
\begin{proof}
See \cite[Lemma~3.1.3]{kedlaya-companions} and onward references.
\end{proof}

\begin{lemma} \label{L:irreducible is pure}
Let $\calE$ be an irreducible coefficient object. Then $\calE$ is $\iota$-pure of some weight.
\end{lemma}
\begin{proof}
See \cite[Theorem~3.1.10]{kedlaya-companions}.
\end{proof}

\subsection{Cohomological rigidity}
\label{subsec:cohom rigidity}

For $\calE$ an absolutely irreducible coefficient object, $H^0(X, \calE^\dual \otimes \calE)$ is equal to the full coefficient field; this means that $\calE$ has no infinitesimal automorphisms. Using weights, we can also assert that $\calE$ has no infinitesimal deformations.

\begin{prop} \label{P:cohomological rigidity}
Let $\calE$ be an absolutely irreducible coefficient object on $X$. 
Let $\calF$ be the trace-zero component of $\calE^\dual \otimes \calE$. Let $\varphi$ denote the action of (geometric) Frobenius on cohomology groups. Then
\[
H^0(X, \calF)_{\varphi} = H^1(X, \calF)^{\varphi} = H^1(X, \calE^\dual \otimes \calE)^{\varphi} = 0.
\]
(Here the superscript indicates invariants while the subscript indicates coinvariants.)
\end{prop}
\begin{proof}
Since $\calE$ is absolutely irreducible,
the multiplicity of 1 as an eigenvalue of $\varphi$ on $H^0(X, \calE^\dual \otimes \calE)$ cannot exceed 1, as contributed by the trace component. Hence $H^0(X, \calF)_{\varphi} = 0$.

Let $\iota$ be an embedding of the full coefficient field of $\calE$ into $\CC$. By Lemma~\ref{L:irreducible is pure}, 
$\calE$ is $\iota$-pure of some weight $w$; then $\calE^\dual$ is $\iota$-pure of weight $-w$
and $\calE^\dual \otimes \calE$ is $\iota$-pure of weight 0. By Lemma~\ref{L:Weil II}, the eigenvalues of $\varphi$
on $H^1(X, \calE^\dual \otimes \calE)$ all have $\iota$-absolute value at least $q^{1/2}$; in particular, none of them is equal to 1. This proves the claim.
\end{proof}

\begin{remark} \label{R:cohomological rigidity}
The space $H^1(X, \calE^\dual \otimes \calE)$ computes infinitesimal deformations of $\calE$ as an object over $X_{\overline{k}}$; compare
\cite{agrawal} in the crystalline case. By the Hochschild--Serre spectral sequence, the infinitesimal deformations of $\calE$ as an object over $X$ form a space sandwiched between $H^0(X, \calF)_{\varphi}$ and $H^1(X, \calE^\dual \otimes \calE)^{\varphi}$; Proposition~\ref{P:cohomological rigidity} thus asserts the vanishing
of this space.
\end{remark}

\subsection{Azumaya algebra objects}

\begin{defn} \label{D:Azumaya algebra}
Let $\calC$ be a category of coefficient objects with full coefficient field $L$.
For $G$ a linear algebraic group over $L$, let $\Rep_L(G)$ denote the category of algebraic representations of $G$ on finite-dimensional $L$-vector spaces.
A \emph{$G$-object} in $\calC$ is an $L$-linear functor $\Rep_L(G) \to \calC$.
In particular, for $G = \GL_n$, a $G$-object is just an object of $\calC$ of rank $n$.

We define an \emph{Azumaya algebra object of degree $n$} in $\calC$ to be a $G$-object for $G = \PGL_n$.
For example, for any object $\calF \in \calC$ of rank $n$, $\calF^\dual \otimes \calF$ is an Azumaya algebra object; we say that such an object is \emph{split} and refer to $\calF$ as a \emph{splitting} thereof.
\end{defn}

\section{Convergent vs. overconvergent isocrystals}
\label{sec:np}

For the construction of crystalline companions, we will not be able to carry out the entire argument in the category of overconvergent $F$-isocrystals; instead, we make some steps in the closely related category of \emph{convergent $F$-isocrystals}. We recall the relevant points here.

\setcounter{equation}{0}
\begin{hypothesis}
Throughout \S\ref{sec:np}, let $L$ denote an algebraic extension of $\QQ_p$.
\end{hypothesis}

\subsection{Convergent isocrystals}

\begin{defn} \label{D:partially overconvergent isocrystals}
Let $\FIsoc(X) \otimes L$ denote the category of \emph{convergent 
$F$-isocrystals} with coefficients in $L$, again in the sense of \cite[Definition~9.2]{kedlaya-isocrystals}.
There is a natural restriction functor from $\FIsoc^\dagger(X) \otimes L$ 
to $\FIsoc(X) \otimes L$, which in light of Lemma~\ref{L:fully faithful to convergent} we will view as a full embedding.

We may view both $\FIsoc(X) \otimes L$ and $\FIsoc^\dagger(X) \otimes L$ as special cases of the category
$\FIsoc(X, Y) \otimes L$ of $F$-isocrystals which are overconvergent along an open immersion $X \to Y$;
namely, one gets $\FIsoc(X) \otimes L$ when $X = Y$ and $\FIsoc^\dagger(X) \otimes L$ when $Y$ is proper over $k$
\cite[Definition~2.4]{kedlaya-isocrystals}.
\end{defn}

\begin{lemma} \label{L:fully faithful to convergent}
For any open immersion $X \to Y$, 
the restriction functor $\FIsoc(X,Y) \otimes L \to \FIsoc(X) \otimes L$ is fully faithful.
In particular, $\FIsoc^\dagger(X) \otimes L \to \FIsoc(X) \otimes L$ is fully faithful.
\end{lemma}
\begin{proof}
See \cite[Theorem~5.3]{kedlaya-isocrystals}.
\end{proof}

We use the following form of Zariski--Nagata purity.
\begin{lemma} \label{L:Zariski purity1}
Let 
\[
\xymatrix{
U \ar[r] \ar[d] & W \ar[d] \\
X \ar[r] & Y
}
\]
be open immersions of smooth $k$-schemes such that:
\begin{itemize}
\item
$Y \setminus X$ is a normal crossings divisor in $X$,
\item
$X \setminus U$ has codimension at least $2$ in $X$, and
\item
$Y \setminus W$ has codimension at least $2$ in $Y$.
\end{itemize}
Then the functor 
$\FIsoc(X, Y) \otimes L \to \FIsoc(U, W) \otimes L$ is an equivalence of categories.
\end{lemma}
\begin{proof}
See \cite[Theorem~5.1]{kedlaya-isocrystals}.
\end{proof}

\begin{defn} \label{D:converget companion}
When $k$ is finite, we extend the definition of the companion relation (Definition~\ref{D:companions})
to the case where one or both objects belongs to $\calE \in \FIsoc(X) \otimes L$. Some care must be taken with this definition because Lemma~\ref{L:companion irreducible}(d) does not extend to this level of generality.
\end{defn}

\begin{defn} \label{D:log isocrystal}
For $(X,Z)$ a smooth pair,
equip $X$ with the logarithmic structure defined by $Z$.
Let $\FIsoc(X^{\log}) \otimes L$ be the category of convergent log-isocrystals on $(X,Z)$ with coefficients in $L$ equipped with a Frobenius structure,
which we insist is an isomorphism \emph{even over $Z$}. Note that this last restriction enforces that the underlying log-isocrystal has nilpotent residues along $Z$: its residues form a finite multiset of a field of characteristic $0$ stable
under multiplication by $p$.
By \cite[Theorem~6.4.5]{kedlaya-semi1},
the restriction functor $\FIsoc(X^{\log}) \otimes L \to \FIsoc(X \setminus Z, X) \otimes L$ is fully faithful and its essential image consists of those objects which are docile along $Z$.
\end{defn}

\subsection{Newton polygons}

\begin{defn} \label{D:NP via DM}
For $x \in X$ (not necessarily closed) and $\calE \in \FIsoc(X) \otimes L$, 
we may give an intrinsic definition of the Newton polygon
$N_x(\calE)$ using the Dieudonn\'e--Manin classification \cite[Definition~3.3]{kedlaya-isocrystals}, \cite[Definition~1.2.3]{kedlaya-companions}.
\end{defn}

In terms of this intrinsic definition of Newton polygons, Lemma-Definition~\ref{L:generic Newton polygon} may be refined as follows.
\begin{lemma} \label{L:convergent Newton polygon}
Let $L$ be a finite extension of $\QQ_p$.
For $\calE \in \FIsoc(X) \otimes L$, the function $x \mapsto N_x(\calE)$ of Definition~\ref{D:NP via DM} has the following properties.
\begin{enumerate}
\item[(a)]
For $k$ finite and $x = X = \Spec(k)$, the construction agrees with Lemma-Definition~\ref{L:generic Newton polygon} (a).
\item[(b)]
For $x \in X$ with Zariski closure $Z$, for all $z \in Z$, $N_z(\calE)$ lies on or above $N_x(\calE)$ with the same right endpoint. In particular, the right endpoint is constant on each component of $X$.
\item[(c)]
In (b), the set $U$ of $z \in Z$ for which equality holds is open and Zariski dense,
and its complement is of pure codimension $1$ in $Z$.
\item[(d)]
In (c), for each curve $C$ in $Z$, the inverse image of $U$ in $C$ is either empty or the complement of a finite set of closed points.
\item[(e)]
The vertices of $N_x(\calE)$ all belong to $\ZZ \times \frac{1}{N} \ZZ$ for 
$N = [L:\QQ_p]$.
\end{enumerate}
Consequently, for $\calE \in \FIsoc^\dagger(X) \otimes L$, Definition~\ref{D:NP via DM} agrees with
Lemma-Definition~\ref{L:generic Newton polygon} (on account of the uniqueness assertion in the latter).
\end{lemma}
\begin{proof}
For (a), see \cite[Lemma~1.2.4]{kedlaya-companions}.
For (b) and (c), see \cite[Theorem~3.12]{kedlaya-isocrystals}. Part (d) is an immediate consequence of (c); we include it only for parallelism with Lemma-Definition~\ref{L:generic Newton polygon}.
For part (e), we may assume that $x = X = \Spec(k)$ and then deduce this directly from the Dieudonn\'e--Manin classification (as in \cite[Lemma~1.2.4]{kedlaya-companions} again).
\end{proof}

\begin{defn}
Let $\FIsoc^{\geq 0}(X) \otimes L$ denote the full subcategory of $\FIsoc(X) \otimes L$ consisting of objects whose Newton polygons have nonnegative slopes everywhere. By Lemma~\ref{L:convergent Newton polygon}(b), it suffices to check the slope condition at generic points.
\end{defn}

\subsection{Slope filtrations}

While the Dieudonn\'e--Manin classification does not extend to the case where $X$ is not a point, a weaker version of the statement does generalize as follows. By  Lemma~\ref{L:convergent Newton polygon}(c), if $X$ is irreducible, then the hypothesis on the constancy of the Newton polygon can always be enforced after restricting from $X$ to a suitable open dense subscheme; this provides a crucial link back to \'etale fundamental groups and related concepts (e.g., Chebotaryov density).

\begin{prop} \label{P:slope filtration1}
Suppose that $\calE \in \FIsoc(X) \otimes L$ is such that 
for some $\mu$ and $e$, for all $x \in X$ the least slope of $N_x(\calE)$ is $\mu$ with multiplicity $e$. Then there exists a unique short exact sequence
\[
0 \to \calE_1 \to \calE \to \calE_2 \to 0
\]
such that for all $x \in X$, $N_x(\calE_1)$ has the unique slope $\mu$ with multiplicity $e$.
\end{prop}
\begin{proof}
Apply \cite[Proposition~4.1]{kedlaya-isocrystals}.
\end{proof}

\begin{cor} \label{C:slope filtration}
Suppose that $\calE \in \FIsoc(X) \otimes L$ is such that the function $x \mapsto N_x(\calE)$ is constant.
Let $\mu_1 < \cdots < \mu_l$ be the slopes of $N_x(\calE)$ for any $x \in X$.
Then $\calE$ admits a unique filtration (the \emph{slope filtration})
\[
0 = \calE_0 \subset \cdots \subset \calE_l = \calE
\]
such that for $i=1,\dots,l$, for all $x \in X$, $N_x(\calE_i/\calE_{i-1})$ consists of the single slope $\mu_i$.
\end{cor}
\begin{proof}
Apply \cite[Corollary~4.2]{kedlaya-isocrystals}.
\end{proof}

\begin{lemma} \label{L:no extension filtration}
Suppose that $X$ is irreducible with generic point $\eta$.
Choose $\calE \in \FIsoc(X) \otimes L$ with least slope $\mu$ occurring with multiplicity $e$.
Let $U$ be the open (by Lemma~\ref{L:convergent Newton polygon}(c)) subset of $x \in X$ on which $N_x(\calE)$ has least slope $\mu$
with multiplicity $e$, then form the exact sequence
\[
0 \to \calE_1 \to \calE|_{U} \to \calE_2 \to 0
\]
in $\FIsoc(U) \otimes L$
as in Proposition~\ref{P:slope filtration1}. Then $\calE_1$ does not extend to $\FIsoc(V) \otimes L$ for any strict superset $V$ of $U$ which is open in $X$.
\end{lemma}
\begin{proof}
Suppose that $\calE_1$ extends across some superset $V$ of $U$ which is open in $X$.
By Lemma~\ref{L:fully faithful to convergent}, the canonical embedding $\calE_1 \to \calE$ extends from $U$ to $V$;
this yields an exact sequence
\[
0 \to \calE_1 \to \calE|_{V} \to \calE'_2 \to 0
\]
in which $\calE'_2|_V \cong \calE_2$.
For $x \in V$, apply Lemma~\ref{L:convergent Newton polygon}(b) and (c) to deduce that the slopes of $N_x(\calE'_2)$ are all strictly greater than $\mu$.
It now follows that $N_x(\calE) = N_x(\calE|_V)$ has least slope $\mu$ with multiplicity $e$; hence $x \in U$.
We conclude that $U = V$, proving the claim.
\end{proof}

\begin{remark} \label{R:slope filtration exterior power}
With notation as in Corollary~\ref{C:slope filtration}, 
for $i=1,\dots,l$,
the first step in the slope filtration of $\wedge^{\rank(\calE_{i-1})+1} \calE$ is
\begin{equation} \label{eq:slope filtration wedge}
\wedge^{\rank(\calE_{i-1})+1} \calE_i \cong 
(\calE_i/\calE_{i-1}) \otimes \wedge^{\rank(\calE_{i-1})} \calE_{i-1}.
\end{equation}
\end{remark}

\begin{remark} \label{R:split filtration}
When $X = \Spec k$, one may combine the Dieudonn\'e--Manin decomposition with Galois descent to see that the slope filtration of any $\calE \in \FIsoc(X) \otimes L$ splits uniquely.
A similar statement holds when $X$ is replaced by an arbitrary perfect $\FF_p$-scheme; see 
for example \cite[Lemma~6.11]{kedlaya-drinfeld-isoc}.
\end{remark}

The individual steps of the slope filtration can in turn be interpreted as representations of \'etale fundamental groups.
\begin{defn}
An object $\calE \in \FIsoc(X) \otimes L$ is \emph{unit-root} if for all $x \in X$,
$N_x(\calE)$ has all slopes equal to $0$. By Lemma~\ref{L:convergent Newton polygon}(b), it suffices to check this at the generic point of each irreducible component of $X$.
\end{defn}

\begin{prop} \label{P:unit-root rep}
Suppose that $X$ is irreducible.
Let $L$ be a finite extension of $\QQ_p$.
\begin{enumerate}
\item[(a)]
There is a functorial (in $X$) equivalence of categories between unit-root objects of $\FIsoc(X) \otimes L$
and continuous representations of  $\pi_1(X)$ 
on finite-dimensional $L$-vector spaces.
\item[(b)]
Under this equivalence, unit-root objects of $\FIsoc^\dagger(X) \otimes L$ correspond to representations of $\pi_1(X)$ which are potentially unramified (i.e., after pullback to some finite \'etale cover they become unramified).
\end{enumerate}
\end{prop}
\begin{proof}
See \cite[Theorem~3.7, Theorem~3.9]{kedlaya-isocrystals}.
\end{proof}

\begin{cor} \label{C:rank 1 companions}
Theorem~\ref{T:companion} holds for coefficient objects of rank $1$.
\end{cor}
\begin{proof}
By Lemma~\ref{L:decompose to finite order}, up to an (algebraic) constant twist, any rank-$1$ coefficient object has finite order.
Using Proposition~\ref{P:unit-root rep}(b) in the crystalline case, we may equate coefficient objects of rank 1 with finite order with characters $\chi\colon \pi_1^{\ab}(X) \to L$ of finite order. The claim is now evident.
\end{proof}

\begin{cor} \label{C:unit-root docile}
Let $L$ be a finite extension of $\QQ_p$.
Let $U$ be an open dense subscheme of $X$. If $\calF \in \FIsoc^\dagger(U) \otimes L$ is unit-root and docile along $X \setminus U$, then $\calF \in \FIsoc^\dagger(X) \otimes L$.
\end{cor}
\begin{proof}
By Proposition~\ref{P:unit-root rep}(b), there exists a finite \'etale cover $\tilde{U}$ of $U$ such that $\calF|_{\tilde{U}}$ is everywhere unramified. Consequently, $\calF$ extends across the smooth locus of the normalization of $X$ in $\tilde{U}$;
by faithfully flat descent,
 we deduce that $\calF$ extends across a subset of $X$ whose complement has codimension at least 2. By
Lemma~\ref{L:Zariski purity1}, this proves the claim.
\end{proof}

\begin{defn} \label{D:unit-root subobject}
In case Proposition~\ref{P:slope filtration1} applies with $\mu = 0$, we refer to $\calE_1$ as the \emph{unit-root subobject} of $\calE$.
We also refer to the corresponding representation
$\rho \colon \pi_1(X) \to \GL_{\rank(\calE_1)}(L)$ from Proposition~\ref{P:unit-root rep}(a)
as the \emph{unit-root representation} of $\calE$.
\end{defn}

\subsection{The minimal slope theorem}

While by Lemma~\ref{L:fully faithful to convergent} the inclusion functor $\FIsoc^\dagger(X) \otimes L \to \FIsoc(X) \otimes L$ is fully faithful, it does not in general reflect subobjects; in particular, the slope filtration given by Corollary~\ref{C:slope filtration} does not in general lift back to $\FIsoc^\dagger(X) \otimes L$ even if the Newton polygon is constant on $X$ (to avoid the obstruction coming from Lemma~\ref{L:no extension filtration}). In some sense, one expects rather the opposite;
a recent result of Tsuzuki \cite{tsuzuki-minimal} addresses a question raised in
\cite[Remark~5.14]{kedlaya-isocrystals} to this effect.

\begin{theorem}[Tsuzuki] \label{T:Tsuzuki minimal slope}
Suppose that $X$ is irreducible and that either $\dim(X) = 1$ or $k$ is finite.
Let $\calE, \calF$ be irreducible objects in $\FIsoc^\dagger(X) \otimes L$.
Let $U$ be an open dense subset of $X$ on which the functions
$x \mapsto N_x(\calE), x \mapsto N_x(\calF)$ are constant (this set exists by Lemma~\ref{L:convergent Newton polygon}(c)). 
Let $\calE_1, \calF_1 \in \FIsoc(U) \otimes L$ be the first steps of the slope filtrations of $\calE, \calF$, respectively, according to
Proposition~\ref{P:slope filtration1}. 
\begin{enumerate}
\item[(a)]
Both $\calE_1$ and $\calF_1$ are irreducible in $\FIsoc(U) \otimes L$.
\item[(b)]
If $\calE_1 \cong \calF_1$ in $\FIsoc(U) \otimes L$, then $\calE \cong \calF$ in 
$\FIsoc^\dagger(X) \otimes L$.
\end{enumerate}
\end{theorem}
\begin{proof}
For (a), see \cite[Proposition~6.2]{tsuzuki-minimal} in the case where $\dim(X) = 1$, then deduce the case where $k$ is finite by a Lefschetz slicing argument (e.g., apply \cite[Lemma~3.2.1]{kedlaya-companions}).
For (b), see \cite[Theorem~1.3]{tsuzuki-minimal} in the case where $\dim(X) = 1$ or \cite[Theorem~7.20]{tsuzuki-minimal} in the case where $k$ is finite.
\end{proof}

\begin{cor} \label{C:minimal slope semistable}
Suppose that $X$ is irreducible and $k$ is finite.
Let $\calE$ be a semisimple object in $\FIsoc^\dagger(X) \otimes L$. 
Let $U$ be an open dense subset of $X$ on which the function
$x \mapsto N_x(\calE)$ is constant (this set exists by Lemma~\ref{L:convergent Newton polygon}(c)). 
Then each successive quotient of the slope filtration of $\calE$  is semisimple in $\FIsoc(U) \otimes L$.
\end{cor}
\begin{proof}
Since $\calE$ is semisimple, $G(\calE)^\circ$ is a semisimple algebraic group per Definition~\ref{D:monodromy}.
For each positive integer $i$, there is a natural surjective morphism $G(\calE) \to G(\wedge^i \calE)$,
so $G(\wedge^i \calE)^\circ$ is a semisimple algebraic group and  so $\wedge^i \calE$ is also semisimple.
In light of Remark~\ref{R:slope filtration exterior power},
we may now deduce the claim from Theorem~\ref{T:Tsuzuki minimal slope}(a).
\end{proof}

Although we do not use it here, we mention an extension of Theorem~\ref{T:Tsuzuki minimal slope} by D'Addezio.
This resolves the \emph{parabolicity conjecture} of Crew \cite[p. 460]{crew-mono} and also naturally includes Lemma~\ref{L:fully faithful to convergent}. 
\begin{theorem}[D'Addezio] \label{T:daddezio}
Suppose that $X$ is irreducible (of any dimension) and $k$ is finite.
Let $\calE^\dagger$ be an irreducible object in $\FIsoc^\dagger(X) \otimes L$.
Let $U$ be an open dense subset of $X$ on which the function
$x \mapsto N_x(\calE^\dagger)$ is constant (this set exists by Lemma~\ref{L:convergent Newton polygon}(c)).
Let $\calE$ be the image of $\calE^\dagger$ in $\FIsoc(U) \otimes L$.
Then $G(\calE)$ is the subgroup of $G(\calE^\dagger)$ preserving the slope filtration (Proposition~\ref{P:slope filtration1}).
Moreover, if $\calE^\dagger$ is semisimple, then $G(\calE)$ is a parabolic subgroup of $G(\calE^\dagger)$.
\end{theorem}
\begin{proof}
Using Lemma~\ref{L:restriction is fully faithful}, we may reduce to the case $X = U$. In this case, apply \cite[Theorem~1.1.1]{daddezio}. (Note that this statement does not require $k$ finite;
we assumed this only because we introduced monodromy groups under this restriction.)
\end{proof}

\section{Rigid and dagger geometry}

In preparation for our later arguments, we recall some constructions in rigid analytic geometry and its sibling theory, dagger analytic geometry, together with some consequences for convergent and overconvergent isocrystals.

\subsection{Local lifting by formal schemes}
\label{subsec:local lifting}

In order to give a sufficiently concrete description
of isocrystals for our work, we need to work with local lifts of varieties from characteristic
$p$ to characteristic 0. We start with a version of the story which is not sufficient for our
purposes, but is a natural starting point. (Our terminology here is not standard.)

\begin{defn} \label{D:smooth lift}
By a \emph{smooth lift} of $X$, we will mean a smooth formal scheme $P$ over $W(k)$ (for the $p$-adic topology) with $P_k \cong X$.
\end{defn}

\begin{lemma}\label{L:compare charts}
Suppose that $X$ is affine. Then $X$ admits a smooth lift, which is unique up to noncanonical isomorphism.
\end{lemma}
\begin{proof}
This may be obtained by the method of Elkik \cite{elkik}, or more precisely by a result of Arabia
\cite[Th\'eor\`eme~3.3.2]{arabia}.
\end{proof}

\begin{lemma} \label{L:functorial smooth lift}
Let $\overline{f}\colon X' \to X$ be an \'etale morphism and let $P$ be a smooth lift of $X$. Then $\overline{f}$ lifts functorially
to an \'etale morphism $f\colon P' \to P$ of formal schemes over $W(k)$, where $P'$ is a certain smooth lift of $X'$ (determined by $P$ and $\overline{f}$).
\end{lemma}
\begin{proof}
This is a consequence of the henselian property of the pair $(W(k), pW(k))$. See
for example \cite[Theorem~5.5.7]{gabber-ramero}, which is written in the more general context of almost commutative algebra, but is nonetheless a good reference for this point.
\end{proof}

\begin{defn}
Let $P$ be a smooth lift of $X$. A \emph{Frobenius lift} on $P$ is a morphism $\sigma\colon P \to P$ which acts on $W(k)$
via the Witt vector Frobenius and lifts the absolute Frobenius morphism on $X$.
\end{defn}

The following construction gives a particular class of Frobenius lifts.

\begin{defn} \label{D:chart}
Let $(X,Z)$ be a smooth pair over $k$. A \emph{smooth chart} for $(X,Z)$ is a sequence $\overline{t}_1,\dots,\overline{t}_n$ of elements of $\calO_X(X)$ such that the induced morphism $\overline{f}\colon X \to \AAA^n_k$ is \'etale
and there exists $m \in \{0,\dots,n\}$ for which the zero loci of $\overline{t}_1,\dots, \overline{t}_m$ on $X$ are
the irreducible components of $Z$ (and in particular are reduced).
\end{defn}

\begin{lemma} \label{L:make smooth chart}
Let $(X,Z)$ be a smooth pair over $k$. Then for each $x \in X$, there exist an open subscheme $U$ of $X$ containing $x$ and a smooth chart for $(U, Z \cap U)$.
\end{lemma}
\begin{proof}
By replacing $x$ with a specialization, we may assume that $x \in X^\circ$.
Since $X$ is smooth, it satisfies the Jacobian criterion; we can thus find elements $\overline{t}_1,\dots,\overline{t}_n \in \calO_{X,x}$ such that $d\overline{t}_1,\dots,d\overline{t}_n$ form a basis of 
$\Omega_{X/k, x}$ over $\calO_{X,x}$. By adjusting the choice of coordinates, we may further ensure that $Z$ is cut out locally at $x$ by $\overline{t}_1 \cdots \overline{t}_m$.
Choose an open affine neighborhood $U$ of $x$ in $X$ omitting every irreducible component of $Z$ not passing through $x$. Then $\overline{t}_1,\dots,\overline{t}_n$ form a smooth chart for $(U, Z \cap U)$.
\end{proof}

\begin{defn} \label{D:lifted smooth chart}
Let $(X,Z)$ be a smooth pair and let $\overline{t}_1,\dots,\overline{t}_n$ be a smooth chart for $(X,Z)$.
Let $P_0$ be the formal completion of $\Spec W(k)[t_1,\dots,t_n]$ along the zero locus of $p$. 
By Lemma~\ref{L:functorial smooth lift}, there exists a unique smooth affine formal scheme $P$ over $W(k)$ equipped with an \'etale morphism $f\colon P \to P_0$ lifting $\overline{f}$;
we refer to $(P,t_1,\dots,t_n)$ as the \emph{lifted smooth chart} associated to the original smooth chart.

Let $\sigma_0\colon P_0 \to P_0$ be the Frobenius lift for which $\sigma_0^*(t_i) = t_i^p$ for $i=1,\dots,n$.
By the functoriality aspect of Lemma~\ref{L:functorial smooth lift}, there exists a unique Frobenius lift $\sigma$ on $P$ making the diagram
\[
\xymatrix{
P \ar^{\sigma}[r] \ar^f[d] &P \ar^f[d] \\
P_0 \ar^{\sigma_0}[r] & P_0
}
\]
commute. We call $\sigma$ the \emph{associated Frobenius lift} of the lifted smooth chart.
\end{defn}

\begin{defn}
For $P$ a smooth lift of $X$, denote by $P_K$ the Raynaud generic fiber of $P$; this is a rigid analytic space
whose points correspond to formal subschemes of $P$ which are integral and finite flat over $W(k)$
(see \cite[\S 7.4]{bosch}).
In particular, there is a \emph{specialization map} taking any such point to the intersection of $X$ with the corresponding formal subscheme. For $S \subseteq P$, let $]S[_P$ denote the inverse image of $S$ under the specialization map; this set is called the \emph{tube} of $S$ within $P_K$.
\end{defn}

\subsection{Weak formal schemes and dagger spaces}

We next upgrade the previous discussion to handle overconvergent isocrystals, by introducing certain analogues of
formal schemes and rigid analytic spaces which directly incorporate the notion of overconvergence.
These can also be handled in terms of \emph{partially overconvergent isocrystals} (compare Definition~\ref{D:partially overconvergent isocrystals}), but as this description will not appear explicitly we will not spell it out.

\begin{defn}
A ring $R$ is \emph{weakly complete} with respect to an ideal $I$ if it is $I$-adically separated and, for any positive real numbers $a,b$  and any elements $x_1,\dots,x_n \in R$, any infinite sum of the form
\[
\sum_{i_1,\dots,i_n=0}^\infty a_{i_1 \cdots i_n} x_1^{i_1}\cdots x_n^{i_n} \qquad
(a_{i_1\cdots i_n} \in I^{\lfloor a(i_1 +\cdots + i_n) - b \rfloor})
\]
converges in $R$. (By contrast, if $R$ is complete with respect to $I$, then the sum also converges under the weaker condition that $a_{i_1 \cdots i_n} \in I^{f(i_1,\dots,i_n)}$ where $f$ is any function for which $f(i_1,\dots,i_n) \to \infty$ as $i_1 + \cdots + i_n \to \infty$.)
By replacing complete rings with weakly complete rings, we obtain Meredith's concept of a \emph{weak formal scheme}
\cite{meredith}; there is an obvious forgetful functor from weak formal schemes to formal schemes.

By a \emph{dagger lift} of $X$, we will mean a smooth weak formal scheme $P^\dagger$ over $\Spf W(k)$ with $P^\dagger_k \cong X$. We write $P$ for the underlying formal scheme of $P^\dagger$.
\end{defn}

\begin{lemma} \label{L:dagger excellent}
Let $P^\dagger$ be a smooth weak formal scheme over $\Spf W(k)$. Then
$P^\dagger[p^{-1}]$ is noetherian and excellent. (The same is true of $P^\dagger$, but we will not need this fact.)
\end{lemma}
\begin{proof}
This reduces at once to the case where $P^\dagger$ is the weak completion of an affine space.
Since $K$ is of characteristic 0, this case is an easy consequence of the 
Nullstellensatz for dagger algebras (see \cite[\S 1.4]{grosse-klonne})
plus the weak Jacobian criterion in the form of \cite[Theorem~102]{matsumura}.
\end{proof}

\begin{lemma} \label{L:dagger faithfully flat}
For $P^\dagger$ a smooth weak formal scheme over $\Spf W(k)$, the morphism $P \to P^\dagger$ is faithfully flat.
\end{lemma}
\begin{proof}
Since we have surjectivity on points, it suffices to check flatness of the morphism on coordinate rings.
This morphism is the direct limit of a family of morphisms, each of which corresponds to an open immersion of affinoid spaces over $K$ and so is flat \cite[Corollary~7.3.2/6]{bgr}.
\end{proof}

\begin{lemma} \label{L:functorial dagger lift}
Let $\overline{f}\colon X' \to X$ be an \'etale morphism and let $P^\dagger$ be a dagger lift of $X$. Then $\overline{f}$ lifts functorially
to an \'etale morphism $f^\dagger\colon P^{\prime\dagger} \to P^\dagger$ of weak formal schemes, 
where $P^{\prime \dagger}$ is a certain dagger lift of $X'$ (determined by $P^\dagger$ and $\overline{f}$).
\end{lemma}
\begin{proof}
Note that if $R$ is weakly complete with respect to $I$, then the pair $(R,I)$ is henselian. We may thus argue as in the proof of Lemma~\ref{L:functorial smooth lift}.
\end{proof}

\begin{defn} \label{D:weak completion}
With notation as in Definition~\ref{D:lifted smooth chart},
let $P_0^\dagger$ be the weak formal completion of $\Spec W(k)[t_1,\dots,t_n]$ along the zero locus of $p$.
By Lemma~\ref{L:functorial dagger lift}, $f$ descends uniquely to an \'etale morphism $f^\dagger\colon P^\dagger \to P_0^\dagger$ of weak formal schemes, and $\sigma$ descends to a morphism $\sigma^\dagger\colon P^\dagger \to P^\dagger$.
We refer to $(P^\dagger,t_1,\dots,t_n)$ as the \emph{lifted dagger chart} associated to the original smooth chart.
\end{defn}

\begin{defn} \label{D:dagger space}
For $P^\dagger$ a dagger lift of $X$, we may again define the \emph{generic fiber} $P^\dagger_K$
as a locally G-ringed space with the same underlying G-topological space as $P_K$, but with a modified structure
sheaf. The space $P^\dagger_K$ lives in the category of \emph{dagger spaces} of Grosse-Kl\"onne
\cite{grosse-klonne}; to summarize, these are built in a fashion analogous to rigid analytic spaces, but with the role of standard Tate algebras (i.e., the coordinate rings of generic fibers of formal completions
of affine spaces over $W(k)$) being played by their overconvergent analogues (in which the formal completions
become weak formal completions).
\end{defn}

\subsection{Relative GAGA for rigid and dagger spaces}

It is well known that Serre's GAGA theorem for varieties over $\CC$ \cite{serre-gaga}
has an analogue over a nonarchimedean field, in which complex analytic spaces are replaced by rigid analytic spaces.
Here we formulate a version of this result in a relative setting.
(At this point we are only using the fact that $K$ is a nonarchimedean field; the discreteness of the valuation plays no role.)

\begin{prop} \label{P:GAGA}
Let $\calC$ denote either the category of rigid analytic spaces over $K$ or the category of dagger spaces over $K$.
Suppose that $S$ is an affinoid space in $\calC$, put $S_0 := \Spec \calO(S)$,
and let $\pi_S\colon S \to S_0$ be the adjunction morphism.
Let $f_0\colon Y_0 \to S_0$ be a proper morphism, let $f\colon Y \to S$ be the analytification of $f_0$ in $\calC$, and let $\pi_Y\colon Y \to Y_0$ be the adjunction morphism.
\begin{enumerate}
\item[(a)]
The morphisms $\pi_S, \pi_Y$ are flat and every closed point has a unique preimage.

\item[(b)]
Pullback along the induced morphism $Y \to Y_0$ defines an equivalence of categories between coherent sheaves on $Y_0$ and on $Y$.

\item[(c)]
Let $\calE_0$ be a coherent sheaf on $Y_0$ and let $\calE$ be the pullback of $\calE_0$ to $Y$. Then the natural morphisms $\pi^* (R^i f_{0*} \calE_0) \to R^i f_* \calE$ of sheaves on $S$ are isomorphisms for all $i \geq 0$.
\end{enumerate}
\end{prop}
\begin{proof}
To prove (a), we need only treat the case $S \to S_0$.
For this, see \cite[Theorem~1.7]{grosse-klonne} for the dagger case, and references therein for the rigid-analytic case.

We next skip to (c). 
Using Chow's lemma \stacktag{0200} as in the complex-analytic case, we may reduce to the case where $f_0$ is projective,
and then further to the case where $Y_0 = \PP^n_{S_0}$.
By (a), pullback of coherent sheaves along $Y \to Y_0$ is an exact functor;
we are thus free to make homological reductions. Using the relative ampleness of $\calO(1)$ with respect to $\PP^n_{S_0} \to S_0$, we may reduce to the cases where $\calE = \calO(m)$ for $m \in \ZZ$.
These cases amount to the fact that the morphism
\[
\calO(S_0)[T_1^{\pm},\dots,T_n^{\pm}] \to \calO(S_0) \langle T_1^{\pm},\dots,T_n^{\pm} \rangle
\]
in the rigid case, and the morphism
\[
\calO(S_0)[T_1^{\pm},\dots,T_n^{\pm}] \to \calO(S_0) \langle T_1^{\pm},\dots,T_n^{\pm} \rangle^\dagger
\]
in the dagger case, induces isomorphisms of associated graded rings for the grading by homogeneous degree.
(This is essentially the method of Serre; see the first paragraph of \cite[Example~3.2.6]{conrad} for another approach that directly handles the case where $f_0$ is proper, without having to add Chow's lemma.)

To prove (b), using (a) and (c) we know that the pullback functor is fully faithful.
To establish essential surjectivity, we may again assume that $Y_0 = \PP^n_{S_0}$;
again using (c), it suffices to check that
for every coherent sheaf $\calF$ on $\PP^n_S$, there exists an integer $m$ such that
$\calF(m)$ is generated by global sections.
In the rigid-analytic case, we may appeal to \cite[Theorem~3.2.4]{conrad} to deduce this immediately.
In the dagger case, we may make the corresponding argument after we note that as in Kiehl's theorem,
one knows that the higher direct images of a coherent sheaf along a proper morphism of dagger spaces
are again coherent \cite[Theorem~3.5]{grosse-klonne}.
\end{proof}

\begin{remark}
Proposition~\ref{P:GAGA}, specialized to the case where $S$ is a point, includes the usual GAGA theorem
for rigid analytic spaces over $K$, which has been known for some time (see the discussion in \cite[Example~3.2.6]{conrad}).
It also includes GAGA for dagger spaces over $K$, but this is not new either because the category of proper dagger spaces over $K$ is equivalent to the category of proper rigid analytic spaces over $K$ \cite[Theorem~2.27]{grosse-klonne}.
\end{remark}

As an application, we record a description of docile overconvergent isocrystals in the case where $X$ has a liftable smooth compactification.

\begin{defn}
By a \emph{smooth lift} of a smooth pair $(\overline{X}, Z)$ over $k$, we will mean a smooth pair
$(\overline{\frakX}, \frakZ)$ of schemes (not formal schemes) over $W(k)$ equipped with compatible identifications
$\overline{\frakX}_k \cong \overline{X}$, $\frakZ_k \cong Z$.
\end{defn}

\begin{prop} \label{P:tame realization}
Let $(\overline{\frakX}, \frakZ)$ be a smooth lift of a smooth pair $(\overline{X},Z)$ with $X \cong \overline{X} \setminus Z$.
Let $L$ be a finite extension of $\QQ_p$. 
\begin{enumerate}
\item[(a)]
There is a fully faithful functor
from the category of objects of $\Isoc(X, \overline{X}) \otimes L$ (isocrystals without Frobenius structure) which are docile along $Z$
to the category of vector bundles on the Raynaud generic fiber $\overline{\frakX}_K$
equipped with an integrable $K$-linear logarithmic
(with respect to $\frakZ_K$) connection with nilpotent residues
and a compatible $\QQ_p$-linear action of $L$.
\item[(b)]
Suppose that $\overline{X}$ is proper over $k$.
Then the functor in (a) factors through the category of vector bundles on the \emph{scheme} $\overline{\frakX}_K$
equipped with an integrable $K$-linear logarithmic
(with respect to $\frakZ_K$) connection with nilpotent residues 
and a compatible $\QQ_p$-linear action of $L$.
\end{enumerate}
\end{prop}
\begin{proof}
We may reduce at once to the case $L = \QQ_p$.
By Proposition~\ref{P:GAGA}, it suffices to prove the corresponding assertion with
$(\overline{\frakX}, \frakZ)$ replaced by its $p$-adic formal completion. In that case,
apply \cite[Theorem~6.4.1]{kedlaya-semi1} to obtain a fully faithful functor from the latter category
to vector bundles on $\overline{\frakX}_K$ equipped with logarithmic
(with respect to $\frakZ_K$) integrable connections with nilpotent residues.
\end{proof}

\begin{remark} \label{R:no Frobenius lift}
In Proposition~\ref{P:tame realization}, one can get an equivalence of categories if $(\overline{\frakX}, \frakZ)$ admits a Frobenius lift, in which case the vector bundles also carry an action of this lift compatible with the connection).
In practice such a lift almost never exists, so the Frobenius structure has to be described locally using smooth lifts of affine subschemes of $\overline{X}$.
\end{remark}

\begin{remark} \label{R:same convergence condition}
The difference between a convergent (resp. overconvergent) isocrystal on $X$ and a vector bundle with integrable connection on $P_K$ for some smooth lift $P$ (resp. on $P^\dagger_K$ for some dagger lift $P^\dagger$) is a certain convergence condition on the formal Taylor isomorphism. For our purposes, the only relevant aspect of this condition is that given a vector bundle with integrable connection on $P^\dagger_K$, the convergence condition can be checked after pullback to $P_K$; this is a consequence of Baldassarri's theorem on continuity of the generic radius of convergence for a $p$-adic differential equation (see for example \cite[Chapter~25]{kedlaya-course}).
\end{remark}

\subsection{Restriction to divisors in a curve fibration}
\label{subsec:restriction to divisors}

We next record some applications of Lemma~\ref{L:faithful restriction Vec} in the context of isocrystals.
\begin{hypothesis}
Throughout \S\ref{subsec:restriction to divisors},
let $f\colon \overline{X} \to S$ be a smooth curve fibration of smooth $k$-schemes with pointed locus $Z$ and unpointed locus $X$.
Let $\eta$ be the generic point of $S$.
Let $D$ be a divisor in $X$ which dominates $S$.
Let $L$ be a finite extension of $\QQ_p$.
\end{hypothesis}

\begin{defn} \label{D:lift smooth fibration}
Suppose that $S$ admits a smooth dagger lift $\frakS^\dagger$
such that $f$ lifts to a smooth curve fibration $\frakf^\dagger\colon \overline{\frakX}^\dagger \to \Spec \calO(\frakS^\dagger)$ with pointed locus $\frakZ^\dagger$.
Let $\frakS$ be the underlying formal scheme of $\frakS^\dagger$.
Let $\frakf\colon \overline{\frakX} \to \Spec \calO(\frakS)$
be the pullback of $\frakf^\dagger$ from $\Spec \calO(\frakS^\dagger)$ to $\Spec \calO(\frakS)$,
as a smooth curve fibration with pointed locus $\frakZ$.

By adapting the proof of Proposition~\ref{P:tame realization},
we obtain a fully faithful realization functor from the full subcategory of $\Isoc^\dagger(X) \otimes L$ consisting of objects docile along $Z$ to the category of vector bundles on the scheme
$\overline{\frakX}_K^\dagger$ equipped with an integrable $K$-linear logarithmic (with respect to $\frakZ_K^\dagger$) connection with nilpotent residues and a compatible $\QQ_p$-linear action of $L$.
Similarly, we obtain a fully faithful realization functor from the full subcategory of $\Isoc(\overline{X}^{\log}) \otimes L$ consisting of objects docile along $Z$ to the category of vector bundles on the scheme
$\overline{\frakX}_K$ equipped with an integrable $K$-linear logarithmic (with respect to $\frakZ_K$) connection with nilpotent residues and a compatible $\QQ_p$-linear action of $L$.

In case $D \to S$ is finite \'etale, by Lemma~\ref{L:functorial dagger lift} we may lift $D$ uniquely to a closed subscheme $\frakD^\dagger$ of $\overline{\frakX}^\dagger$ which is finite \'etale over $\Spec \calO(\frakS^\dagger)$ with $\deg(\frakD^\dagger \to \Spec \calO(\frakS^\dagger)) = \deg(D \to S)$.
\end{defn}

\begin{prop} \label{P:descent using divisor}
For any positive integer $n$, there exists a positive integer $d$ (depending only on $f$ and $n$)
such that whenever $\deg(D \times_S \eta \to \eta) \geq d$,
the restriction functor from the full subcategory of $\FIsoc^\dagger(X) \otimes L$ consisting of objects docile along $Z$ to
\[
(\FIsoc(\overline{X}^{\log}) \otimes L) \times_{\FIsoc(D) \otimes L}
(\FIsoc^\dagger(D) \otimes L)
\]
restricts to an equivalence of categories on objects of rank $\leq n$.
\end{prop}
\begin{proof}
By Lemma~\ref{L:fully faithful to convergent}, we may check the claim \'etale-locally on $S$. 
Moreover, by Lemma~\ref{L:Zariski purity1}, we may check the claim after replacing $S$ by an open dense subscheme.
Hence by Proposition~\ref{P:stable reduction}, we may assume
$f$ lifts as in Definition~\ref{D:lift smooth fibration}.
We may also check the claim after shrinking $S$ using \cite[Corollary~5.9]{kedlaya-isocrystals}; we may thus assume that $D \to S$ is finite \'etale.

Since the map $\calO(\frakS^\dagger) \to \calO(\frakS)$ is faithfully flat (Lemma~\ref{L:dagger faithfully flat}),
we may apply Lemma~\ref{L:faithful restriction Vec}(a) (to reflect faithfully flat descent data in the sense of Remark~\ref{R:reflect fpqc descent data})
and Remark~\ref{R:same convergence condition}
to show that the restriction functor from the full subcategory of $\Isoc^\dagger(X) \otimes L$ consisting of objects docile along $Z$ to
\[
(\Isoc(\overline{X}^{\log}) \otimes L) \times_{\Isoc(D) \otimes L}
(\Isoc^\dagger(D) \otimes L)
\]
restricts to an equivalence of categories on objects of rank $\leq n$.
By adding Frobenius structures, we may formally promote this assertion to the desired result.
\end{proof}

\begin{prop} \label{P:split Azuyama algebra using divisor}
For any positive integer $n$, there exists a positive integer $d$ (depending only on $f$ and $n$)
such that whenever $\deg(D \times_S \eta \to \eta) \geq d$,
the following statement holds.
Let $\calG \in \FIsoc(\overline{X}^{\log}) \otimes L$ be an Azumaya algebra object of degree $n$ in the sense of 
Definition~\ref{D:Azumaya algebra}, i.e., the image of the standard representation under some $L$-linear tensor functor
$\Rep_L(\PGL_n) \to \FIsoc(\overline{X}^{\log}) \otimes L$.
Suppose further that $\calG|_D$ admits a splitting $\calF_D$ in $\FIsoc(D) \otimes L$
(that is, $\calG|_D \cong \calF_D^\dual \otimes \calF_D$).
Then there exists $\calF \in \FIsoc(\overline{X}^{\log}) \otimes L$ such that $\calG \cong \calF^\dual \otimes \calF$ and $\calF|_D \cong \calF_D$; in fact, $\calF$ is unique up to unique isomorphism.
\end{prop}
\begin{proof}
With notation as in Definition~\ref{D:lift smooth fibration} (and invoking Remark~\ref{R:tannakian to azumaya}),
we may realize $\calF$ as an Azumaya algebra on the scheme $\overline{\frakX}_K \times_{\QQ_p} L$.
We first apply Corollary~\ref{C:faithful restriction split Azumaya}
to lift the splitting of this Azumaya algebra from $\frakD_K \times_{\QQ_p} L$ to
$\overline{\frakX}_K \times_{\QQ_p} L$, then again to promote the splitting to the level of vector bundles with connection.
Finally, we repeat the argument with $\calF$ and $\calG$ replaced by their Frobenius pullbacks 
to see that the Frobenius structure on $\calF|_D$ extends to $\calF$.
\end{proof}

We will need at one point the following variant of  Proposition~\ref{P:descent using divisor}.
\begin{defn}
For any $\FF_p$-scheme $T$, let $T^{\perf}$ denote the perfect closure of $T$,
i.e., the inverse limit of
\[
\cdots \to T \stackrel{\varphi_T}{\to} T
\]
where $\varphi_T$ is the absolute Frobenius on $T$.
\end{defn}

\begin{prop} \label{P:descent using divisor perfect}
For $n,d,D$ as in Proposition~\ref{P:descent using divisor},
the restriction functor from $\FIsoc(\overline{X}^{\log}) \otimes L$ to
\begin{equation} \label{eq:fiber using divisor2}
(\FIsoc(\overline{X}^{\log} \times_S \eta^{\perf}) \otimes L) \times_{\FIsoc(D \times_S \eta^{\perf}) \otimes L}
(\FIsoc(D) \otimes L)
\end{equation}
restricts to an equivalence of categories on objects of rank $\leq n$.
\end{prop}
\begin{proof}
As the claim is again local on $S$ and may be checked after replacing $S$ with an open dense subscheme
(by Lemma~\ref{L:fully faithful to convergent} and Lemma~\ref{L:Zariski purity1}), we may assume that $D \to S$ is finite \'etale and
that $\frakS$ admits a smooth lifted chart as in Definition~\ref{D:lifted smooth chart}, and hence a corresponding Frobenius lift $\sigma$.
We then obtain a Frobenius-equivariant map of $p$-adic formal schemes $W(S^{\perf}) \to \frakS$
by making the canonical identification of the $p$-adic completion of 
$\varinjlim_\sigma \calO(\frakS)$ with $W(\calO(S^{\perf}))$.
The $\calO(\frakS)$-module $\varinjlim_\sigma \calO(\frakS)$ is free on the basis
\[
\{\sigma^{-j}(t_i)\colon i=1,\dots,n; j=0,1,\dots\};
\]
it follows directly from this that $W(S^{\perf}) \to \frakS$ is faithfully flat.

We define an \emph{ad hoc} category $\FIsoc(\overline{X}^{\log}/S) \otimes L$ of ``convergent $F$-isocrystals on $\overline{X}^{\log}$ relative to $S$.'' We first define $\Isoc(\overline{X}^{\log}/S)$ to be the category of vector bundles on $\frakX_K$ equipped with a $\frakS_K$-linear logarithmic (with respect to $\frakZ_K$) connection with nilpotent residues
such that the pullback along $W(\eta^{\perf})  \to W(S^{\perf}) \to \frakS$ gives an object of 
$\Isoc(\overline{X}^{\log} \times_S \eta^{\perf})$; this last condition imposes a convergence
condition as in Remark~\ref{R:same convergence condition}. 
We then observe that the convergence condition allows us to define a Frobenius pullback
endofunctor on $\Isoc(\overline{X}^{\log}/S)$ using Frobenius lifts extending $\sigma$
(compare Lemma~\ref{L:lattice change of Frobenius lift} and Definition~\ref{D:Frobenius pullback});
we then define $\FIsoc(\overline{X}^{\log}/S)$ 
by equipping objects of $\Isoc(\overline{X}^{\log}/S)$ with Frobenius structures.
We finally tensor with $L$ in the sense of   \cite[Definition~9.2]{kedlaya-isocrystals}.

We first verify that the functor from $\FIsoc(\overline{X}^{\log}) \otimes L$ to 
\begin{equation} \label{eq:fiber using divisor3}
(\FIsoc(\overline{X}^{\log}/S) \otimes L) \times_{\FIsoc(D/S) \otimes L} (\FIsoc(D) \otimes L)
\end{equation}
is an equivalence.
(Note that $\Isoc(D/S)$ is just the category of vector bundles on $\frakD_K$.)
Given an object of the target, we must show that the $\frakS_K$-linear connection extends uniquely to an integrable $K$-linear connection; that is, we must provide actions of derivations on the base and check that these are compatible with each other and with the Frobenius and $L$-linear structures. 
All of this follows by applying Lemma~\ref{L:faithful restriction Vec}(a).

We next define the category $\FIsoc((\overline{X}^{\log} \times_S S^{\perf})/S^{\perf}) \otimes L$
by analogy with $\FIsoc(\overline{X}^{\log}/S) \otimes L$.
By Lemma~\ref{L:faithful restriction Vec}(a) again plus the fact that $W(S^{\perf}) \to \frakS$ is faithfully flat, the functor from $\FIsoc(\overline{X}^{\log}/S) \otimes L$ to 
\[
(\FIsoc((\overline{X}^{\log} \times_S S^{\perf})/S^{\perf}) \otimes L) \times_{\FIsoc(D^{\perf}/S^{\perf}) \otimes L} (\FIsoc(D/S) \otimes L)
\]
is an equivalence. (Note that $D^{\perf} =D \times_S S^{\perf}$ because $D \to S$ is finite \'etale.)
Consequently, the functor from
\eqref{eq:fiber using divisor3} to
\begin{equation} \label{eq:fiber using divisor4}
(\FIsoc((\overline{X}^{\log} \times_S S^{\perf})/S^{\perf}) \otimes L) \times_{\FIsoc(D^{\perf}/S^{\perf}) \otimes L} (\FIsoc(D) \otimes L)
\end{equation}
is an equivalence. 

Finally, by Lemma~\ref{L:faithful restriction Vec}(b), the functor from $\FIsoc((\overline{X}^{\log} \times_S S^{\perf})/S^{\perf}) \otimes L$ to
\[
(\FIsoc(\overline{X}^{\log} \times_S \eta^{\perf}) \otimes L) \times_{\FIsoc(D \times_S \eta^{\perf}) \otimes L}
(\FIsoc(D^{\perf}/S^{\perf}) \otimes L)
\]
is fully faithful, as then is
the functor from \eqref{eq:fiber using divisor4} to \eqref{eq:fiber using divisor2}.
As for essential surjectivity, we deduce something slightly weaker: given an object of  \eqref{eq:fiber using divisor2}, we obtain from Lemma~\ref{L:faithful restriction Vec}(b) a lift to $\FIsoc(\overline{X}^{\log}) \otimes L$ after replacing $S$ by the complement of a subset of codimension $\geq 2$.
However, as noted above we are already free to make such a replacement.
\end{proof}

\subsection{Crystalline lattices}
As our approach to constructing crystalline companions involves building these objects up as an inverse limit over truncations modulo powers of $p$, we need some statements about lattices in $F$-isocrystals.

\begin{defn} \label{D:lattice}
Let $L$ be a finite extension of $\QQ_p$.
Let $(X,Z)$ be a smooth pair with $X$ affine 
admitting a smooth chart, let $(P, t_1,\dots,t_n)$ be the resulting lifted smooth chart, and let $\sigma$ be the associated Frobenius lift.
For $\calE \in \FIsoc(X^{\log}) \otimes L$, as per Proposition~\ref{P:tame realization} we realize $\calE$
as a vector bundle on $P_K$ equipped with a logarithmic $K$-linear integrable connection with nilpotent residues and compatible actions of $L$ and $\sigma$.

By a \emph{lattice} in $\calE$, we will mean a vector bundle $\frakE$ on $P$ equipped with an isomorphism of $\calE$ with the pullback of $\frakE$ to $P_K$. Note that this concept is not functorial in $X$, as it depends on the choice of the lifted smooth chart.
\end{defn}

The interaction of differentiation with the Frobenius structure
has the following consequence (compare Remark~\ref{R:same convergence condition}).
\begin{lemma} \label{L:spectral norm}
With notation as in Definition~\ref{D:lattice},
suppose that $X$ is irreducible with generic point $\eta$.
Let $\calE_\eta$ be the generic point of $X$,
as a finite-dimensional vector space over the fraction field of a Cohen ring for the residue field of $\eta$.
Equip this vector space with the supremum norm defined by some basis (or equivalently by some lattice). Then
for $i=1,\dots,n$, the spectral norm of $\frac{\partial}{\partial t_i}$ on $\calE_\eta$ equals $p^{-1/(p-1)} < 1$.
\end{lemma}
\begin{proof}
See \cite[Theorem~17.2.1]{kedlaya-course}.
\end{proof}

\begin{defn} \label{D:crystalline lattice}
With notation as in Definition~\ref{D:lattice},
a lattice $\frakE$ in $\calE$ is \emph{crystalline} if it is stable under the connection, the action of $\frako_L$, and the action of $\sigma$. This concept by contrast is functorial in $X$
by Lemma~\ref{L:lattice change of Frobenius lift} below;
consequently, we can extend the definition of crystalline lattices to the case where $X$ is not necessarily affine.
\end{defn}

\begin{lemma} \label{L:lattice change of Frobenius lift}
With notation as in Definition~\ref{D:crystalline lattice}, suppose that $\frakE$ is a crystalline lattice in $\calE$. Let $\sigma'$ be any other Frobenius lift on $P$ (not necessarily associated to a lifted smooth chart). Then there is a Frobenius structure on $\frakE$ with respect to $\sigma'$ given by the Taylor series
\[
\sigma'(\bv) = \sum_{i_1,\dots,i_n=0}^\infty \frac{(\sigma'(t_1) - t_1^p)^{i_1} \cdots (\sigma'(t_n) - t_n^p)^{i_n}}{i_1! \cdots i_n!} 
\sigma \left( \frac{\partial^{i_1}}{\partial t_1^{i_1}} \cdots 
\frac{\partial^{i_1}}{\partial t_n^{i_n}} \bv \right).
\]
\end{lemma}
\begin{proof}
By Lemma~\ref{L:spectral norm},
 the sum converges in $\calE$ and defines a Frobenius structure (compare \cite[Proposition~17.3.1]{kedlaya-course}), so it suffices to check that each summand belongs to $\frakE$. This follows from the fact that $\sigma'(t_j) - t_j^p$ is divisible by $p$, 
so
$(\sigma'(t_j) - t_j^p)^{i_j}/(i_j!)$ is $p$-adically integral (and even divisible by $p$).
\end{proof}

\begin{remark} \label{R:use crystalline lattice}
The existence of a crystalline lattice in $\calE \in \FIsoc(X^{\log}) \otimes L$
clearly implies that $\calE \in \FIsoc^{\geq 0}(X^{\log}) \otimes L$. It is not known whether the reverse implication holds except when $\dim(X) = 1$; see Lemma~\ref{L:crystalline lattice1} and Lemma~\ref{L:crystalline lattice} below.
\end{remark}

\begin{lemma} \label{L:crystalline lattice1}
For any $\calE \in \FIsoc^{\geq 0}(X^{\log}) \otimes L$,
there exists an open dense subspace $U$ of $X$ such that the restriction of $\calE$ to $\FIsoc(U) \otimes L$ admits a crystalline lattice.
\end{lemma}
\begin{proof}
We may assume that $X$ admits a smooth chart. With notation as above, there exists a nonnegative integer $c$ such that each of the following operations carries $\frakE$ into $p^{-c} \frakE$:
\begin{enumerate}
\item[(i)]
multiplication by every element of $\frako_L$;
\item[(ii)]
every power of the Frobenius structure $\varphi$ with respect to $\sigma$;
\item[(iii)]
every power of $\frac{\partial}{\partial t_i}$ for $i=1,\dots,n$.
\end{enumerate}
Namely, (i) is obvious because $\frako_L$ is finite over $\ZZ_p$; (ii) is a consequence of the hypothesis on $N_\eta(\calE)$ (see \cite[Sharp Slope Estimate 1.5.1]{katz-slope});
and (iii) follows from Lemma~\ref{L:spectral norm}.

At this point, we can construct a new lattice $\frakE'$ in $\calE_\eta$ which is stable under multiplication by $\frako_L$, $\sigma$, and $\frac{\partial}{\partial t_i}$ for $i=1,\dots,n$:
namely, the images of $\frakE$ under all combinations of these operations span a sublattice of $p^{-(n+2)c} \frakE$ with the desired effect. Over some $U$, we may spread $\frakE'$ out to a lattice in $\calE$ with the desired effect.
\end{proof}

\begin{remark} \label{R:constant twist of dual lattice}
With notation as in Definition~\ref{D:crystalline lattice}, suppose that $\frakE$ is a crystalline lattice in $\calE$.
Let $c$ be an integer greater than or equal to the sum of the Newton slopes of $N_x(\calE)$ for every $x \in X$, or equivalently by Lemma~\ref{L:convergent Newton polygon}(b) for the generic point of every irreducible component of $X$. Then the constant twist by $p^c$ of the induced Frobenius structure on $\calE^\dual$ defines a Frobenius structure on $\frakE^\dual$.
\end{remark}

\section{Uniformities for isocrystals on curves}
\label{sec:uniformities}

We next assemble some standard statements about vector bundles on curves, then use these to control the Harder--Narasimhan polygons of vector bundles arising from isocrystals on curves.
This yields some key uniformity statements needed
for the construction of crystalline companions.

\subsection{Harder--Narasimhan polygons}
\label{subsec:HN}

We recall some standard properties of vector bundles on curves, using the language of Harder--Narasimhan (HN) polygons.

\begin{hypothesis} \label{H:curve over field}
Throughout \S\ref{subsec:HN}--\ref{subsec:gaps between slopes}, let $C$ be a proper curve of genus $g$ over an arbitrary field $L$ (of arbitrary characteristic).
\end{hypothesis}

\begin{defn} 
We define the \emph{degree} of a nonzero vector bundle on $C$, denoted $\deg(E)$, as the degree of its top exterior power,
setting the degree of the zero bundle to be 0. 
Define the \emph{slope} of a nonzero vector bundle $E$ on $C$, denoted $\mu(E)$, as the degree divided by the rank.
We record some basic properties.
\begin{itemize}
\item
If $F$ is a subbundle of $E$ of the same rank, then $\deg(F) \leq \deg(E)$. (By taking top exterior powers, this reduces to the case where
the common rank is 1.)
\item
If $0 \to F \to E \to G \to 0$ is a short exact sequence of vector bundles,
then 
\[
\deg(F) + \deg(G) = \deg(E), \qquad \rank(F) + \rank(G) = \rank(E).
\]
In particular,
the inequalities
\[
\mu(F) < \mu(E), \qquad \mu(F) = \mu(E), \qquad \mu(F) > \mu(E)
\]
are equivalent respectively to
\[
\mu(G) > \mu(E), \qquad \mu(G) = \mu(E), \qquad \mu(G) < \mu(E).
\]
\item
For $E^\dual$ the dual bundle of $E$, we have
$\deg(E^\dual) = -\deg(E)$ and (if $E$ is nonzero) $\mu(E^\dual) = -\mu(E)$.
\end{itemize}
\end{defn}

\begin{defn} \label{D:HN filtration}
Let $E$ be a nonzero vector bundle over $C$.
We say that $E$ is \emph{semistable} if $\mu(E)$ is maximal among the slopes of nonzero subbundles of $E$.
We say that $E$ is \emph{stable} if ($E$ is semistable and) $\mu(F) < \mu(E)$ for every nonzero strict subbundle $F$ of $E$.

By a standard argument (e.g., see \cite[\S 3.4]{kedlaya-aws}), 
every vector bundle $E$ on $C$ admits a unique filtration
\[
0 = E_0 \subset \cdots \subset E_l = E
\]
by subbundles such that each quotient $E_i/E_{i-1}$ is a vector bundle which is semistable of some slope $\mu_i$
and $\mu_1 > \cdots > \mu_l$; this is called the \emph{Harder--Narasimhan filtration}, or \emph{HN filtration}, of $E$. The associated \emph{Harder--Narasimhan polygon}, or \emph{HN polygon}, of $E$ is the graph of the function $[0,\rank(E)] \to \RR$ which interpolates linearly between the points
$(\rank(E_i), \deg(E_i))$ for $i=0,\dots,l$.
\end{defn}

\begin{remark} \label{R:upper convex hull}
A useful characterization of the HN polygon of $E$ is as the boundary of the upper convex hull of the 
set of points $(\rank(F), \deg(F)) \in \RR^2$ as $F$ varies over all subobjects $F$ of $E$.
The proof is easy; see for example \cite[Lemma~3.4.15]{kedlaya-aws}.
\end{remark}

\begin{remark}
Observe that the dual bundle $E^\dual$ is (semi)stable if and only if $E$ is.
Consequently, the slopes of the HN polygon of $E^\dual$ are the negations of the HN slopes of $E$.
\end{remark}

\begin{lemma} \label{L:semistable base change}
Let $F$ be an algebraic closure of $L$. Then a vector bundle on $C$ is semistable if and only if its pullback to $C_F$ is semistable.
\end{lemma}
\begin{proof}
See \cite[Proposition~3]{langton}.
\end{proof}

\begin{lemma} \label{L:homs by slope}
Let $E_1, E_2$ be vector bundles on $C$. Suppose that every slope of the HN polygon of 
$E_1$ is strictly greater than every slope of the HN polygon of $E_2$. Then $\Hom(E_1, E_2) = 0$.
\end{lemma}
\begin{proof}
Observe first that if $E_1$ and $E_2$ are semistable and $E_1 \to E_2$ is a nonzero morphism with image $F$, then $F$ is both a quotient of $E_1$ and a subobject of $E_2$, so $\mu(E_1) \leq \mu(F) \leq \mu(E_2)$. In particular, if $\mu(E_1) > \mu(E_2)$, then $\Hom(E_1, E_2) = 0$.

We use this as the base case of an induction on $\rank(E_1) + \rank(E_2)$; we may also assume that $E_1, E_2$ are both nonzero as else there is nothing to check. Away from the base case, one of $E_1$ and $E_2$ is not semistable. If $E_1$ is not semistable, then let $F$ be the first step of its HN filtration; given any homomorphism $E_1 \to E_2$, the induction hypothesis implies first that the restriction to $F$ is zero, and second that the induced homomorphism $E_1/F \to E_2$ is zero. A similar argument applies if $E_2$ is not semistable.
\end{proof}

\begin{remark} \label{R:RR to generators}
By comparison with Lemma~\ref{L:homs by slope}, one may ask whether if every slope of the HN polygon of 
$E_1$ is strictly less than every slope of the HN polygon of $E_2$, then $\Hom(E_1, E_2) \neq 0$.
This is false; for example, if $E$ is a semistable vector bundle on $C$ of positive slope, it is not necessarily true that $H^0(C, E) \neq 0$. 

That said, using Riemann--Roch, one can prove some weaker statements in this direction. For example,
suppose that $\rank(E) = 1$. Then $\dim_L H^0(C, E) \geq \deg(E) - g + 1$, so if $\deg(E) \geq g$ then $H^0(C, E) \neq 0$.
See also Lemma~\ref{L:HN to sections} below.
\end{remark}

\begin{lemma} \label{L:HN to sections}
Let $E$ be a vector bundle on $C$ such that every slope of the HN polygon of $E$ is greater than $2g-1$. Then 
$E$ is generated by global sections.
\end{lemma}
\begin{proof}
By Lemma~\ref{L:semistable base change},
we may reduce to the case where $L$ is algebraically closed.
By Serre duality,
\begin{equation} \label{eq:Serre duality1}
\dim_L H^0(C, E) - \dim_L H^0(C, \Omega_{C/L} \otimes E^\dual) = \rank(E)(\mu(E) + 1-g).
\end{equation}
Since $\Omega_{C/L} \otimes E^\dual$ has all slopes less than $-1$, Lemma~\ref{L:homs by slope} implies that
$H^0(C, \Omega_{C/L} \otimes E^\dual(P)) = 0$ for each $P \in C^\circ$ (here we use that $L$ is now algebraically closed). By comparing \eqref{eq:Serre duality1} with the analogous equation
\[
\dim_L H^0(C, E(-P)) - \dim_L H^0(C, \Omega_{C/L} \otimes E^\dual(P)) = \rank(E)(\mu(E)-g),
\]
we deduce that the trivial inequality
\[
\dim_L H^0(C, E) \leq \dim_L H^0(C, E(-P)) + \rank(E)
\]
is an equality; this is only possible if the fiber of $E$ at $P$ is generated by global sections of $E$. This proves the claim.
\end{proof}

\begin{lemma} \label{L:HN in extension}
Let $0 \to E_1 \to E \to E_2 \to 0$ be an exact sequence of vector bundles on $C$.
\begin{enumerate}
\item[(a)]
The HN polygon of $E$ lies on or above the union of the HN polygons of $E_1$ and $E_2$ 
(that is, the polygon in which each slope occurs with multiplicity equal to the sum of its multiplicities in the
HN polygons of $E_1$ and $E_2$), with the same endpoint.
\item[(b)]
Suppose that every slope of the HN polygon of $E_1$ is strictly less than every slope of the HN polygon of $E_2$.
Then equality holds in (a) if and only if the sequence splits.
\end{enumerate}
\end{lemma}
\begin{proof}
It is easy to see that if $E = E_1 \oplus E_2$, then the HN polygon of $E$ is equal to the union
(see \cite[Lemma~3.4.13]{kedlaya-aws}); one then deduces (a) using Remark~\ref{R:upper convex hull}
as in \cite[Lemma~3.4.17]{kedlaya-aws}. Note that this also proves the ``if'' direction of (b).
As for the ``only if'' direction of (b), 
suppose that equality holds; then the HN filtration of $E$ admits a step $E'_1$ of rank equal to $\rank(E_2)$. Put $E'_2 := E/E'_1$; by Lemma~\ref{L:homs by slope}, the induced map
$E'_1 \to E_1$ is zero, so we must have $E'_1 \cong E_2$. This gives the desired splitting.
\end{proof}

\begin{remark} \label{R:semicontinuity}
Given a curve $\frakC$ over a scheme $S$ and a vector bundle $\frakE$ over $\frakC$,
for each point $s \in S$ we may pull back $\frakE$ to the fiber of $\frakC$ over $s$ and compute the HN polygon of the resulting bundle; this defines a function from $S$ to the space of polygons which is upper semicontinuous \cite[Theorem~3]{shatz}. For example, taking $S = \Spec W(k)$, this statement can be used to transfer bounds on HN polygons from the special fiber to the generic fiber.
(Note that the right endpoint of the HN polygon is preserved under specialization, and so is locally constant.)
\end{remark}

\subsection{Gaps between slopes}
\label{subsec:gaps between slopes}

In many cases, one can use extra structure on, or special properties of, a vector bundle on a curve to deduce constraints on the gaps between consecutive HN slopes, and thus on the HN polygon. We collect some statements of this form here; see \S\ref{subsec:isocrystal slopes} for a further argument in this vein.
(Recall that Hypothesis~\ref{H:curve over field} remains in effect.)

\begin{defn}
For $E$ a nonzero vector bundle on $C$, 
we define the \emph{width} of $E$,
denoted $\width(E)$,
 to be the maximum of $\left| \mu(E) - \mu_i \right|$ as $\mu_i$ varies over the slopes of the HN polygon of $E$.
 (This is not standard terminology.)
\end{defn}

\begin{remark} \label{R:width in filtration}
Let $E$ be a nonzero vector bundle on $C$, consider a strictly increasing filtration $0 = E_0 \subset \cdots \subset E_l = E$ of $E$ by saturated subbundles, and put $F_i := E_i/E_{i-1}$ for $i=1,\dots,l$. Then 
using Lemma~\ref{L:HN in extension}, we see that
\begin{align*}
\width(E) &\leq \max\{\width(F_i)\colon i\in \{1,\dots,n\}\}  \\
&\qquad 
+ \max\{ \left| \mu(F_i) - \mu(F_j) \right|\colon i,j \in \{1,\dots,n\}\}.
\end{align*}
In particular, if 
\[
0 \to E_1 \to E \to E_2 \to 0
\]
is an exact sequence of nonzero vector bundles on $C$ and $\mu(E_1) = \mu(E_2)$, then
$\width(E) \leq \max\{\width(E_1), \width(E_2)\}$.
\end{remark}

\begin{remark} \label{R:consecutive slopes}
Let $E$ be a vector bundle of degree $d$ and rank $r$ on $C$. Let $c$ be a real number with the property that no two consecutive slopes of the HN polygon of $E$ differ by more than $c$.
Then 
\[
\frac{d}{r} = \mu(E) \geq \frac{1}{r} \sum_{j=0}^{r-1} (\mu_1 - jc)
= \mu_1 - (r-1)c/2;
\]
by this plus the corresponding argument for $E^\dual$, we see that $\width(E) \leq (r-1)c/2$.
This remark is applicable in a variety of situations where one can bound the difference between consecutive slopes;
see Lemma~\ref{L:indecomposable} and Lemma~\ref{L:Frobenius pullback slopes}
for examples, and Proposition~\ref{P:isocrystal reduction slopes} for a closely related argument.
\end{remark}

\begin{lemma} \label{L:indecomposable}
For any indecomposable vector bundle $E$ on $C$, the differences between consecutive slopes of $E$ are bounded by $2g-2$.
\end{lemma}
\begin{proof}
See  \cite[Proposition~2.1]{schiffmann}.
\end{proof}

\begin{lemma} \label{L:Frobenius pullback slopes}
Suppose that $L$ is of characteristic $p$ and let $\varphi_C$ denote the absolute Frobenius on $C$.
Then for any semistable vector bundle $E$ on $C$, the differences between consecutive slopes of $\varphi_C^* E$ are bounded by $2g-2$.
\end{lemma}
\begin{proof}
See  \cite[Corollary~2]{shepherd-barron}.
\end{proof}

\begin{cor} \label{C:Frobenius pullback slopes}
Suppose that $L$ is of characteristic $p$ and let $\varphi_C$ denote the absolute Frobenius on $C$.
Let $E$ be a vector bundle of rank $r$ on $C$.
Then 
\[
\width(\varphi_C^* E) \leq 2p \width(E) + (r-1)(g-1).
\]
\end{cor}
\begin{proof}
Let $0 = E_0 \subset \cdots \subset E_l = E$ be the HN filtration of $E$
and put 
\[
F_i := \varphi_C^*(E_i)/\varphi_C^*(E_{i-1}) \cong \varphi_C^*(E_i/E_{i-1}).
\]
For each $i$, by Remark~\ref{R:consecutive slopes} and Lemma~\ref{L:Frobenius pullback slopes}
we have 
\[
\width(F_i) \leq (r-1)(g-1).
\]
On the other hand, we have
\[
\max\{|\mu(F_i) - \mu(F_j)|\} \leq
\mu(F_1) - \mu(F_l) = 
p \mu(E_1) - p \mu(E_{l-1}) \leq 2p \width(E).
\]
By Remark~\ref{R:width in filtration}, this yields the desired result.
\end{proof}

\begin{lemma} \label{L:connection degree 0}
 Let $D$ be a finite union of closed points of $C$.
Let $E$ be a vector bundle on $C$ admitting a logarithmic connection with singularities in $D$ for which all residues are nilpotent.
\begin{enumerate}
\item[(a)]
If $L$ is of characteristic $0$, then $\deg(E) = 0$.
\item[(b)]
If $L$ is of characteristic $p$, then $\deg(E) \equiv 0 \pmod{p}$.
\item[(c)]
Suppose that $L = k$.
Let $n$ be a positive integer and let $\frakC$ be a smooth proper curve over $W_n(k)$ with $\frakC_k \cong C$.
Let $\frakD$ be a closed subscheme of $\frakC$, smooth over $W_n(k)$, with $\frakD_k = D$ as subschemes of $C$.
Suppose that $E$ is the pullback to $C$ of a vector bundle on $\frakC$ admitting a logarithmic connection with singularities in $\frakD$ for which all residues are nilpotent. Then $\deg(E) \equiv 0 \pmod{p^n}$.
\end{enumerate}
\end{lemma}
\begin{proof}
By taking the top exterior power, we may reduce all three cases to the setting where $\rank(E) = 1$; in this case,
the connection has no singularities. By choosing a rational section $\bv$ of $E$, we may write the connection as 
$\bv \mapsto \bv \otimes \omega$ where $\omega$ is a meromorphic differential. By a local computation,
the image in $L$ of the degree of $E$ equals the sum of the residues of $\omega$, which equals 0 by the residue theorem. This implies (a) and (b); we may similarly deduce (c) by repeating the computation on $\frakC$.
\end{proof}

\begin{remark}
Part (a) of Lemma~\ref{L:connection degree 0} also follows from a formula of Ohtsuki \cite[Theorem~3]{ohtsuki}: the existence of the connection implies that $E$ has vanishing first Chern class, and hence is of degree $0$. For $D = \emptyset$, the corresponding statement dates back to Atiyah \cite[Theorem~4]{atiyah}.

More precisely, the Chern class formula of Ohtsuki involves a sum over poles of the connection, in which each summand involves the Chern polynomial evaluated at the residue. The Chern polynomial depends on its argument only up to semisimplification, so nilpotent residues contribute as if they were zero; the Chern class formula thus returns  zero, just as it would in the case of an everywhere holomorphic connection.

What this shows is that some restriction on residues in Lemma~\ref{L:connection degree 0} is essential to obtaining any meaningful restriction on $E$.
In fact, without the residue restriction, \emph{every} vector bundle on $C$ admits a logarithmic connection; see for example \cite[Theorem~5.3]{biswas-heu}.

On the other hand, we can still obtain a bound of the desired form if we insist that the exponents be not necessarily zero, but belong to some bounded interval such as $[0,1)$. 
\end{remark}

\begin{cor} \label{C:connection degree 0}
With notation as in Lemma~\ref{L:connection degree 0}(a),
the slopes of the HN polygon of $E$ are bounded in absolute value by $(\rank(E)-1)(g-1)$.
\end{cor}
\begin{proof}
We may assume at once that $E$ is indecomposable. By
Remark~\ref{R:consecutive slopes} and Lemma~\ref{L:indecomposable},
$\width(E) \leq (\rank(E) - 1)(g-1)$. By Lemma~\ref{L:connection degree 0}, this yields the desired result.
\end{proof}

\begin{lemma} \label{L:width bound generation}
Fix an ample line bundle $\calO(1)$ on $C$. Let $E$ be a vector bundle on $C$. Then for every integer 
$n > \width(E) - \mu(E) + 2g-1$, $E(n)$ is generated by global sections.
\end{lemma}
\begin{proof}
By definition, $\mu(E) - \width(E) + n$ is a lower bound on the slopes of $E(n)$.
Since this is greater than $2g-1$, Lemma~\ref{L:HN to sections} implies that
$E(n)$ is generated by global sections.
\end{proof}

\begin{remark}
In Definition~\ref{D:HN filtration}, suppose that $L$ is of characteristic $0$.
Using transcendental methods, specifically the Narasimhan--Seshadri theorem \cite{narasimhan-seshadri} relating stable vector bundles on a compact Riemann surface to irreducible unitary representations of the fundamental group),
one sees that the tensor product of two semistable vector bundles on $C$ is again semistable
(compare \cite[Proposition~2.2]{hartshorne} for the corresponding argument for symmetric powers).

By contrast, if $L$ is of positive characteristic, then this assertion fails in general (unless one factor is a line bundle); the first counterexamples are due to Gieseker
\cite[Corollary~1]{gieseker}. There are various ways to control the damage caused by this fact; see
Lemma~\ref{L:tensor product width} for a simple argument that suffices for our purposes.
\end{remark}

\begin{lemma} \label{L:tensor product width}
For any two vector bundles $E_1, E_2$ on $C$, 
\[
\width(E_1 \otimes E_2) \leq \width(E_1) + \width(E_2) + 4g.
\]
\end{lemma}
\begin{proof}
By Lemma~\ref{L:semistable base change}, we may assume that $L$ is algebraically closed, so that there exists an ample
line bundle $\calO(1)$ of degree $1$ on $C$.
Put $n_i = \lfloor \width(E_i) - \mu(E_i) + 2g \rfloor$.
By Lemma~\ref{L:width bound generation}, $E_i(n_i)$ is generated by global sections, as then is
$(E_1 \otimes E_2)(n_1 + n_2)$; the latter bundle therefore has all slopes nonnegative.
This gives the claimed upper bound on the slopes of $E_1 \otimes E_2$; the claimed lower bound follows by a similar argument applied to the duals of $E_1, E_2$.
\end{proof}

\begin{remark}
While here we only need to consider the Frobenius pullback $\varphi_C^*$, we mention in passing that the pushforward $\varphi_{C*}$ preserves semistability
\cite[Theorem~1.1]{mehta-pauly} and stability \cite[Theorem~2.2]{sun}.
\end{remark}

\subsection{Isocrystals on curves and slopes}
\label{subsec:isocrystal slopes}

We now specialize to curves and combine with the theory of vector bundles to obtain the crucial uniformity.
In passing, we mention that some related results have been obtained by Esnault--Shiho \cite{esnault-shiho2}
and Bhatt \cite{bhatt-lurie} from a somewhat different point of view.

\begin{hypothesis} \label{H:isocrystals slopes}
Let $(\overline{X},Z)$ be a smooth pair in which $\overline{X}$ 
is a proper curve over $k$, and put $X := \overline{X} \setminus Z$.
Let $\eta$ be the generic point of $X$. Let $g$ be the genus of $\overline{X}$ and let $m$ be the $k$-length of $Z$.
Let $\varphi\colon X \to X$ denote the absolute Frobenius morphism on $X$.

Fix a smooth lift $(\overline{\frakX}, \frakZ)$ of $(\overline{X},Z)$
(which exists because of the smoothness of the moduli stack of curves;
see Proposition~\ref{P:stable reduction})
and put $\frakX := \overline{\frakX} \setminus \frakZ$.
Fix a finite extension $L$ of $\QQ_p$ of degree $d$ and an object $\calE \in \FIsoc^\dagger(X) \otimes L$ of rank $r$
which is docile along $Z$. 
As per Proposition~\ref{P:tame realization}, we realize $\calE$
as a vector bundle on $\overline{\frakX}_K$
 equipped with a logarithmic integrable connection with nilpotent residues and a compatible action of $L$.
\end{hypothesis}

\begin{defn} \label{D:generic lattice on curve}
Denote by $W(\eta)$ the completed local ring of 
$\overline{\frakX}$ at $\eta$; this is not the ring of Witt vectors over $\eta$,
but rather a Cohen ring with residue field $\eta$.
Let $\calE_\eta$ denote the pullback of $\calE$ to $\Spec W(\eta)[p^{-1}]$; this is a finite-dimensional vector space over $W(\eta)[p^{-1}]$ equipped with a connection and a $\QQ_p$-linear action of $L$.
Moreover, for any Frobenius lift $\varphi$ on $W(\eta)$, $\calE_\eta$ admits a semilinear $\varphi$-action
on $\calE_\eta$ compatible with the connection and the $L$-action.

By a \emph{lattice} in $\calE$, we will mean a vector bundle $\frakE$ on $\overline{\frakX}$
equipped with an isomorphism of $\calE$ with the pullback of $\frakE$ to $\overline{\frakX}_K$. 
There is a pullback functor from lattices in $\calE$ to lattices in $\calE_\eta$.
There is also a pullback functor from lattices in $\calE$ to vector bundles on $X$, which we call the 
\emph{reduction} functor.
Since $\calE$ has degree $0$ by Lemma~\ref{L:connection degree 0}, any reduction of $\calE$ also has degree $0$ by Remark~\ref{R:semicontinuity}.
\end{defn}

\begin{lemma} \label{L:crystalline lattice}
There exists a crystalline lattice in $\calE$ if and only if $N_\eta(\calE)$ has all slopes nonnegative (that is, $\calE \in \FIsoc^{\dagger,\geq 0}(X) \otimes L$).
\end{lemma}
\begin{proof}
If $\calE$ admits a crystalline lattice, then it is clear that $N_x(\calE)$ has all slopes nonnegative for
each $x \in X^\circ$; by Lemma~\ref{L:convergent Newton polygon}(c), $N_\eta(\calE)$ has all slopes
nonnegative.
Conversely, suppose that $N_\eta(\calE)$ is nonnegative. 
To produce a crystalline lattice in $\calE$, 
apply Lemma~\ref{L:crystalline lattice1} to produce a lattice in $\calE_\eta$ stable under the connection, Frobenius, and $\frako_L$-action; this lattice extends uniquely to a reflexive sheaf on $\frakX$, which is a vector bundle because $\frakX$ is regular of dimension 2 \stacktag{0B3N}.
\end{proof}

\begin{remark} \label{R:reduction degree}
By Corollary~\ref{C:connection degree 0}, the slopes of the vector bundle $\calE$ on $\overline{\frakX}_K$
are bounded in absolute value by $([L:\QQ_p] \rank(\calE)-1)(g-1)$. However, this does not directly imply anything
about the slopes of a reduction of $\calE$, because the semicontinuity property of HN polygons goes in the wrong direction for this (see Remark~\ref{R:semicontinuity}). In fact, we do not claim anything about the slopes of an arbitrary reduction; rather, we will produce specific lattices whose reductions we can control. 
\end{remark}

\begin{prop} \label{P:isocrystal reduction slopes}
Suppose that $N_\eta(\calE)$ has all slopes nonnegative.
Then $\calE$ admits a crystalline lattice whose reduction has width at most $dr(3g+m)$.
Moreover, since the reduction has degree $0$ (see Lemma~\ref{L:connection degree 0}), this is also a bound on the absolute value of every HN slope of the reduction.
\end{prop}
\begin{proof}
To simplify notation, we assume at once that $d=1$ (at the expense of replacing $r$ with $rd$).
In light of Remark~\ref{R:width in filtration} and
Lemma~\ref{L:connection degree 0}, 
we may also assume that $\calE$ is irreducible of rank $r>1$.
Under these assumptions, we define a sequence of vector bundles
$\frakE_0, \frakE_1, \dots$ on $\overline{\frakX}$ corresponding to crystalline lattices of $\calE$; show that this sequence must terminate; and show that the terminal value has the desired property. 

We start by applying Lemma~\ref{L:crystalline lattice} to construct $\frakE_0$.
Given $\frakE_i$ for some $i$, form the HN filtration of its reduction $\frakE_{i,k}$. 
Suppose that there exists a pair of consecutive slopes $\mu_j, \mu_{j+1}$ of the HN polygon of $\frakE_{i,k}$
such that
\[
\mu_j - \mu_{j+1} > 2g-2+m, \qquad p\mu_j - \mu_{j+1} > (r-1)(g-1).
\]
By Lemma~\ref{L:homs by slope} plus the first inequality, the Kodaira--Spencer morphism 
$F_j \to \frakE_{i,k}/F_j \otimes \Omega_{X/k}(Z)$ must vanish; thus $F_j$ is stable under the induced connection on $\frakE_{i,k}$. 
Similarly, by Corollary~\ref{C:Frobenius pullback slopes} plus the second inequality, the morphism $\varphi^* F_j \to \frakE_{i,k}/F_j$ must vanish, so $F_j$ is stable under the induced Frobenius action on $\frakE_{i,k}$.
We may thus extend the sequence by taking $\frakE_{i+1}$ to be the inverse image of $F_j$ under the projection
$\frakE_i \to \frakE_{i,k}$, and this will again be a crystalline lattice of $\calE$.

By construction, there exists an exact sequence
\[
0 \to \frakE_{i,k}/F_j \to \frakE_{i+1,k} \to F_j \to 0.
\]
By Lemma~\ref{L:HN in extension}, the HN polygon of $\frakE_{i+1,k}$ is bounded below by the HN polygon of $\frakE_{i,k}$ with the same endpoints, and equality only holds if this sequence splits.
In particular, the HN polygons of the $\frakE_{i,k}$ form a discrete, monotone, bounded sequence; this sequence must therefore either terminate or stabilize. 
If the sequence stabilizes, then the intersection of the $\frakE_i$ forms
a strict subbundle of $\frakE$ stable under the Frobenius and the connection, contrary to our hypothesis that $\calE$ is irreducible.
We thus have the desired termination.

It remains to show that termination implies the indicated bound. Let $\frakE_i$ be the terminal value of the sequence.
Let
$\mu_1 \geq \dots \geq \mu_r$ be the slopes of the HN polygon of $\frakE_{i,k}$ listed with multiplicity. 
Put 
\[
c := \frac{(r-1)(g-1) - p(2g-2+m)}{p-1}.
\]
If $\mu_{j+1} \geq c$, then
\[
p\mu_j - \mu_{j+1} \leq (r-1)(g-1) \Longrightarrow
p\mu_j - p\mu_{j+1}  \leq p(2g-2+m)
\]
and so termination implies  $\mu_j - \mu_{j+1} \leq 2g-2+m$.
If instead $\mu_{j+1} < c$, then
\[
p\mu_j - p\mu_{j+1}  \leq p(2g-2+m)
\Longrightarrow
p\mu_j - \mu_{j+1} \leq (r-1)(g-1) 
\]
and so termination implies $p\mu_j - \mu_{j+1} \leq (r-1)(g-1)$;
in particular, 
\[
\mu_j < \frac{c + (r-1)(g-1)}{p} = \frac{(r-3)(g-1) - m}{p-1}.
\]
If there exists $j \in \{1,\dots,r-1\}$ such that $\mu_{j+1} < c$, then for the smallest such $j$ 
we obtain
\begin{align*}
\mu_1 &= \mu_j + \sum_{h=1}^{j-1} (\mu_h - \mu_{h+1}) \\
&< \frac{(r-3)(g-1) - m}{p-1} + (r-2)(2g-2+m) \\
&\leq \frac{rg}{p-1} + r(2g+m) \leq r(3g+m).
\end{align*}
Otherwise $0 \geq \mu_r \geq c$ (because $\mu_1 + \cdots + \mu_r = 0$ and $r>1$)
and
\begin{align*}
\mu_1 &= \mu_r + \sum_{i=1}^{r-1} (\mu_i - \mu_{i+1}) \leq (r-1)(2g-2+m) \leq r(3g+m).
\end{align*}
We thus have an upper bound on $\mu_1$ of the desired form.
Repeating the argument for the dual of $\calE$ gives a lower bound on $\mu_r$ of the same form.
\end{proof}

\begin{remark}
In Proposition~\ref{P:isocrystal reduction slopes},
the space of lattices of $\calE$ can be interpreted as a Bruhat--Tits building for $\GL_r(L)$; the space of crystalline lattices is a closed bounded subspace of the building; and the proof of Proposition~\ref{P:isocrystal reduction slopes} amounts to finding a lattice ``near the center'' of the closed bounded subspace.
\end{remark}

\begin{remark}
Note that in Proposition~\ref{P:isocrystal reduction slopes}, for fixed $p,r,d$, the slope bound is linear in $g$ and $m$.
For a similar bound on the length of the jumping locus, see Corollary~\ref{C:uniform bound on jumping}.

For our present purposes, there is no need to optimize the bound (and so we instead optimized for simplicity of the derivation). We only use that the bound is a function of $p,r,g,m,d$ alone
(and in fact uniform in $p$, but this plays no role here). 
\end{remark}

\subsection{Jumping loci on curves}
\label{subsec:jumping loci}

As a complement to Proposition~\ref{P:isocrystal reduction slopes}, we derive a statement about the jumping loci of Newton polygons on curves by adapting an argument of Tsuzuki \cite{tsuzuki}. 

\begin{hypothesis} \label{H:jumping1}
Throughout \S\ref{subsec:jumping loci}, 
assume that $k$ is finite,
and define $X, \overline{X}, Z, g, m, L$ as in Hypothesis~\ref{H:isocrystals slopes}.
Let $\calE$ be an $E$-algebraic coefficient object for some number field $E$, fix a place of $E$ above $p$,
and let $q$ be the order of the residue field of this place.
\end{hypothesis}

We derive an analogue of \cite[Lemma~3.1]{tsuzuki} (but see Remark~\ref{R:jump bound discrepancy} for the differences).

\begin{prop} \label{P:uniform jump bound1}
Suppose that $\calE$ is tame of rank $r$
and that $N_\eta(\calE)$ has least slope $0$ with multiplicity $1$.
Let $W$ be the finite set of $x \in X$ for which $N_x(\calE)$ has nonzero least slope.
Then
\[
\deg_k(W) \leq (q-1) r (2g- 2 + m)-m+1.
\]
\end{prop}
\begin{proof}
By Corollary~\ref{C:algebraic companion}, we may reduce to the crystalline case.
Let $L_0$ be the $p$-adic completion of $E$.

Put $U = X \setminus W$; by Proposition~\ref{P:slope filtration1}, the restriction of $\calE$ to
$\FIsoc(U) \otimes L$ admits a unit-root subobject $\calE_1$ of rank $1$.
For each $x \in X^\circ$, Hensel's lemma implies that $P_x(\calE,T)$ admits a slope factorization over $L_0$;
consequently, $P_x(\calE_1, T) \in L_0[T]$. By this fact plus Proposition~\ref{P:unit-root rep}(a),
$\calE_1$ corresponds to a character $\chi\colon \pi_1^{\ab}(U) \to L_0^\times$,
which by compactness must take values in $\frako_{L_0}^\times$.
Let $\varpi$ be a uniformizer of $L_0$; we then have $\chi^{q-1} \equiv 1 \pmod{\varpi}$.

For $x \in X^\circ$, we have
\[
P(\calF_x, T) \equiv \begin{cases} 1 - T^{[\kappa(x):k]} \pmod{\varpi} & (x \in U) \\ 1 \pmod{\varpi} & (x \in W); \end{cases}
\]
consequently, at the level of $L$-functions we have
\[
L(\calE^{\otimes (q-1)}, T) \equiv L(\calO_U, T) \pmod{\varpi}.
\]
We will derive the desired bound by analyzing both sides of this congruence. On one hand,
by the Grothendieck--Ogg--Shafarevich formula \cite[Proposition~2.4.2]{kedlaya-companions}, $L(\calE^{\otimes (q-1)}, T)$ is a rational function of degree
$(q-1)r (2g-2+m)$; its reduction modulo $\varpi$ thus has degree at most $(q-1)r (2g-2+m)$. On the other hand, 
\[
L(\calO_U, T) = L(\calO_{\overline{X}}, T)\prod_{x \in W \cup Z} (1 - T^{[\kappa(x):k]})
\]
and the first factor equals $(1-T)^{-1}(1-\#(k)T)^{-1}$ times a polynomial with constant term 1;
after reducing modulo $\varpi$, we obtain a rational function of degree at least $m-1 + \deg_k(W)$. Comparing these estimates yields the desired result.
\end{proof}

\begin{remark} \label{R:jump bound discrepancy}
Proposition~\ref{P:uniform jump bound1} differs from \cite[Lemma~3.1]{tsuzuki} in two key respects.
One is that Tsuzuki does not make explicit the uniformity in the genus (this being irrelevant in his use case), although his method yields it immediately. The other is that Tsuzuki does not allow logarithmic singularities (these being immaterial for his purposes), but again this makes no serious difference to the argument.
\end{remark}

We now recover a bound as in \cite[Theorem~3.3]{tsuzuki}.

\begin{cor} \label{C:uniform bound on jumping}
Suppose that $\calE$ is tame of rank $r$.
Let $W$ be the finite set of $x \in X$ for which $N_x(\calE) \neq N_\eta(\calE)$.
Then $\deg_k(W) \leq (q-1)(2^r-2) (2g-2+m)-m+1$.
\end{cor}
\begin{proof}
For each $s \in \{1,\dots,r-1\}$ for which $N_\eta(\calE)$ has a vertex with $x$-coordinate $s$, apply Proposition~\ref{P:uniform jump bound1} to $\wedge^s \calE$ to deduce that
the set of $x \in X$ for which $N_x(\calE)$ omits this vertex
has $k$-length at most $(q-1) \binom{r}{s} (2g-2+m)-m+1$. Summing over $s$ yields the claim.
\end{proof}

\begin{remark} \label{R:Igusa polynomial}
We have not made any effort to optimize the upper bound in Corollary~\ref{C:uniform bound on jumping}. Some useful test cases come from Shimura varieties,
where the variation of Newton polygons can be studied quite explicitly. 

For a concrete example, assume $p>2$, take $X := \PP^1_k \setminus \{0,1,\infty\}$ with the coordinate $\lambda$, and let $\calE$ be the middle cohomology of the Legendre pencil $y^2 = x(x-1)(x-\lambda)$ of elliptic curves; then $r=2$, $g=0$, $m=3$, $q = p$,
and $\deg_k(W) = \frac{p-1}{2}$ since the Igusa polynomial is squarefree of this degree
e.g., see \cite[Theorem~4.1]{silverman}).
\end{remark}

\begin{remark}
We expect that Proposition~\ref{P:uniform jump bound1} can be extended to the case where $k$ is not necessarily finite (and there is no algebraicity condition).
This will yield a corresponding extension of Corollary~\ref{C:uniform bound on jumping} as a corollary.

A result of the same type, but with a worse upper bound,
can be obtained using Proposition~\ref{P:isocrystal reduction slopes}.
Since we will not need this here, we take the liberty of sketching the argument only when $L = \QQ_p$.
Choose a crystalline lattice $\frakE$ whose reduction $\frakE_k$ has width bounded as in Proposition~\ref{P:isocrystal reduction slopes}.
We can then interpret the jumping locus as the support of the torsion subsheaf of $\coker((\varphi^r)^* \frakE_k \to \frakE_k)$,
so its degree is bounded by $\deg((\varphi^r)^* \frakE_k^\dual \otimes \frakE_k)$.
We bound the latter in terms of $\width(\frakE_k)$
using Lemma~\ref{C:Frobenius pullback slopes} and Lemma~\ref{L:tensor product width}.
\end{remark}

\begin{remark}
It is also a question of great interest to give lower bounds on the degree of the jumping locus. One difficulty with this is to control multiplicities; this requires generalizing the proof that the Igusa polynomials are squarefree (Remark~\ref{R:Igusa polynomial}). A context where that difficulty does not arise is trying to prove that the jumping locus is nonempty; in this case any positive lower bound on the degree would suffice.
\end{remark}

\section{Moduli stacks and generic crystalline companions}
\label{sec:companions in fibration}

In this section, we construct a family of candidates for the generic fiber of a crystalline companion of an \'etale coefficient object on the total space of a curve fibration; this uses the uniformity estimates from \S\ref{sec:uniformities} in order to construct some ``small'' (i.e., noetherian) moduli stacks of relative mod-$p^n$ connections.
The relationship with the original \'etale coefficient object will be expressed in terms of the $p$-adic local systems associated to the succcessive quotients of the slope filtrations. One key complication is that we end up not with a single candidate for the generic fiber, but a family indexed by a profinite set; we will upgrade the conclusion (under a more restrictive hypothesis) in \S\ref{sec:companion points}.

\subsection{Setup}
\label{sec:companion setup}

We begin by describing the geometric setup in which we will work. 
In \S\ref{sec:companion points} we will impose some additional conditions; see  Hypothesis~\ref{H:companion points new}.

\begin{hypothesis} \label{H:companion points1}
Throughout \S\ref{sec:companions in fibration}:
\begin{itemize}
\item
Assume that $k$ is finite.
\item
Let $S$ be a smooth, affine
geometrically irreducible scheme of finite type over $k$.
\item
Let $\eta$ be the generic point of $S$.
\item
Let $f\colon \overline{X} \to S$ be a smooth curve fibration 
with pointed locus $Z$ and unpointed locus $X$.
\item
Let $s\colon S \to X$ be an additional section.
\item
Let $\eta_X$ be the generic point of $X$ (or equivalently $\overline{X}$).
\end{itemize}
\end{hypothesis}

\begin{defn} \label{D:density}
For any subset $W$ of $S^\circ$, the \emph{lower density} $\underline{\delta}(W)$ and the \emph{upper density} $\overline{\delta}(W)$ are given by 
\[
\underline{\delta}(W) = \liminf_{n \to \infty} \frac{\#W(k_n)}{\#S(k_n)},
\qquad
\overline{\delta}(W) = \limsup_{n \to \infty} \frac{\#W(k_n)}{\#S(k_n)},
\]
where $k_n$ denotes the degree-$n$ extension of $k$.
When $\underline{\delta}(W) = \overline{\delta}(W)$, we denote the common value by $\delta(W)$ and call it the \emph{density} of $W$.
Note that any nowhere dense closed subset of $S$ has density 0.

We may similarly define upper and lower densities for subsets of $\overline{X}^{\circ}$.
As in the proof of the Lang--Weil bound \cite{lang-weil}, we have
\[
\#f^{-1}(W)(k_n) = \#k_n \#W(k_n) + O((\#k_n)^{1/2}\#S(k_n))
\]
where the implied constant is uniform in $W$ and $n$; this shows that
\[
\underline{\delta}(f^{-1}(W)) = \underline{\delta}(W), \qquad
\overline{\delta}(f^{-1}(W)) = \overline{\delta}(W).
\]
\end{defn}

\begin{remark} \label{R:interpret connection}
For any perfect point $y$ mapping to $S$,
we may interpret $\overline{X} \times_S y \to y$ as a smooth curve fibration,
and so we may consider the category $\FIsoc^{\dagger}(X \times_S y)$
or the full subcategory of objects which are docile along $Z \times_S y$.
\end{remark}

\begin{hypothesis} \label{H:companion points}
Throughout \S\ref{sec:companions in fibration}, 
fix a category $\calC$ of \'etale coefficient objects on $X$ and a normalized $p$-adic valuation $v$ on the algebraic closure of $\QQ$ of the full coefficient field of $\calC$, and consider all Newton polygons (Lemma{-}Definition~\ref{L:generic Newton polygon}) and crystalline companions
with respect to $v$. Choose $\calE \in \calC$ such that:
\begin{itemize}
\item
$\calE$ is docile along $Z$ and absolutely irreducible of rank $r$;
\item
$\calE$ is $E$-algebraic for some number field $E$;
\item
$\wedge^r \calE$ is constant;
\item
the pullback of $\calE$ along $s$ is constant;
\item
the generic Newton polygon $N_{\eta_X}(\calE)$
(Lemma-Definition~\ref{L:generic Newton polygon}) has 
least slope $0$ occurring with some multiplicity $e$;
\item
we have $N_{s(x)}(\calE) = N_{\eta_X}(\calE)$ for all $x \in S$.
\end{itemize}
We also fix data as follows.
\begin{itemize}
\item
A field $L$ as in Lemma~\ref{L:bound coefficient extension} for $r$ as above and $L_0$ the $v$-adic completion of $E$.
\item
For $i=1,\dots,r$, an element $\lambda_i \in E^\times$ such that the constant twist of $\wedge^i \calE$ by $\lambda_i^{-1}$ has least generic Newton slope 0 occurring with some multiplicity $e_i$, taking $\lambda_1 = 1$.
\end{itemize}
\end{hypothesis}

\begin{lemma}  \label{L:etale fibration}
The sheaves
$R^j f_{\et*} \calE$ are \'etale coefficient objects on $S$ for $j \geq 0$, and vanish for $j > 2$;
moreover, the formation of these commutes with arbitrary base change on $S$. 
\end{lemma}
\begin{proof}
These statements reduce to the corresponding statements for a torsion \'etale sheaf $E$ on $X_{\overline{k}}$
of order prime to $p$. In this context, the sheaves $R^j f_{\et *}  E$ are constructible and vanish for $j>2$.
To check that they are lisse, using the proper base change theorem \stacktag{095S} it suffices to check that the cohomology groups
of $E$ on geometric fibers have locally constant dimension. For $j=0$, this follows from the tame specialization theorem
\cite[Expos\'e XIII, Corollaire~2.12]{grothendieck-sga1-tame}; this implies the case $j=2$ by Poincar\'e duality
\cite[Theorem~V.2.1]{milne}. Given these cases, the case $j=1$ follows from the local constancy of the Euler characteristic, which is implied by the
Grothendieck--Ogg--Shafarevich formula \cite{grothendieck-ogg-shafarevich}.
\end{proof}

We next restrict $\calE$ to the fibers of the morphism $f$ and apply what we know about the existence of crystalline companions on curves.

\begin{defn} \label{D:fiber restriction}
For $x \in S^\circ$, write $\calE_x$ for the restriction $\calE|_{X \times_S x}$.
Since $N_{s(x)}(\calE) = N_{\eta_X}(\calE)$, by Lemma-Definition~\ref{L:generic Newton polygon}(b) and (c), the generic Newton polygon of $\calE_x$ equals $N_{\eta_X}(\calE)$.
\end{defn}

The following is a variant of \cite[Lemma~3.2.1]{kedlaya-companions}.

\begin{lemma} \label{L:contracted section2}
The following statements hold.
\begin{enumerate}
\item[(a)]
The coefficient object $\calE_x$ is absolutely irreducible.
\item[(b)]
For $i=1,\dots,r$, $\wedge^i \calE$ is semisimple and absolutely semisimple, and each of its (absolutely) irreducible components restricts to an (absolutely) irreducible object on $X \times_S x$.
\item[(c)]
The object $\calE^\dual \otimes \calE$ is semisimple and absolutely semisimple, and each of its (absolutely) irreducible components restricts to an (absolutely) irreducible object on $X \times_S x$.
\end{enumerate}
\end{lemma}
\begin{proof}
By Lemma~\ref{L:etale fibration}, $f_{\et *} (\calE^\dual \otimes \calE)$ is an \'etale coefficient object.
Since it is also a subobject of the constant coefficient object $s_{\et}^* (\calE^\dual \otimes \calE)$, it must itself be constant. By Schur's lemma, this yields (a).

For (b), there are surjective homomorphisms $G(\calE) \to G(\wedge^i \calE)$ and $\overline{G}(\calE) \to \overline{G}(\wedge^i \calE)$, so we may apply Lemma~\ref{L:monodromy} to see that $\wedge^i \calE$ is semisimple and absolutely semisimple; we then apply the previous logic to each irreducible constituent
of $\wedge^i \calE$. The proof of (c) is similar.
\end{proof}

\begin{defn} \label{D:companions on fibers}
For $x \in S^\circ$, 
by Theorem~\ref{T:companion dimension 1} there exists $\calF_x \in \FIsoc^\dagger(X \times_S x) \otimes \overline{\QQ}_p$ which is a companion of $\calE_x$ with respect to the place $v$. By Lemma~\ref{L:companion irreducible}(b) and (d),
$\calF_x$ is absolutely irreducible and unique up to (noncanonical) isomorphism, and $\wedge^r \calF_x$ is constant.
By Corollary~\ref{C:companions are tame}, $\calF_x$ is docile along $Z \times_S x$. 
\end{defn}

We next set some notation regarding Newton polygons, keeping in mind the temporary discrepancy between Lemma{-}Definition~\ref{L:generic Newton polygon}(c) and Lemma~\ref{L:convergent Newton polygon}(c) (compare Corollary~\ref{C:open constant NP}).
\begin{defn} \label{D:companions on fibers twist}
For $x \in S^\circ$, let $V_x$ be the set of $y \in X \times_S x$ for which $N_y(\calF_x) = N_{\eta_X}(\calE)$
(or equivalently $N_y(\calE_x) = N_{\eta_X}(\calE)$).
By Hypothesis~\ref{H:companion points}, $s(x) \in V_x$ and so $V_x$ is nonempty. By Lemma{-}Definition~\ref{L:generic Newton polygon}(c),
$V_x$ is open in $X \times_S x$.

For $i=1,\dots,r$, let $e_i$ be the multiplicity of the least slope of the generic Newton polygon of $\wedge^i \calE$.
For $x \in S^\circ$, let $\tilde{\calF}_{x,i}$ be the constant twist of $\wedge^i \calF_x$ by $\lambda_i^{-1}$. By the previous paragraph, for every $y \in V_x$, $N_y(\tilde{\calF}_{x,i})$ has least slope $0$ with multiplicity $e_i$;
let $\rho_{x,i}\colon \pi_1(V_x) \to \GL_{e_i}(L)$ be the resulting unit-root representation  (Definition~\ref{D:unit-root subobject}).
\end{defn}

\begin{remark} \label{R:correspondence preserving NP}
For $i=1,\dots,r$, 
Lemma~\ref{L:companion irreducible}(a) produces a bijection between
the irreducible constituents of $\wedge^i \calF_x$ and $\wedge^i \calE_x$ preserving the generic Newton polygon.
By Lemma~\ref{L:contracted section2}, the latter are the restrictions of irreducible constituents of $\wedge^i \calE$.
\end{remark}

\subsection{Moduli stacks of $F$-crystals}
\label{subsec:crystals on log}

We now pick up the thread from \S\ref{subsec:moduli conn} to introduce some moduli stacks of $F$-crystals. The uniformity estimates from \S\ref{sec:uniformities} will allow these to be of finite type over $\ZZ$ while still being large enough to contain points reflecting the existence of the fiberwise crystalline companions $\calF_x$ of $\calE$.

The development of moduli stacks here is limited to what is needed for our application to constructing crystalline companions. A more robust development, in the language of formal and adic stacks, can be found in \cite{oh-shimizu}. In the other direction, many of the key ideas appear already in \cite{esnault-groechenig}.

\begin{hypothesis} \label{H:relative connections}
Throughout \S\ref{subsec:crystals on log}, 
fix a smooth affine scheme $\frakS$ over $W(k)$ lifting $S$ and a lift $\sigma_S \colon \frakS \to \frakS$ of the absolute Frobenius on $S$.
Assume the existence of, and fix the choice of, a smooth curve fibration $\frakf\colon \overline{\frakX} \to \frakS$ 
lifting $f$ with pointed locus $\frakZ$ and unpointed locus $\frakX$,
and a line bundle $\calL$ on $\overline{\frakX}$ which is very ample relative to $\frakf$.
\end{hypothesis}

\begin{defn}
For $n$ a positive integer and let
$\frakS_n, \overline{\frakX}_n, \frakX_n, \frakZ_n$ be the base extensions of $\frakS, \overline{\frakX}, \frakX, \frakZ$ along $\Spec \ZZ/p^n \ZZ \to \Spec \ZZ_p$.
We also allow $n = \infty$, corresponding to taking $\frakS_\infty = \frakS$ and so on.
\end{defn}

When $n=1$, we have schemes over $\FF_p$ and thus the notion of $p$-curvature (Definition~\ref{D:conn}).
In the general case, we impose this condition modulo $p$.

\begin{defn} \label{D:Frobenius pullback}
Define
\[
\Conn^{p}_{(\overline{\frakX}_n,\frakZ_n)/\frakS_n}  := 
\Conn_{(\overline{\frakX}_n,\frakZ_n)/\frakS_n} \times_{\Conn_{(\overline{\frakX}_1,\frakZ_1)/\frakS_1}} \Conn^{p}_{(\overline{\frakX}_1,\frakZ_1)/\frakS_1} 
\]
By Lemma~\ref{L:quasicompact connection stack2}, $\Conn^{p}_{(\overline{\frakX}_n,\frakZ_n)/\frakS_n}$ is qcqs and of finite type over $\frakS_n$.

Following \stacktag{07JH}, \cite[Theorem~7.3]{oh-shimizu}, we obtain a \emph{Verschiebung}
morphism 
\[
V\colon \Conn^{p}_{(\overline{\frakX}_n,\frakZ_n)/\frakS_n} \to \Conn^{p}_{(\overline{\frakX}_n,\frakZ_n)/\frakS_n}
\]
lying over $\sigma_S\colon \frakS_n \to \frakS_n$, whose effect can be interpreted locally on $\overline{\frakX}_n$ as pullback along a local Frobenius lift.
(For $p>2$ this morphism
extends to $\Conn_{(\overline{\frakX}_n,\frakZ_n)/\frakS_n}$;
see for example \cite[Proposition~8.1]{esnault}.)
\end{defn}

\begin{defn}
Let us temporarily write
\[
C := \Conn^p_{(\overline{\frakX}_n,\frakZ_n)/\frakS_n} \times \Conn^p_{(\overline{\frakX}_n,\frakZ_n)/\frakS_n}.
\]
Let $\pi_1, \pi_2\colon C \to \Conn^{p}_{(\overline{\frakX}_n,\frakZ_n)/\frakS_n}$ be the first and second projection maps.
Let $H$ be the stack over 
$C$ parametrizing morphisms from the first connection to the second;
in other words,
\[
H := C \times_{\Conn_{(\overline{\frakX}_n,\frakZ_n)/\frakS_n} \times \Conn_{(\overline{\frakX}_n,\frakZ_n)/\frakS_n}} \ConnHom_{(\overline{\frakX}_n,\frakZ_n)/\frakS_n}.
\]
Let $\Gamma_V^T$ be the stack over $C$ parametrizing isomorphisms of the first connection with the $V$-image of the second (i.e., the transpose of the graph of $V$).
Then set
\[
\FConn_{(\overline{\frakX}_n,\frakZ_n)/\frakS_n} := H \times_{C} \Gamma_F^T
\]
and
\[
\FConn^{P,\calL,m}_{(\overline{\frakX}_n,\frakZ_n)/\frakS_n} := \FConn_{(\overline{\frakX}_n,\frakZ_n)/\frakS_n} \times_{\pi_1, \Conn_{(\overline{\frakX}_n,\frakZ_n)/\frakS_n}} \Conn^{P,\calL,m}_{(\overline{\frakX}_n,\frakZ_n)/\frakS_n}.
\]
\end{defn}

\begin{lemma} \label{L:noetherian stack}
For any positive integer $n$, the stack $\FConn^{P, \calL, m}_{(\overline{\frakX}_n, \frakZ_n)/\frakS_n} \otimes_{\ZZ_p} \frako_L$ is qcqs and of finite type over $\ZZ$. In particular, it is noetherian.
\end{lemma}
\begin{proof}
By two applications of Proposition~\ref{P:conn hom} (plus the fact that $\Spec \frakS_n \to \Spec \ZZ$ is of finite type), we reduce to the fact that
$\Conn^{p}_{(\overline{\frakX}_n,\frakZ_n)/\frakS_n}$ is qcqs and of finite type over $\frakS_n$
(Definition~\ref{D:Frobenius pullback}).
\end{proof}

We can now translate the uniformity statements from \S\ref{sec:uniformities} into a geometric uniformity statement about the points in moduli associated to crystalline lattices.
\begin{lemma} \label{L:uniformity for crystalline lattices}
There exist a polynomial $P \in \QQ[t]$ and a positive integer $m$ with the following property:
for every $x \in S^\circ$, the object $\calF_x \in \FIsoc^\dagger(X \times_S x) \otimes L$
from Definition~\ref{D:companions on fibers}
admits a crystalline lattice represented by a coherent system of morphisms
\[
W_n(x) \to \FConn^{P, \calL, m}_{(\overline{\frakX}_n, \frakZ_n)/\frakS_n} \otimes_{\ZZ_p} \frako_L
\]
lying over the Frobenius-equivariant maps $W_n(x) \to \frakS_n$.
\end{lemma}
\begin{proof}
This is a direct consequence of Proposition~\ref{P:isocrystal reduction slopes},
using Lemma~\ref{L:HN to sections} to translate the bound on HN slopes into a choice of $P$ and $m$, and pulling back once along Frobenius to enforce zero $p$-curvature.
\end{proof}

\subsection{Companion points and their deformations}
\label{subsec:universal unit root}

We next study the moduli points arising from fiberwise crystalline companions in more detail.
This study crucially includes a weak integral adaptation of cohomological rigidity
(Proposition~\ref{P:cohomological rigidity}).

\begin{hypothesis} \label{H:uniform bounds}
Throughout \S\ref{subsec:universal unit root},
retain Hypothesis~\ref{H:relative connections};
fix a choice of $P,m$ as in Lemma~\ref{L:uniformity for crystalline lattices};
and let $n$ be an arbitrary positive integer.
\end{hypothesis}

\begin{lemma} \label{L:uniform H2 torsion}
There exists a nonnegative integer $n_0$ with the following property:
for every perfect field $\kappa$ of characteristic $p$ and every morphism $x \cong \Spec W(\kappa) \to \Conn_{(\overline{\frakX},\frakZ)/\frakS}^{P,\calL,m}$, for $\frakF_x$ the corresponding vector bundle with logarithmic connection on $(\overline{\frakX},\frakZ) \times_{\frakS} W(\kappa)$, the torsion subgroups of the de Rham cohomology groups $H^i(\frakF_x)$ for $i=1,2$ are killed by $p^{n_0}$.
\end{lemma}
\begin{proof}
Since $\Conn_{(\overline{\frakX},\frakZ)/\frakS}^{P,\calL,m}$ is quasicompact by Lemma~\ref{L:quasicompact connection stack},
we can choose a surjective morphism $U \to \Conn_{(\overline{\frakX},\frakZ)/\frakS}^{P,\calL,m}$ with $U$ a noetherian scheme. Let $\frakF$ be the corresponding vector bundle with logarithmic $U$-linear connection on
$(\overline{\frakX},\frakZ) \times_{\frakS} U$. Using the ``generic flatness stratification'' lemma from \stacktag{0ASY}, we can then construct a locally closed stratification of $U$ with the property that for any stratum $V \subset U$, for any morphism $g\colon W \to V$ of noetherian schemes,
we have $g^* H^i(\frakF|_V) \cong H^i(\frakF|_W)$ for all $i$. 
Consequently, for any morphisms $x_j \cong \Spec W(\kappa_j) \to U$ for $j=1,2$ which factor through the same
irreducible stratum $V$, the smallest value of $n_0$ which works for $x_1$ 
(which exists because the groups $H^i(\frakF_{x_1})$ are finitely generated $W(\kappa_1)$-modules)
also works for $x_2$;
the claim is now evident.
\end{proof}

\begin{defn} \label{D:truncation}
For each $x \in S^\circ$, we obtain from Lemma~\ref{L:uniformity for crystalline lattices}
a $W_n(x)$-valued point of the stack $\FConn^{P, \calL, m}_{(\overline{\frakX}_n, \frakZ_n)/\frakS_n} \otimes_{\ZZ_p} \frako_L$. Let $M_n$ be the scheme-theoretic image (Definition~\ref{D:scheme-theoretic image}) of the disjoint union of these maps;
we refer to the resulting point of finite type $W_n(x) \to M_n$ as a \emph{companion point} of $M_n$ lying over $x$.
\end{defn}

\begin{remark} \label{R:closed conditions}
We may read off various properties of the $M_n$ from the construction.
\begin{itemize}
\item
We have a natural projection $M_{n+1} \to M_n$ lying over $\frakS_{n+1} \to \frakS_n$.
\item
On some open dense subset of $M_n$ through which all of the companion points factor,
 the Frobenius structure on the universal connection is injective with cokernel killed by $\lambda_r$.
\item
The $r$-th exterior power of the universal connection on $M_n$ is locally trivial.
In particular, we get a covering of $M_n$ by imposing a fixed global trivialization; see 
Proposition~\ref{P:carryover data}(d) for a related point.
\end{itemize}
\end{remark}

\begin{prop} \label{P:Crys stack applied1}
The stack $M_n$ is qcqs and of finite type over $\frakS_n$ via $\pi_1$.
In particular, it is noetherian.
\end{prop}
\begin{proof}
As $M_n$ is defined as a closed substack of $\FConn^{P, \calL, m}_{(\overline{\frakX}_n, \frakZ_n)/\frakS_n} \otimes_{\ZZ_p} \frako_L$,
this follows directly from Lemma~\ref{L:noetherian stack}.
\end{proof}

\begin{defn} \label{D:only dominant components}
By Proposition~\ref{P:Crys stack applied1}, $M_n$ is a noetherian stack and hence has finitely many irreducible components. Here, as in \stacktag{0DR4}, we are defining irreducible components at the level of underlying topological spaces and then taking the corresponding reduced closed substacks.

For each irreducible component $T_n$ of $M_n$, let $W(T_n)$ be the set of $x \in S^\circ$ for which $T_n$ contains a companion point lying over $x$. From Definition~\ref{D:truncation}, $\bigcup_{T_n \in M_n} W(T_n) = S^\circ$.

Let $M'_n$ be the union of the irreducible components $T_n$ of $M_n$
for which $\overline{\delta}(W(T_n)) > 0$;
then  $\delta\left(\bigcup_{T_n \in M'_n} W(T_n) \right) = 1$.
Consequently, $M'_n \neq \emptyset$; the projection $M_{n+1} \to M_n$ induces a dominant map
$M'_{n+1} \to M'_n$; and the maps $M'_n \to \frakS_n$ are also dominant.
\end{defn}

\begin{lemma} \label{L:infinitesimal deformation}
For each $n$ and each $\epsilon \in (0,1)$, 
there exists an integer $n' \geq n$
and a subset $W$ of $S^\circ$ with $\underline{\delta}(W) \geq 1-\epsilon$
with the following property. 
Let $W_{n'}(x) \to M'_{n'}$ be a companion point lying over some $x \in W$.
Let $Y$ be the connected component of $M'_{n'} \times_{\frakS_{n'}} W_{n'}(x)$
through which this companion point factors. Then the 
map $Y_{\red} \to M'_n$ is ``constant'' in the sense that it factors through the companion point $W_n(x) \to M'_n$ lying under $W_{n'}(x) \to M'_{n'}$.
\end{lemma}
\begin{proof}
Fix a value $n_0$ as in Lemma~\ref{L:uniform H2 torsion} and a positive integer $c$ greater than or equal to the slope of $\wedge^r \calE$.
Let $\tilde{\calE}^\dual$ be the constant twist of $\calE^\dual$ by $p^c$. Similarly, let $\tilde{\calF}_x^\dual$ be the constant twist of $\calF_x^\dual$ by $p^c$, and let $\tilde{\frakF}_x^\dual$ be the crystalline lattice in $\tilde{\calF}_x^\dual$ induced by $\frakF_x^\dual$
(Remark~\ref{R:constant twist of dual lattice}).
Finally, let $\frakH_x$ be the quotient of $\tilde{\frakF}_x^\dual \otimes \frakF_x$ by the trace component.
For each integer $i \geq 0$ and each pair of integers $m' \geq m > 0$, we have a commutative diagram
\begin{equation} \label{eq:torsion in de Rham}
\xymatrix{
0 \ar[r] & H^i(\frakH_x)/(p^{m'}) \ar[r] \ar[d] & H^i(\frakH_{x,m'}) \ar[r] \ar[d] & H^{i+1}(\frakH_x)[p^{m'}] \ar[r] \ar^{\times p^{m'-m}}[d] \ar[r] & 0 \\
0 \ar[r] & H^i(\frakH_x)/(p^{m}) \ar[r] & H^i(\frakH_{x,m}) \ar[r] & H^{i+1}(\frakH_x)[p^{m}] \ar[r] & 0
}
\end{equation}
with exact rows.

By Lemma~\ref{L:etale fibration}, $R^1 f_* (\tilde{\calE}^\dual \otimes \calE)$ is an \'etale coefficient object on $S$, which by  Hypothesis~\ref{H:companion points1} admits a crystalline companion $\calH \in \FIsoc^\dagger(S) \otimes L'$ for some finite extension $L'$ of $\QQ_p$ (which need not equal $L$). 
By Corollary~\ref{C:slope filtration}, $\calH$ admits a slope filtration in $\FIsoc(S_0) \otimes L'$ for some open dense subscheme $S_0$ of $S$.
After twisting (and possibly enlarging $L'$), we may apply Proposition~\ref{P:unit-root rep}(a) to each successive quotient of the filtration; this yields a sequence of continuous representations
$\{\rho_i\colon \pi_1(S_0) \to \GL_{m_i}(\frako_{L'})\}_{i=1}^l$.

By Proposition~\ref{P:cohomological rigidity}, 1 never occurs as a Frobenius eigenvalue on $\calH$. From the continuity of the $\rho_i$ and Chebotaryov density, we obtain a positive integer $n_1$ and a subset $W$ of $W_0 \cap S_0^\circ$ with $\underline{\delta}(W) \geq 1-\epsilon$ such that the Frobenius characteristic polynomial 
of $(R^1 f_* (\tilde{\calE}^\dual \otimes \calE))_x$ evaluates at $p^c$ to a value not divisible by $p^{n_1}$.
We will prove the claim for this choice of $W$ with $n' := n + 3n_0+ n_1$; we thus assume hereafter that $x \in W$.

Since $M'_n \to S$ is of finite type, for $Y_{\red} \to M'_n$ not to be constant there would have to exist a map $C \to Y$ through which $W_{n'}(x) \to M'_{n'}$ factors such that $C$ is a smooth curve over $W_{n'}(x)$. To rule this out, it will suffice to produce an integer $n'' \geq 0$ for which any infinitesimal deformation of $\frakF_{x,n'}$ with trivialized determinant becomes trivial upon multiplication by $p^{n''}$ followed by restriction to $\frakF_{x,n+n''}$. We prove this with $n'':= 2n_0$.

The infinitesimal deformations of $\frakF_{x,n'}$ with trivialized determinant form a group which by Hochschild--Serre sits in the middle of an exact sequence
\begin{equation} \label{eq:integral Hochschild-Serre}
H^0(\frakH_{x,n'})_{\varphi-p^c} \to * \to H^1(\frakH_{x,n'})^{\varphi-p^c}.
\end{equation}
Since $\calE$ is absolutely irreducible and $H^0(\frakH_x)$ is torsion-free, 
we have $H^0(\frakH_x) = 0$.
We may thus apply Lemma~\ref{L:uniform H2 torsion} and \eqref{eq:torsion in de Rham} to see that $H^0(\frakH_{x,n'})$ is killed by
$p^{n_0}$, as then is $H^0(\frakH_{x,n'})_{\varphi-p^c}$.

Meanwhile, by Lemma~\ref{L:uniform H2 torsion}, the torsion subgroup of $H^2(\frakH_x)$ is killed by $p^{n_0}$. By \eqref{eq:torsion in de Rham}, any element of $H^1(\frakH_{x,n'})^{\varphi-p^c}$ 
projects to an element of $H^1(\frakH_{x,n'-n_0})^{\varphi-p^c}$ which lifts to $(H^1(\frakH_x)/(p^{n'-n_0}))^{\varphi-p^c}$. By Lemma~\ref{L:uniform H2 torsion} again, multiplication by $p^{n_0}$ maps $H^1(\frakH_x)$ onto a torsion-free subgroup $N$ which is a complement of the torsion subgroup;
in particular, this complement is stable under $\varphi - p^c$. We can thus decompose our given element of $(H^1(\frakH_x)/(p^{n'-n_0}))^{\varphi-p^c}$ uniquely as the sum of an element of 
the group $H^1(\frakH_x)[p^n]^{\varphi-p^c}$
plus an element of the group $(N/(p^{n'-n_0}))^{\varphi-p^c}$.
The former element is killed by $p^{n_0}$;
by Cayley--Hamilton, the latter element is killed by $p^{n_1}$ and hence 
projects to zero in $N/(p^{n'-n_0-n_1})$.

Combining the previous discussion with \eqref{eq:integral Hochschild-Serre} and a quick diagram chase, we see that 
multiplying the original deformation class by $p^{2n_0} =  p^{n''}$
yields a class that reduces to zero modulo $p^{n'-n_0-n_1} = p^{n+n''}$.
\end{proof}

\begin{cor} \label{cor:generically finite moduli}
Fix an integer $n'$ as in Lemma~\ref{L:infinitesimal deformation}.
For each positive integer $n$, let $M''_n$ be the scheme-theoretic image (see Definition~\ref{D:scheme-theoretic image}) of the map $M'_{n+n'} \to M'_{n}$. Then the map $M''_n \to \frakS_n$ is generically finite.
\end{cor}
\begin{proof}
We may check the claim after pulling back along a dominant \'etale morphism to $S$.
By forming a Stein factorization \stacktag{0A1B}, after such a pullback we can ensure that every component of $M''_n$ maps to $\frakS_n$ with  connected geometric fibers of constant dimension.
By Lemma~\ref{L:infinitesimal deformation} (applied with any $\epsilon > 0$), this
constant dimension must be at most $0$.
\end{proof}

\subsection{Candidates for the generic fiber}
\label{subsec:companions multiplicity 1}

We now deliver the final result of our stack-theoretic analysis: 
a family of candidates for the generic fiber of the crystalline companion.
The relationship with the original object $\calE$ will be expressed in terms of
the unit-root representations of 
the successive quotients of the slope filtration (Corollary~\ref{C:slope filtration}).

\begin{prop} \label{P:carryover data}
There exists a collection of data as follows
(where $n$ runs over all positive integers and $i$ runs over $\{1,\dots,r\}$):
\begin{itemize}
\item
a family of dominant \'etale morphisms $\cdots \to T_2 \to T_1 \to S$
in which each $T_n$ is irreducible with generic point $\eta_n$;
\item
for each $n$, an open dense subset $V_n$ of $X \times_S T_n$;
\item
a surjective sequential inverse system of finite sets $\varprojlim_n J_n$ with limit $J$;
\item
for each $n$, each $i$, and each $j_n \in J_n$,
a continuous representation $\rho_{n,i,j_n}\colon \pi_1(V_n) \to \GL_{e_i}(\frako_L/p^n \frako_L)$;
\item
a perfect point $\tilde{\eta}$ lying over $\varprojlim_n \eta_n$;
\item
for each $j \in J$, 
an object $\calF^{(j)}_{\tilde{\eta}} \in \FIsoc^\dagger(X \times_S \tilde{\eta}) \otimes L$;
\end{itemize}
subject to the following compatibilities.
\begin{itemize}
\item[(a)]
For each $j \in J$, $\calF^{(j)}_{\tilde{\eta}}$ is docile along $Z \times_S \tilde{\eta}$, absolutely irreducible,  and its generic Newton polygon is $N_{\eta_X}(\calE)$.
\item[(b)]
Let $\tilde{\eta}_X$ be the generic point of $X \times_S \tilde{\eta}$.
For each $i$, for $j \in J$,
let $\tilde{\calF}^{(j)}_{\tilde{\eta},i}$ be the constant twist of $\wedge^i \calF^{(j)}_{\tilde{\eta}}$ by $\lambda_i^{-1}$.
Define the representation $\rho^{(j)}_i \colon \pi_1(\tilde{\eta}_X) \to \GL_{e_i}(L)$ by
forming  the constant twist of $\wedge^i \calF^{(j)}_{\tilde{\eta}}$ by $\lambda_i^{-1}$
and restricting the unit-root representation (Definition~\ref{D:unit-root subobject}). Then we can choose a stable lattice for $\rho^{(j)}_i$ in such a way that for each $n$, the mod-$p^n$ reduction of $\rho^{(j)}_i$ factors through $\rho_{n,i,j_n}$.

\item[(c)]
For each $n$, for each $i$, 
for each $x$ outside of some subset of $S^\circ$ of density $0$,
we can find $j_n \in J_n$ \emph{depending on $x$} such that
 for all $y \in (V_n \times_X V_x)^\circ$,
the characteristic polynomial of $\rho_{x,i}(\Frob_y)$ reduces modulo $p^n$ to the characteristic polynomial of $\rho_{n,i,j_n}(\Frob_y)$.

\item[(d)]
Let $I$ be any subset of $\{1,\dots,r\}$.
Suppose that for some open dense subset $W$ of $X$, 
for each $i \in I$, there exists a continuous representation $\rho_i\colon \pi_1(W) \to \GL_{e_i}(L)$ such that for each $x \in S^\circ$ with $W \times_X V_x \neq \emptyset$, the restrictions of $\rho_i$ and $\rho_{x,i}$ to $\pi_1(W \times_X V_x)$ are isomorphic. Then for each $i \in I$, for all $j \in J$, $\rho_i^{(j)}$ is isomorphic to the restriction of $\rho_i$.

\end{itemize}
\end{prop}
\begin{proof}
By invoking Proposition~\ref{P:stable reduction} and then shrinking $S$ as needed, we may put ourselves in a situation where Hypothesis~\ref{H:relative connections} applies.
Let $J'_n$ be the set of irreducible components of the stack $M''_n$ (Corollary~\ref{cor:generically finite moduli}) and set $J := \varprojlim_n J'_n$.
We can then choose a family of dominant \'etale morphisms $\cdots \to T'_2 \to T'_1 \to S$ so that $T'_n$ dominates each component of $M''_n$.
Let $\eta'_n$ be the generic point of $T'_n$ and let $\tilde{\eta}$ be a perfect point lying over
$\varprojlim_n \eta_n$. Set $J := \varprojlim_n J'_n$;
for $j = (j'_n) \in J$, use the compositions 
\begin{equation} \label{eq:truncation of limit crystal}
\tilde{\eta} \to \eta'_n \to T'_n \to j'_n \to M''_n \to 
M_n \to \FConn^{P, \calL, m}_{(\overline{\frakX}_n, \frakZ_n)/\frakS_n} \otimes_{\ZZ_p} \frako_L
\end{equation}
to construct $\calF^{(j)}_{\tilde{\eta}}$ together with a crystalline lattice.
This does not by itself guarantee that the Frobenius structure becomes an isomorphism after inverting $p$; for this, we use Remark~\ref{R:constant twist of dual lattice} (and specifically the uniformity of the constant $c$, which depends only on $N_{\eta_X}(\calE)$) to construct a Frobenius structure on the dual connection.

At this point we check (a) as follows.

\begin{itemize}
\item
To check that $\calF^{(j)}_{\tilde{\eta}}$ is absolutely irreducible,
let $\calH$ be the quotient of $(\calF^{(j)}_{\tilde{\eta}})^\dual \otimes \calF^{(j)}_{\tilde{\eta}}$
by the trace component; let $\frakH$ be a crystalline lattice in $\calH$; and let $\frakH_n$ be the mod-$p^n$ reduction of $\frakH$. It will suffice to check that $H^0(\frakH) = 0$, or equivalently that $H^0(\frakH_n) \to H^0(\frakH_1)$ is the zero map for some $n$; this may be deduced from Lemma~\ref{L:uniform H2 torsion} and 
Lemma~\ref{L:infinitesimal deformation}.

\item
To check that the generic Newton polygon of $\calF^{(j)}_{\tilde{\eta}}$  is $N_{\eta_X}(\calE)$, 
it suffices to check that for $i=1,\dots,r$, 
the least generic slopes of $\wedge^i \calF^{(j)}_{\tilde{\eta}}$ and its dual are bounded below by the least generic slopes of $\wedge^i \calE$ and its dual, respectively.
For this, we apply the relationship between Hodge and Newton polygons in an $F$-crystal
\cite[Basic Slope Estimate 1.3.4]{katz-slope}.
\end{itemize}

We can now define $\tilde{F}_{\tilde{\eta},i}^j$ and $\rho_i^{(j)}$ as in (b).
Using  Remark~\ref{R:constant twist of dual lattice} again,
we see that for some strictly monotone function $f\colon \{1,2,\dots\} \to \{1,2,\dots\}$,
we can construct stable lattices for the $\rho_i^{(j)}$ in such a way that the mod-$p^n$ reduction of $\rho_i^{(j)}$ depends only on the image $j_n$ of $j$ in $J_n := J'_{f(n)}$. 
We now set $T_n := T'_{f(n)}$ and $\eta_n := \eta'_{f(n)}$.
Using \eqref{eq:truncation of limit crystal} with $n$ replaced by $f(n)$,
we obtain $V_n$ and $\rho_{n,i,j_n}$ satisfying (b).

We obtain the remaining compatibilities as follows.

\begin{itemize}

\item
To check (c), we apply Lemma~\ref{L:infinitesimal deformation} repeatedly with $n$ fixed and $\epsilon$ tending to 0; the intersection of the exceptional sets then has density 0 (but see Remark~\ref{R:density 0 subset}).

\item
To check (d), we argue as in Remark~\ref{R:closed conditions}:
we replace $M_n$ with a covering by fixing maps to the universal connection from the appropriate truncations of the unit-root connections corresponding to the $\rho_i$ for $i \in I$ via Proposition~\ref{P:unit-root rep}(a).
 \qedhere
\end{itemize}
\end{proof}

\begin{remark} \label{R:density 0 subset}
In Proposition~\ref{P:carryover data}(c), the exceptional set has density 0 but can still be ``large'' from other points of view. For instance, it may be Zariski dense.

More seriously, as $n$ varies, there is nothing preventing the union of the exceptional sets from filling up $S^\circ$ (as density is only finitely additive). We thus cannot directly infer any comparison between the representations $\rho_i$ and $\rho_{i,x}$ for any particular $x \in S^\circ$.

In some sense these are both symptoms of the fact that the proof of Proposition~\ref{P:carryover data}
uses only data from the existence of ``vertical'' crystalline companions (on fibers over $S$), whereas we also have information available from ``horizontal'' crystalline companions. We will incorporate that information in \S\ref{sec:companion points}.
\end{remark}

\section{Crystalline companions on curve fibrations}
\label{sec:companion points}

With some key stack-theoretic input in hand (Proposition~\ref{P:carryover data}),
we now proceed to construct the crystalline companion of a suitable \'etale coefficient on the total space of a smooth curve fibration over a smooth base scheme, under some technical restrictions which will be lifted later.
See Hypothesis~\ref{H:companion points new} for running notation that will be in effect thereafter.

\subsection{Rank-\texorpdfstring{$1$}{1} objects from curves}
\label{subsec:rank 1 objects}

In order to improve upon Proposition~\ref{P:carryover data},
we need to construct some convergent unit-root isocrystals out of their putative restrictions to curves;
using Proposition~\ref{P:unit-root rep}(a), this translates into constructing $p$-adic representations of \'etale fundamental groups out of their putative restrictions to curves.
We describe here an adaptation of the technique of Wiesend \cite{wiesend1} (compare  \cite[Lemme~7.3]{cadoret}) without any tameness restriction, but limited in practice to the abelian case: the matching condition is at the level of elements of the target group, not conjugacy classes
(see Remark~\ref{R:rank 1 limitation} for further discussion).
One minor adaptation required in the $p$-adic setting is that we are only building a representation of the fundamental group of a dense open subset which we cannot identify \emph{a priori}; we must thus introduce a bit of extra bookkeeping in order to give the argument enough power to pin down this subset.

\begin{hypothesis}
Throughout \S\ref{subsec:rank 1 objects},
assume that $k$ is finite and $X$ is irreducible.
We do not need the running hypotheses of \S\ref{sec:companions in fibration} here.
\end{hypothesis}

\begin{lemma} \label{L:rank-1 assembly finite level}
Let $G$ be a finite group written multiplicatively.
Let $f\colon X^\circ \to G \cup \{0\}$ be a function (not identically zero) satisfying the following conditions.
\begin{enumerate}
\item[(i)]
For every curve $C$ in $X$, there is a (possibly empty) open subset $C_0$ of $C$ such that for $x \in C^\circ$, $f(x) \neq 0$ if and only if $x \in C_0^\circ$.
\item[(ii)]
With notation as in (i), if $C_0 \neq \emptyset$, then there exists a homomorphism $\chi_{C}\colon \pi_1^{\ab}(C_0) \to G$ such that $f(x) = \chi_{C}(\Frob_x)$.
\item[(iii)]
There exists a dominant \'etale morphism $Y \to X$ such that for every $x \in X^\circ$ in the image of $Y$,
either $f(x) = 0$ or there exists $y \in Y \times_X x$
with $f(x)^{\deg(y \to x)} = 1$.
\end{enumerate}
Then there exist an open dense subspace $U$ of $X$
and a character $\chi_n\colon \pi_1^{\ab}(U) \to G$ such that for each $x \in U^\circ$,
either $f(x) = 0$ or $f(x) = \chi(\Frob_x)$.
\end{lemma}
\begin{proof}
We may assume that $Y$ is connected and finite \'etale over some open dense subspace $U$ of $X$.
We may further assume that $\pi_1(Y)$ is normal in $\pi_1(U)$; then $\pi_1(U)/\pi_1(Y)$ is finite, so the set 
\[
H := \Hom(\pi_1(U)/\pi_1(Y), G)
\]
is also finite. Moreover, the statement of (iii) is now formally stronger: if $f_n(x) \neq 0$, then \emph{every} $y \in Y \times_X x$
satisfies $f_n(x)^{\deg(y \to x)} = 1$.

Suppose that $C$ is a curve in $X$ such that $C_0 \neq \emptyset$ and
\[
\pi_1(C \times_X U) \to \pi_1(U) \to \pi_1(U)/\pi_1(Y)
\]
is surjective.
The preimage of $\pi_1(Y)$ in $\pi_1(C \times_X U)$
may be identified with $\pi_1(C')$ where $C'$ is a connected component of $C \times_X Y$; by (ii) and (iii), for each $y' \in C' \times_C C_0$ we have $\chi_{C}(\Frob_{y'}) = 1$.
By Chebotaryov density, the kernel of $\chi_{C}|_{C \times_X U}$ contains $\pi_1(C')$, from which it follows that $\chi_{C}|_{C \times_X X} \in H$.

Let $W$ be an arbitrary nonempty finite subset of $X^\circ$
such that $f(x) \neq 0$ for all $x \in W$.
By \cite[Theorem~2.15]{drinfeld-deligne}, we can find a curve $C$ in $X$ such that $W \subseteq C$ (so $C_0 \neq \emptyset$)
and $\pi_1(C \times_X U) \to \pi_1(U) \to \pi_1(U)/\pi_1(Y)$ is surjective.
By the previous paragraph, there exists $\chi \in H$ for which $f(x) = \chi(\Frob_x)$ for all $x \in W \cap U$. By a compactness argument, we can choose a single $\chi \in H$ that works for all $W$,
which implies the desired result.
\end{proof}

\begin{remark} \label{R:no density 1 version}
One cannot strengthen Lemma~\ref{L:rank-1 assembly finite level} so that condition (iii) only holds for $x$ in some subset of $X^\circ$ of density 1 in the sense of Definition~\ref{D:density}:
the complement of this subset can be big enough to meet \emph{every} curve in $X$ in a cofinite set.
We will work around this point in the proof of Proposition~\ref{P:putative unit root}.
\end{remark}

\begin{cor} \label{C:rank-1 assembly infinite level}
Let $L$ be a finite extension of $\QQ_p$ and let $\varpi$ be a uniformizer of $L$.
Let $f\colon X^\circ \to \frako_L^\times \cup \{0\}$ be a function
(not identically zero) satisfying the following conditions.
\begin{enumerate}
\item[(i)]
For every curve $C$ in $X$, there is a (possibly empty) open subset $C_0$ of $C$ such that for $x \in C^\circ$, $f(x) \neq 0$ if and only if $x \in C_0^\circ$.
\item[(ii)]
With notation as in (i), if $C_0 \neq \emptyset$, then there exists a character $\chi_C\colon \pi_1^{\ab}(C_0) \to \frako_L^\times$ such that $f(x) = \chi_C(\Frob_x)$. Moreover, $\chi_C$ is ramified at every point of $C \setminus C_0$.
\item[(iii)]
For each $n>0$, there exists a dominant \'etale morphism $Y_n \to X$ such that for every $x \in X^\circ$ in the image of $Y_n$,
either $f(x) = 0$ or there exists $y \in Y_n \times_X x$
with $f(x)^{\deg(y \to x)} \equiv 1 \pmod{\varpi^n}$. 
\end{enumerate}
Then there exist a decreasing sequence of open dense subsets $X_n$ of $X$ and a sequence of characters $\chi_n\colon \pi_1^{\ab}(X_n) \to (\frako_L/\varpi^n \frako_L)^\times$ for which the following statements hold.
\begin{enumerate}
\item[(a)]
For each $x \in X^\circ$, $f(x) = 0$ if and only if $x \notin X_n$ for some $n$.
\item[(b)]
For each $x \in \bigcap_n X_n^\circ$, 
for each $n$, $f(x) \equiv \chi_n(\Frob_x) \pmod{\varpi^n}$.
\item[(c)]
For each $n> 0$, $X \setminus X_{n}$ has pure codimension $1$ in $X$.
\end{enumerate}
\end{cor}
\begin{proof}
For each $n>0$, we may apply Lemma~\ref{L:rank-1 assembly finite level} to produce an open dense subset $X_{n,0}$ of $X$ and a character $\chi_n\colon \pi_1^{\ab}(X_{n,0}) \to (\frako_L/\varpi^n \frako_L)^\times$ such that for each $x \in X_{n,0}^\circ$,
either $f(x) = 0$ or $f(x) \equiv \chi_n(\Frob_x) \pmod{\varpi^n}$.
For each curve $C$ in $X$ for which $C_0 \neq \emptyset$, 
combining the previous assertion with (ii) and Chebotaryov density yields
\begin{equation} \label{eq:compare restriction to curve1}
\chi_{n}|_{C_0 \times_X X_{n,0}} = \chi_{C,n}|_{C_0 \times_X X_{n,0}}.
\end{equation}
Let $X_n$ be the maximal open subset of $X$ to which $\chi_n$ extends; this satisfies (c) by Zariski--Nagata purity. From \eqref{eq:compare restriction to curve1}, we also have
\begin{equation} \label{eq:compare restriction to curve}
\chi_{n}|_{C_0 \times_X X_{n}} = \chi_{C,n}|_{C_0 \times_X X_{n}}.
\end{equation}

By \eqref{eq:compare restriction to curve},
$\chi_n|_{C_0 \times_X X_{n}}$ extends across $C \times_X X_n$.
In particular,  if $x \in X_n^\circ$ for all $n$, then $\chi_C$ is unramified at $x$, and (i) and (ii) together imply that $f(x) \neq 0$. This yields the ``only if'' implication of (a).

For $x \in \bigcap_n X_n^\circ$, for each $n$ we may choose a point $y \in X_{n,0}^\circ$ and a curve $C$ in $X$ containing $x$ and $y$, then apply (ii) and \eqref{eq:compare restriction to curve} to deduce (b).

By \eqref{eq:compare restriction to curve} and Chebotaryov again,
\[
\chi_{n+1}|_{X_n \times_X X_{n+1}} \equiv \chi_n|_{X_n \times_X X_{n+1}} \pmod{\varpi^n}.
\]
In particular, $\chi_n$ extends over $X_{n+1}$, so $X_{n+1} \subseteq X_n$.

It only remains to check the ``if'' implication of (a). 
Suppose the contrary; then there exist some $n > 0$ and some $x \in (X \setminus X_n)^\circ$ for which $f(x) \neq 0$.
For any irreducible component $Z$ of $X \setminus X_{n}$ containing $x$, the set $U$ of $y \in Z^\circ$ for which $f(y) \neq 0$ is Zariski dense in $Z$: for any closed strict subset $W$ of $Z$, we may choose a curve $C$ in $Z$ containing $x$ and some point $z \in (Z \setminus W)^\circ$ and then apply (i) to see that $f(y) \neq 0$ for some $y \in (Z \setminus W)^\circ$. For each $y \in U^\circ$,
choose a point $z \in X_{n}^\circ$
and a curve $C$ in $X$ containing $\{x,y,z\}$;
the inclusion of $x$ ensures that $C_0$ is a nonempty open subset of $C$, while 
the inclusion of $z$ ensures that $C \times_X X_n \neq \emptyset$,
so $C_0 \times_X X_n \neq \emptyset$.
By \eqref{eq:compare restriction to curve}, $\chi_{n}|_{C_0 \times_X X_{n}}$ extends across $y$; since $y$ runs over a Zariski dense subset of $Z$, this is enough to deduce that $\chi_{n}$ is unramified along $Z$, which contradicts our choice of $X_{n}$.
\end{proof}

\begin{cor} \label{C:rank-1 assembly infinite level stable}
With notation as in Corollary~\ref{C:rank-1 assembly infinite level}, suppose also that the sequence $\{X_n\}_n$ stabilizes at some value $X_\infty$.
Then there exists a character $\chi\colon \pi_1^{\ab}(X_\infty) \to \frako_L^\times$ such that 
\[
f(x) = \begin{cases} \chi(\Frob_x) & x \in X_\infty^\circ \\
0 & x \in (X \setminus X_\infty)^\circ.
\end{cases}
\]
\end{cor}
\begin{proof}
This is immediate from Corollary~\ref{C:rank-1 assembly infinite level}.
\end{proof}

\begin{remark} \label{R:rank 1 limitation}
It is natural to try to formulate analogous results for representations of $\pi_1$ of dimension greater than 1, taking $f$ to be a function mapping points of $X^\circ$ to Frobenius characteristic polynomials; a successful such effort would simply our proof by making the restriction $e=1$ redundant. 
However, this runs into some difficulties stemming from the fact that in the Brauer--Nesbitt theorem, identifying the Frobenius characteristic polynomials of a group representation (over a field) only pins down the representation up to semisimplification. In particular, 
even if we manage to produce representations $\rho_n$ with the correct Frobenius characteristic polynomials modulo $\varpi^n$, we cannot guarantee the existence of a coherent sequence of such representations.

A version of this issue has been addressed in \cite{drinfeld-deligne}
by introducing some additional hypotheses on the covers $Y_n$ in order to control $\bigcap_n \pi_1(Y_n)$; the idea is to ensure that for each $n$ the mod-$p^n$ representation is limited to some fixed finite set, which then allows for a compactness argument. (Note that to even form $\bigcap_n \pi_1(Y_n)$, we must insist that the $Y_n$ all be finite \'etale over some fixed open subspace $X_0$ of $X$, and then the intersection can be taken within $\pi_1(X_0)$.)

While the resulting argument is for $\ell$-adic rather than $p$-adic representations, under some tameness hypotheses the argument can be adapted to $p$-adic representations (Lemma~\ref{L:isoclinic companions}). By contrast, the preceding results will be applied to representations that do include some wild ramification, and we do not have an analogue of Drinfeld's argument in mind in that setting.
\end{remark}

\subsection{Setup}

We now establish running hypotheses for the remainder of \S\ref{sec:companion points},
and set some additional notation. We also bring in some new information from the uniformity of jumping loci (\S\ref{subsec:jumping loci}). 

\begin{hypothesis} \label{H:companion points new}
For the remainder of \S\ref{sec:companion points},
fix a positive integer $d$; assume that Theorem~\ref{T:companion} holds in all cases when $\dim(X) \leq d$; assume Hypothesis~\ref{H:companion points1} with $\dim(S) = d$;
and assume Hypothesis~\ref{H:companion points} with $e=1$.

We retain all notation set in \S\ref{sec:companion setup} and 
in the statement of Proposition~\ref{P:carryover data}.
Note that we will not make any further reference to the algebraic stacks considered in
\S\ref{subsec:crystals on log}--\S\ref{subsec:universal unit root}.
\end{hypothesis}

\begin{defn} \label{D:jumping divisor}
Define a \emph{jumping divisor} to be a divisor $W$ in $X$ for which $N_y(\calE) \neq N_{\eta_X}(\calE)$ for all $y \in W$;
any such divisor is quasifinite over $S$.
By Proposition~\ref{P:uniform jump bound1}, $\deg(W \to S)$ is bounded over all jumping divisors $W$. Let $V$ be the (open) complement in $X$ of the maximal jumping divisor; by definition, $V_x \subseteq V \times_S x$ for all $x \in S^\circ$.
\end{defn}

\begin{defn} \label{D:wedge constant twist}
For any locally closed subscheme $D$ of $X$ of dimension at most $d$, let $\calF_D \in \FIsoc^\dagger(D) \otimes \overline{\QQ}_p$ be the semisimple crystalline companion of $\calE|_D$ (Hypothesis~\ref{H:companion points new})
Since $\calF_D$ is $E$-algebraic, we may apply Lemma~\ref{L:bound coefficient extension}
to deduce that $\calF_D \in \FIsoc^\dagger(D) \otimes L$.
For $i=1,\dots,r$, let $\tilde{\calF}_{i,D}$ be the constant twist of $\wedge^i \calF_D$ by $\lambda_i^{-1}$.

Now suppose that there exists $y \in D$ such that $N_y(\calE) = N_{\eta_X}(\calE)$. 
For some open dense subscheme $V_{i,D}$ of $D$, we can associate to $\wedge^i \calF_D$ a unit-root representation $\rho_{i,D}\colon \pi_1(V_{i,D}) \to \GL_{e_i}(L)$ (Definition~\ref{D:unit-root subobject}); once Corollary~\ref{C:open constant NP} is available, we will be able to take $V_{i,D} = V \times_X D$.
By Theorem~\ref{T:Tsuzuki minimal slope} and Lemma~\ref{L:contracted section2}(b), $\rho_{i,D}$ is semisimple.

An important special case is when $D = X \times_S x$ for some $x \in S^\circ$ (in which case $N_{s(x)}(\calE) = N_{\eta_X}(\calE)$). In this case, we write $x$ in place of $D$ in the previous notations for consistency with Definition~\ref{D:companions on fibers twist}.
\end{defn}

\subsection{Construction of the unit-root representation}

Using the technique of \S\ref{subsec:rank 1 objects} and the results of the stack-theoretic analysis, we construct a candidate for the unit-root representation of the crystalline companion of $\calE$. 
At this point we must already take advantage of working in the context of an induction on dimension (Hypothesis~\ref{H:companion points new}): while Lemma~\ref{L:rank-1 assembly finite level} is formulated in terms of restrictions to curves, verifying its input hypotheses will require data from restrictions to divisors in $X$.

\begin{prop} \label{P:putative unit root}
There exist an open dense subset $V_1$ of $X$ containing $V$ and a character
$\chi\colon \pi_1^{\ab}(V_1) \to \frako_L^\times$ for which the following statements hold.
\begin{enumerate}
\item[(a)]
For $z \in X^\circ$, $z \in V_1$ if and only if $N_z(\calE)$ has least slope $0$ with multiplicity $1$.
\item[(b)]
For $z \in V_1^\circ$, $\chi(z)$ equals the unit-root Frobenius eigenvalue of $\calE_z$ (identified with an element of $L$ using the chosen place $v$).
\end{enumerate}
\end{prop}
\begin{proof}
Let $f\colon X^\circ \to \frako_L^\times \cup \{0\}$ be the function such that $f(z)$ equals the unit-root Frobenius eigenvalue of $\calE_z$ when the latter exists and $0$ otherwise. 
By Lemma-Definition~\ref{L:generic Newton polygon}, the function $f$ satisfies condition (i) of Corollary~\ref{C:rank-1 assembly infinite level}.
To check condition (ii), set notation as in Definition~\ref{D:wedge constant twist} with $D := C$; set $C_0 := V_{1,D}$ and $\chi_C := \rho_{1,C}\colon \pi_1^{\ab}(C_0) \to L^\times$;
note that $\chi_C$ has image in $\frako_L^\times$ since $\pi_1^{\ab}(C_0)$ is compact;
and apply Lemma~\ref{L:no extension filtration} to see that $\chi_C$ is ramified at every point of $C \setminus C_0$.

To check condition (iii), we must work around the issue raised in Remark~\ref{R:no density 1 version}.
We choose a dominant \'etale morphism $Y_n \to X$ so that each of the representations
$\rho_{n,1,j_n}$ in Proposition~\ref{P:carryover data}(c) restricts trivially to $\pi_1(Y_n)$
(this is possible because $J_n$ is finite).
For each $y \in Y_n^\circ$ lying over $z \in X$ with $f(z) \neq 0$,
choose an irreducible divisor $D$ in $X$ dominating $S$ and containing $z$.
By Proposition~\ref{P:carryover data}(c), we have $f(y') = 1$ for every $y' \in (Y_n \times_X D)^\circ$
lying over an element of $S$ lying in some subset of $S^\circ$ of density 1;
these points form a subset of $(Y_n \times_X D)^\circ$ of density 1.
By Chebotaryov, this implies the same conclusion for all $y' \in (Y_n \times_X D)^\circ$,
and hence in particular for $y$.

At this point, the hypotheses of Corollary~\ref{C:rank-1 assembly infinite level} are satisfied.
Now note that for any irreducible divisor $Z$ on $X$, either $Z$ dominates $S$ or $Z$ is the preimage in $X$ of some divisor $W$ in $S$. If $f(z) = 0$ for all $z \in Z^\circ$, then we cannot be in the latter case because $N_{s(x)}(\calE) = N_{\eta_X}(\calE)$ for $x \in W$; hence $Z$ must be a jumping divisor.
 We conclude that $\bigcup_n (X \setminus X_n) \subseteq X \setminus V$ is a finite union of irreducible divisors of $X$; together with the previous paragraph, this verifies the hypotheses of Corollary~\ref{C:rank-1 assembly infinite level stable} and thus yields the desired result.
\end{proof}

At this point, we close the gap between Lemma{-}Definition~\ref{L:generic Newton polygon}(c) and Lemma~\ref{L:convergent Newton polygon}(c). See \S\ref{subsec:Newton} for further discussion.
\begin{cor} \label{C:open constant NP}
The set of $y \in X^\circ$ for which $N_y(\calE) = N_{\eta_X}(\calE)$ is equal to $V^\circ$. In particular, $V_x = V \times_S x$ for all $x \in S^\circ$.
\end{cor}
\begin{proof}
By Definition~\ref{D:jumping divisor}, 
if $y \in (X \setminus V)^\circ$ then $N_y(\calE) \neq N_{\eta_X}(\calE)$. Conversely, for $y \in V^\circ$, for each $i \in \{1,\dots,r-1\}$ which is the $x$-coordinate of a vertex of $N_{\eta_X}(\calE)$,
we may apply Proposition~\ref{P:putative unit root} with $\calE$ replaced by $\wedge^i \calE$ (since the condition $e=1$ in Hypothesis~\ref{H:companion points new} remains valid) 
to deduce that
$N_y(\wedge^i \calE)$ and $N_{\eta_X}(\wedge^i \calE)$ have the same least slope (with multiplicity $1$).
By Lemma{-}Definition~\ref{L:generic Newton polygon}(b),
this implies that $N_y(\calE) = N_{\eta_X}(\calE)$.
\end{proof}

\subsection{Descending the generic fiber}

By combining the construction of a candidate for the unit-root representation of the crystalline companion with the minimal slope theorem, we improve upon Proposition~\ref{P:carryover data}
to give a \emph{single} candidate for the generic fiber of the crystalline companion of $\calE$.

\begin{prop} \label{P:generic companion}
For $\eta^{\perf}$ the perfect closure of $\eta$,
there exists an absolutely irreducible object $\calF_{\eta^{\perf}} \in \FIsoc^{\dagger}(X \times_S \eta^{\perf}) \otimes L$ docile along $Z \times_S \eta^{\perf}$ with the following properties.
\begin{enumerate}
\item[(a)]
The object $\calF_{\eta^{\perf}}$ is docile along $Z \times_S \eta^{\perf}$,
absolutely irreducible, and its  generic Newton polygon is $N_{\eta_X}(\calE)$.
\item[(b)]
The Newton polygon of $\calF_{\eta^{\perf}}$ is constant on $V \times_S \eta^{\perf}$.
\item[(c)]
For $i=1,\dots,r$, let
$\tilde{\calF}_{i,\eta^{\perf}}$ be the constant twist of $\wedge^i \calF_{\eta^{\perf}}$ by $\lambda_i^{-1}$.
Let $\rho_i\colon \pi_1(V \times_S \eta^{\perf}) \to \GL_{e_i}(L)$
be the unit-root representation of $\tilde{\calF}_{i,\eta^{\perf}}$ (Definition~\ref{D:unit-root subobject}).
Then $\rho_i$ is semisimple.
\item[(d)]
For $i=1,\dots,r$, choose a stable lattice in $\rho_i$ and,
for each positive integer $n$,
let $\rho_{i,n}\colon \pi_1(V \times_S \eta^{\perf}) \to \GL_{e_i}(\frako_L/p^n \frako_L)$ be the reduction modulo $p^n$.
Then for every open dense subset $V_{i,n}$ of $V$ to which $\rho_{i,n}$ extends, for every $z \in V_{i,n}^\circ$ lying over $x \in S^\circ$, the characteristic polynomial of $\rho_{i,x}(\Frob_z)$ reduces modulo $p^n$ to the characteristic polynomial of $\rho_{i,n}(\Frob_z)$.
\end{enumerate}
\end{prop}
\begin{proof}
Let $I$ be the set of $i \in \{1,\dots,r\}$ which occur as the $x$-coordinate of a vertex of $N_{\eta_X}$.
We start with the objects $\calF^{(j)}_{\tilde{\eta}}$ for $j \in J$ 
from Proposition~\ref{P:carryover data}; we may invoke Proposition~\ref{P:carryover data}(d) thanks to Proposition~\ref{P:putative unit root} (applied to $\wedge^i \calE$ for each $i \in I$).
Since these objects are absolutely irreducible and $1 \in I$, 
we may apply Theorem~\ref{T:Tsuzuki minimal slope} to deduce that all of the
$\calF^{(j)}_{\tilde{\eta}}$ are restrictions of a single object $\calF_{\eta^{\perf}} \in \FIsoc^{\dagger}(X \times_S \eta^{\perf}) \otimes L$;
by the same token, for each $i \in \{1,\dots,r\}$,
the representations $\rho_i^{(j)}$ are all restrictions of a single representation $\rho_i$.
By Proposition~\ref{P:carryover data}(a), we obtain (a).

For each $i \in I$, the application of Proposition~\ref{P:putative unit root} in the previous paragraph shows that the representation $\rho_i$ extends across $V$. By Lemma~\ref{L:no extension filtration} and Proposition~\ref{P:unit-root rep}(a), we deduce (b).
By Corollary~\ref{C:minimal slope semistable} and
Proposition~\ref{P:generic companion}(a), we deduce (c).

Applying Proposition~\ref{P:carryover data}(b,c), we obtain (d) for some choice of $V_{i,n}$,
modulo the point that $\rho_i$ is currently only defined on some open dense subset of $V \times_S \eta^{\perf}$. To upgrade the conclusion to a general choice of $V_{i,n}$ (with the same proviso),
it will suffice to check the claim for a specified $z \in V^\circ$
lying in a specified superset $V'_{i,n}$ of $V_{i,n}$ to which $\rho_{i,n}$ extends.
To do this, choose a curve $C$ in $V'_{i,n}$ passing through $z$ and some closed point in $V_{i,n}$.
Choose a stable lattice in $\rho_{i,C}$ and let $\rho_{i,C,n}\colon \pi_1(C) \to \GL_{e_i}(\frako_L/p^n \frako_L)$ be the reduction modulo $p^n$. By comparing $\rho_{i,C,n}$ with the restriction of $\rho_{i,n}$ to $\pi_1(C)$ using Chebotaryov density, we deduce the claim.
\end{proof}

\subsection{Matching with divisors}
\label{subsec:extend slope steps}

We next extend of the successive quotients of the slope filtration of the candidate generic companion $\calF_{\eta^{\perf}}$ across the fibration. The strategy here is to use the hypothesis that crystalline companions exist on divisors in $X$ to enlarge the subsets $V_{i,n}$ in Proposition~\ref{P:generic companion}(b). 

\begin{remark} \label{R:density for comparison}
Let $W$ be an open dense subset of $X$.
Let $\rho_1, \rho_2\colon \pi_1(W) \to \GL_m(\frako_L/p^n \frako_L)$ be two continuous representations.
Let $G \subseteq \GL_m(\frako_L/p^n \frako_L) \times \GL_m(\frako_L/p^n \frako_L)$ be the image of the product representation $\rho_1 \times \rho_2$.
By Chebotaryov density, the following statements are equivalent.
\begin{itemize}
\item
For every $g \in G$, the characteristic polynomials of the two components of $g$ are congruent modulo $p^n$.
\item
For all $y \in W^\circ$, the characteristic polynomials of $\rho_1(\Frob_y)$
and $\rho_2(\Frob_y)$ are congruent modulo $p^n$.
\item
For all $y$ in a subset of $W^\circ$ of density 1, the characteristic polynomials of $\rho_1(\Frob_y)$
and $\rho_2(\Frob_y)$ are congruent modulo $p^n$.
 (It is even sufficient to take a subset of lower density at least $1-(\#\GL_m(\frako_L/p^n \frako_L))^{-2}$.)
\end{itemize}
\end{remark}

The following can be viewed as a variant of Lemma~\ref{L:pin monodromy restriction}.
\begin{lemma} \label{L:subset to pin restriction}
For $\rho_i$ as in Proposition~\ref{P:generic companion},
for any irreducible divisor $D$ in $X$,
the restrictions of $\rho_i$ and $\rho_{i,D}$ to $\pi_1(D\times_S \eta^{\perf} )$ have the same semisimplification.
\end{lemma}
\begin{proof}
For each $n$, Proposition~\ref{P:generic companion}(d) and Remark~\ref{R:density for comparison}  imply that for every $g \in \pi_1( V_{i,n} \times_X D)$, the characteristic polynomial of $\rho_{i,D}(g)$ reduces modulo $p^{n}$ to the characteristic polynomial of $\rho_{i,n}(g)$. In the limit, this implies that for every $g \in \pi_1(D \times_S \eta^{\perf})$, the characteristic polynomials of $\rho_i(g)$ and $\rho_{i,D}(g)$ coincide.
We deduce the claim from Brauer--Nesbitt.
\end{proof}

\begin{lemma} \label{L:preserve semisimple}
Let $G$ be a profinite group, let $m$ be a positive integer, and let $\rho\colon G \to \GL_m(\frako_L)$ be a continuous homomorphism such that the action of $G$ on $L^m$ via $\rho$ is semisimple. Then there exists a positive integer $n_0$ such that for all $n \geq n_0$,
for $\rho_n\colon G \to \GL_m(\frako_L/p^n \frako_L)$ the mod-$p^n$ reduction of $\rho$,
for any closed subgroup $G'$ of $G$ for which $\rho_n(G') = \rho_n(G)$, the action of $G'$ on $L^m$ via $\rho$ is again semisimple.
\end{lemma}
\begin{proof}
The action of $G$ on $L^m$ is semisimple if and only if $\End(L^m)^G$ is a semisimple $L$-algebra. It will thus suffice to ensure that $\End(L^m)^G = \End(L^m)^{G'}$.

We define a sequence $U_1, \dots, U_l$ of closed-open subsets of $G$ and a strictly decreasing sequence of $L$-vector spaces
\[
\End(L^m) = E_0 \supset \cdots \supset E_l = \End(L^m)^G
\]
such that for $j \leq i$, $E_i$ is fixed by every $g \in U_j$.
Given $E_i \neq \End(L^m)^G$, pick some $v \in E_i \setminus \End(L^m)^G$; then there exists some closed-open subset $U_{i+1}$ of $G$ such that $v$ is not fixed by any $g_{i+1} \in U_{i+1}$.
Since $\dim_L \End(L^m) < \infty$, the construction must terminate.

We now choose $n_0$ so that each of $U_1,\dots,U_l$ maps to a single element under $\rho_{n_0}$. This has the desired effect.
\end{proof}

\begin{prop} \label{P:uniform factorization}
For $i \in \{1,\dots,r\}$, 
the representation $\rho_i$ defined in Proposition~\ref{P:generic companion} factors through $\pi_1(V)$.
\end{prop}
\begin{proof}
Choose $n_0$ as in Lemma~\ref{L:preserve semisimple}; it suffices to check that
for $n \geq n_0$, we may take $V_{i,n} = V$ in Proposition~\ref{P:generic companion}(d).
By Proposition~\ref{P:generic companion}(b), 
we need only show that
for each irreducible divisor $W$ in $S$, $\rho_{i,n}$ is unramified along the divisor $X \times_S W$ in $X$. 

Suppose to the contrary that this fails for some $W$.
Choose a finite subset $H$ of $(W \cap V_{i,n})^\circ$ such that the image of $\rho_{i,n}$ is generated by the images of $\Frob_y$ for $y \in H$.
We can then choose an irreducible divisor $D$ in $X$ containing $H$ such that the restriction of $\rho_{i,n}$ to $\pi_1(V_{i,n} \times_X D)$ is ramified along some component of $V_{i,n} \times_X D \times_S W$.

By Proposition~\ref{P:generic companion}(c) and Lemma~\ref{L:preserve semisimple}, the restriction of $\rho_i$ to $\pi_1(D \times_S \eta^{\perf})$ is semisimple. 
Meanwhile, $\rho_{i,D}$ is semisimple (see Definition~\ref{D:wedge constant twist}),
so its restriction along the surjection $\pi_1(D \times_S \eta^{\perf}) \to \pi_1(V \times_X D)$
is also semisimple. By Lemma~\ref{L:subset to pin restriction}, we conclude that the restrictions of $\rho_i$ and $\rho_{i,D}$ to $\pi_1(D \times_S \eta^{\perf})$ are isomorphic;
in particular, the restriction of $\rho_i$ to $\pi_1(D \times_S \eta^{\perf})$ factors through $\pi_1(V \times_X D)$.
This in turn implies that the restriction of $\rho_{i,n}$ to $\pi_1(D \times_S \eta^{\perf})$ factors through $\pi_1(V \times_X D)$, contradicting our choices of $W$ and $D$.
\end{proof}

\begin{cor} \label{C:extend slope steps}
Let $\calG_{1,\eta^{\perf}},\dots,\calG_{l,\eta^{\perf}}$ be the successive quotients of the slope filtration of $\calF_{\eta^{\perf}}$.
\begin{enumerate}
\item[(a)]
The objects $\calG_{1,\eta^{\perf}},\dots,\calG_{l,\eta^{\perf}}$ extend to semisimple objects $\calG_1,\dots,\calG_l$ of $\FIsoc(V) \otimes L$.
\item[(b)]
For each $x \in S^\circ$,
the restrictions of $\calG_1,\dots,\calG_l$ to $\FIsoc(V_x) \otimes L$ are isomorphic to the successive quotients of the slope filtration of $\calF_x$. 
\item[(c)]
Let $\calG \in \FIsoc(V) \otimes L$ be the direct sum of the $\calG_i$. Then for every $y \in V^\circ$, $\calG_y$ is a crystalline companion of $\calE_y$.
\end{enumerate}
\end{cor}
\begin{proof}
By Proposition~\ref{P:uniform factorization} and Proposition~\ref{P:unit-root rep}(a),
for $i=1,\dots,r$,
the unit-root subobject of $\tilde{\calF}_{i,\eta^{\perf}}$ extends to an object of $\FIsoc(V) \otimes L$, which by Proposition~\ref{P:generic companion}(c) is semisimple.
Using Remark~\ref{R:slope filtration exterior power}, we deduce (a).
We then deduce (b) from Proposition~\ref{P:generic companion}(d), which in turn implies (c).
\end{proof}

\subsection{Final glueing}

We are now ready to glue together our candidate generic fiber with the object $\calG$
to obtain the desired crystalline companion.
\begin{prop} \label{P:full companion on fibration} 
The object $\calF_{\eta^{\perf}}$ extends to an object $\calF$ of
$\FIsoc^\dagger(X) \otimes L$ which is a crystalline companion of $\calE$.
\end{prop}
\begin{proof}
Define the semisimple object
$\calG \in \FIsoc(V) \otimes L$ as in Corollary~\ref{C:extend slope steps}(c).
For any divisor $D$ in $X$ meeting $V$, $D \times_S \eta^{\perf} = (D \times_S \eta)^{\perf}$ is a finite disjoint union of perfect points. 
By Corollary~\ref{C:extend slope steps}(b) and Remark~\ref{R:split filtration}, the restrictions of $\calF_{\eta^{\perf}}$ and $\calG$ to $\FIsoc(D \times_S \eta^{\perf}) \otimes L$ are isomorphic.

Meanwhile, we have a semisimple object $\calF_D \in \FIsoc^\dagger(D) \otimes L$ from Definition~\ref{D:wedge constant twist}. 
Since $\calF_D$ has constant Newton polygon on $V$, the restriction of $\calF_D$ to 
$\FIsoc(V \times_X D) \otimes L$ admits a slope filtration by 
Corollary~\ref{C:slope filtration},
whose successive quotients are semisimple by Corollary~\ref{C:minimal slope semistable}.
By combining Proposition~\ref{P:unit-root rep}(a), Chebotaryov density, and Brauer--Nesbitt, 
we see that each successive quotient of $\calF_D$ is \emph{isomorphic} to a corresponding
component of $\calG$ (since both are known to be semisimple).
By Remark~\ref{R:split filtration}, the slope filtration of $\calF_D$ splits uniquely in $\FIsoc(D \times_S \eta^{\perf}) \otimes L$;
consequently, the restrictions of $\calF_D$ and $\calG$ to $\FIsoc(D \times_S \eta^{\perf}) \otimes L$ are in fact isomorphic.

We now have an object of
\[
(\FIsoc(\overline{X}^{\log} \times_S \eta^{\perf}) \otimes L)
\times_{\FIsoc(D \times_S \eta^{\perf}) \otimes L}
(\FIsoc^\dagger(D) \otimes L);
\]
for a suitable choice of $D$, we may apply Proposition~\ref{P:descent using divisor} and Proposition~\ref{P:descent using divisor perfect} to obtain an object $\calF \in \FIsoc^\dagger(X) \otimes L$ which is docile along $Z$.
By Corollary~\ref{C:extend slope steps}(c), $\calF|_V$
is a crystalline companion of $\calE|_V$;
using Lemma~\ref{L:extension of companions}, we may upgrade this conclusion to the desired result.
\end{proof}

\section{Companions and corollaries}
\label{sec:corollaries}

With the key construction in hand, we complete the construction of crystalline companions, and record some corollaries.

\begin{hypothesis}
Throughout \S\ref{sec:corollaries}, assume that $k$ is finite.
\end{hypothesis}

\subsection{Proof of the main theorems}
\label{subsec:main results}

We complete the proofs of Theorem~\ref{T:deligne} and Theorem~\ref{T:companion}. While the heavy lifting was already done in \S\ref{sec:companion points}, there remains a bit of work to eliminate the restriction that the least generic Newton slope occurs with multiplicity 1 (Hypothesis~\ref{H:companion points new}).

\begin{theorem} \label{T:algebraic companions}
Any algebraic \'etale coefficient object on $X$ admits all crystalline companions. 
\end{theorem}
\begin{proof}
We proceed by induction on $\dim(X)$, the case $\dim(X) = 1$ being
Corollary~\ref{C:algebraic companion}.
Let $\calE$ be an algebraic \'etale coefficient object of rank $r$ on $X$;
as per Definition~\ref{D:algebraic}, $\calE$ is $E$-algebraic for some number field $E$.
By Lemma~\ref{L:alteration reduction}, we may check the claim after replacing $X$ with an alteration or an open dense subspace;
we may thus assume that $\calE$ is absolutely irreducible and
(by Lemma~\ref{L:alter to docile}) docile.
By making a constant twist, we may ensure that the Newton polygon $N_\eta(\calE)$ (Lemma{-}Definition~\ref{L:generic Newton polygon})
at the generic point $\eta$ of $X$ has least slope 0 with some multiplicity $e$. 
By Lemma~\ref{L:monodromy}, we may reduce to the case where $\overline{G}(\calE)$ is connected; then $\overline{G}(\wedge^r \calE)$ is connected and semisimple, hence trivial, meaning that $\wedge^r \calE$ is constant.

At this point, we can use Corollary~\ref{C:stable curve fibration}
and Remark~\ref{R:contracted section} to pull back along a dominant generically \'etale morphism so as to put ourselves in the situation of Hypothesis~\ref{H:companion points1}; that is, $X$ is the unpointed locus of some smooth curve fibration $f\colon \overline{X} \to S$, and there exists a section $s\colon S \to X$ such that $s^* \calE$ is constant.
Using the last sentence of Remark~\ref{R:contracted section},
we may further ensure that the Newton polygon of $s^* \calE$
is identically equal to $N_\eta(\calE)$.
By shrinking $S$, we may ensure that the fibers of $f$ have constant $p$-rank (e.g., by pushing forward the unit crystalline coefficient from $\overline{X}$ to $S$ 
and applying Lemma~\ref{L:convergent Newton polygon}(c) to make the Newton polygon of the result constant).

We now distinguish cases based on the value of $e$.
\begin{itemize}
\item
If $e=1$, then Hypothesis~\ref{H:companion points new} holds, so we may apply Proposition~\ref{P:full companion on fibration}
to obtain the desired crystalline companion. 

\item
If $1 < e < r$, then $\wedge^e \calE$ is semisimple and  $N_\eta(\wedge^e \calE)$ has least slope 0 with multiplicity 1; 
consequently, $\wedge^e \calE$ admits a unique irreducible subobject $\calE'$ which has least slope 0 with multiplicity 1
(so in particular $\overline{G}(\calE')$ is not finite).
By the previous paragraph, $\calE'$ admits a crystalline companion; per Remark~\ref{R:same Tannakian category2}, $\calE^\dual \otimes \calE$ also admits a crystalline companion. We may then conclude using Lemma~\ref{L:deduce from exterior power}.

\item
If $e = r$, then $\calE$ is unit-root and we may produce a crystalline companion using a version of Drinfeld's argument in the \'etale case. See Lemma~\ref{L:isoclinic companions}.
\qedhere
\end{itemize}
\end{proof}

\begin{lemma} \label{L:deduce from exterior power}
In the context of the proof of Theorem~\ref{T:algebraic companions}, if $1<e<r$ and $\calE^\dual \otimes \calE$ admits a crystalline companion, then so does $\calE$.
\end{lemma}
\begin{proof}
Per Remark~\ref{R:same Tannakian category2}, the homomorphism $\overline{G}(\calE) \to \overline{G}(\calE^\dual \otimes \calE)$ is surjective with kernel equal to the center of $\overline{G}(\calE)$, which is a finite cyclic group of some order $n$.
Let $\calG$ be a semisimple crystalline companion of $\calE^\dual \otimes \calE$. 

At this point, we may work \'etale-locally on $S$.
For any specified $d$, we may thus assume that there exists a divisor $D$ in $X$ which is finite \'etale over $S$ of degree $\geq d$
and which passes through some finite subset of $X^\circ$ as in Lemma~\ref{L:pin monodromy restriction}; the latter restriction ensures that the maps $G(\calE|_D) \to G(\calE)$ and $\overline{G}(\calE|_D) \to \overline{G}(\calE)$ are isomorphisms.
In particular, the restriction map $H^0(X, \calE^\dual \otimes \calE) \to H^0(D, \calE^\dual \otimes \calE)$ is an isomorphism,
so every irreducible constituent of $\calE^\dual \otimes \calE$ remains irreducible upon restriction to $D$.
By Lemma~\ref{L:companion irreducible}(a) the same is true of $\calG$; in particular, $\calG|_{D}$ is again semisimple.

In the context of the proof of
Theorem~\ref{T:algebraic companions}, $\calE|_D$ is known to have a crystalline companion $\calF_D$, which by Lemma~\ref{L:companion irreducible}(a) is again irreducible.
By Lemma~\ref{L:companion irreducible}(d) again (since both objects are semisimple), $\calG|_D$ is isomorphic to $\calF_D^\dual \otimes \calF_D$.

By results of Chin and D'Addezio (see Theorem~\ref{T:same component group} and Theorem~\ref{T:same connected monodromy})
the algebraic groups $G(\calE^\dual \otimes \calE)$ and $G(\calG)$ have isomorphic component groups,
and their neutral components may be viewed as base extensions of a common group defined over some number field in a manner compatible with their natural representations.
Moreover, the corresponding statement for $\overline{G}(\calE^\dual \otimes \calE)$ and $\overline{G}(\calG)$ also holds, and both of these statements also apply after restriction to $D$.
Consequently, the map $G(\calG) \to G(\calG|_D)$ is an isomorphism;
since $\calG|_D$ visibly carries the structure of an Azumaya algebra object in its category (Definition~\ref{D:Azumaya algebra}), so then does $\calG$.
We may now apply Proposition~\ref{P:descent using divisor}
and Proposition~\ref{P:split Azuyama algebra using divisor} to obtain an isomorphism $\calG \cong \calF^\dual \otimes \calF$ for some $\calF \in \FIsoc^\dagger(X) \otimes L$.

Since we already have Theorem~\ref{T:companion} for $\ell' \neq p$ by \cite[Theorem~3.5.2]{kedlaya-companions}, $\calF$ itself admits a companion $\calE'$ in the category of $\calE$,
which is again docile along $Z$. 
By construction, $\calE^\dual \otimes \calE$ and $(\calE')^\dual \otimes \calE'$
are both companions of $\calG$ and hence are companions of each other in the same category;
in particular, their trace-zero components are also companions in the same category.
In the context of the proof of Theorem~\ref{T:algebraic companions},
$\overline{G}(\calE)$ is connected and hence the trace-zero component of $\calE^\dual \otimes \calE$ is absolutely irreducible; we may now apply Lemma~\ref{L:companion irreducible}(d) to see that $\calE^\dual \otimes \calE$ and $(\calE')^\dual \otimes \calE'$ are isomorphic.

Let $W_X$ denote the Weil group of $X$.
Since $\calE$ and $\calE'$ are \'etale coefficient objects, they correspond to continuous representations
$\rho\colon \pi_1(W_X) \to \GL(V)$ and 
$\rho'\colon \pi_1(W_X) \to \GL(V')$ for some finite-dimensional $\overline{\QQ}_\ell$-vector spaces $V, V'$. By the previous paragraph, the induced representations on $V^\dual \otimes V$ and $(V')^\dual \otimes V$ are isomorphic, and moreover the isomorphism is unique up to scaling. Consequently,
the fixed subspaces of
\[
(V^\dual \otimes V)^\dual \otimes ((V')^\dual \otimes V')
\cong (V^\dual \otimes V')^\dual \otimes (V^\dual \otimes V')
\]
are both one-dimensional, so this subspace must be idempotent for the multiplication on the right-hand side;
it thus contains a projector onto a one-dimensional $W_X$-stable subspace of $V^\dual \otimes V'$.
If we denote this subspace by $V_0$, then $V_0 \otimes V \cong V'$; that is, $V_0$
corresponds to a coefficient object $\calL$ of rank 1 in the category of $\calE$ such that
$\calE \otimes \calL \cong \calE'$. By Corollary~\ref{C:rank 1 companions}, $\calL$ admits a crystalline companion; twisting $\calF$ by the dual of the latter yields a crystalline companion of $\calE$.
\end{proof}

\begin{lemma} \label{L:surjective pro-p}
Let $G \to H$ be a surjective homomorphism of topologically 
finitely presented pro-$p$ groups. Then $G \to H$ is an isomorphism if and only if
the maps $H^i(H, \FF_p) \to H^i(G, \FF_p)$ are isomorphisms for $i=1,2$.
\end{lemma}
\begin{proof}
The $\FF_p$-dimension of $H^1(G, \FF_p)$ equals the minimum number of generators of $G$,
while the $\FF_p$-dimension of $H^2(G, \FF_p)$ equals the minimum number of relations among a minimal system of generators of $G$
(see \cite[\S 4.2, 4.3]{serre-galois}). 
Consequently, the isomorphism for $i=1$ ensures that a minimal system of generators of $G$ remains minimal as a system of generators of $H$, while the isomorphism for $i=2$ ensures that no extra relations are needed.
\end{proof}

\begin{lemma} \label{L:isoclinic companions}
In the context of the proof of Theorem~\ref{T:algebraic companions}, if $e=r$, then
$\calE$ admits a crystalline companion.
\end{lemma}
\begin{proof}
Set notation, including the finite extension $L$ of $\QQ_p$, as in
\S\ref{sec:companion setup}.
Let $\varpi$ be a uniformizer of $L$.

For each curve $C$ in $\overline{X}$, the crystalline companion $\calF_{C \times_{\overline{X}} X}$
of $\calE|_{C \times_{\overline{X}} X}$ is docile (by Corollary~\ref{C:companions are tame}).
Corollary~\ref{C:unit-root docile} implies that $\calF_{C \times_{\overline{X}} X}$ extends over $C$,
as then does the semisimplification of $\calE|_{C \times_{\overline{X}} X}$ (Lemma~\ref{L:extension of companions}).

On one hand, we may apply \cite[Theorem~2.5]{drinfeld-deligne} to produce a companion of $\calE$ in the same category which extends across $\overline{X}$; by Lemma~\ref{L:extension of companions} again, $\calE$ itself extends across $\overline{X}$.
That is, we may assume hereafter that $Z = \emptyset$.

On the other hand, applying Proposition~\ref{P:unit-root rep}(a)
produces a skeleton sheaf (see Definition~\ref{D:skeleton sheaf}) of $L$-local systems on $X$. 
This will allow us to apply Drinfeld's method as exposed in \cite[\S 7]{cadoret}, with only minor modifications to accommodate the case $\ell=p$.

For $\overline{s} \to S$ a geometric point, let $\pi_1^{p}(X \times_S \overline{s})$ denote the maximal pro-$p$ quotient 
of $\pi_1(X \times_S \overline{s})$. We have natural identifications 
\[
H^i(\pi_1^{p}(X \times_S \overline{s}), \FF_p) \cong H^i_{\et}(X \times_S \overline{s}, \underline{\FF_p}).
\]
Since we are assuming that the fibers of $f$ have constant $p$-rank, Lemma~\ref{L:surjective pro-p} implies that for $\overline{\eta} \to S$ a geometric point lying over the generic point of $S$, the surjective maps
\begin{equation} \label{eq:tame specialization at p}
\pi_1^{p}(X \times_S \overline{\eta}) \to \pi_1^{p}(X\times_S \overline{s})
\end{equation}
provided by the specialization theorem \cite[Expos\'e X, Corollaire~2.4]{grothendieck-sga1-tame} are always isomorphisms.
(This is the only substantial accommodation needed to handle $\ell=p$.)

We first argue as in \cite[Lemme~7.4]{cadoret} that after shrinking $S$, we can find a connected finite \'etale cover $X_1$ of $X$ which trivializes all of the Frobenius characteristic polynomials modulo $\varpi$. 
Let $\eta$ be the generic point of $S$ and let $\overline{\eta}$ be a geometric point above $\eta$.
Since $\pi_1(X \times_S \overline{\eta})$ is topologically finitely generated \cite[Expos\'e X, Th\'eor\`eme~2.9]{grothendieck-sga1-tame},
all continuous homomorphisms $\pi_1(X \times_S \overline{\eta}) \to \GL_r(\frako_L/\varpi\frako_L)$ factor through some characteristic open subgroup $N$. We may thus choose a cover of $X$ whose pullback along $\overline{\eta} \to S$ corresponds to a subgroup of $N$ and which trivializes the pullback of $\calE$ along $s$ modulo $\varpi$ (note that in our setup the latter only requires a constant field extension).

At this point we may assume that all of the Frobenius characteristic polynomials of $\calE$ are trivial modulo $\varpi$
(that is, all of their roots reduce to 1 modulo $\varpi$).
We then argue as in \cite[Lemme~7.4]{cadoret} that with \emph{no} further shrinking of $S$,
for each $n$ we can find a connected finite \'etale cover $X_n$ of $X$ which trivializes all of the Frobenius characteristic polynomials modulo $\varpi^n$. Namely, let $Z_n$ be the kernel of $\GL_r(\frako_L/\varpi^n \frako_L) \to \GL_r(\frako_L/\varpi \frako_L)$, which is a $p$-group.
Since $\pi_1(X \times_S \overline{\eta})$ is topologically finitely generated and  \eqref{eq:tame specialization at p} is an isomorphism for all $s \in S$, the joint kernel of all continuous homomorphisms $\pi_1(X \times_S \overline{\eta}) \to Z_n$ corresponds to a finite \'etale cover that spreads out over all of $S$.

We next argue as in \cite[Lemme~7.3]{cadoret} to recover an $L$-local system on $X$.
Let $\Pi$ be the intersection of the subgroups $\pi_1(X_n)$ of $\pi_1(X)$; the group $\pi_1(X)/\Pi$ is topologically finitely generated and virtually pro-$p$. In particular, the Frattini subgroup $F$ of $\pi_1(X)/\Pi$ is open
and the set
\[
H := \Hom(\pi_1(X)/\Pi, \ker(\GL_r(\frako_L) \to \GL_r(\frako_L/\varpi \frako_L)) = \varprojlim_n \Hom(\pi_1(X)/\Pi, Z_n)
\]
is compact for the inverse limit of the discrete topologies. For each $x \in X^\circ$, let $H_{x} \subseteq H$ be the closed subset of representations whose Frobenius characteristic polynomial at $x$ matches that of $\calE_x$; by compactness, it will suffice to check that $\bigcap_{x \in W} H_{x} \neq \emptyset$ for every finite subset $W$ of $X^\circ$.
By \cite[Theorem~2.15]{drinfeld-deligne}, there exists a curve $C$ in $X$ containing $W$ such that $\pi_1(C) \to \pi_1(X) \to (\pi_1(X)/\Pi)/F$ is surjective, as then is $\pi_1(C) \to \pi_1(X)/\Pi$; let $K_\Pi$ be the kernel of the latter map. 
A semisimple crystalline companion of $\calE|_C$ corresponds to a semisimple representation $\rho_C\colon \pi_1(C) \to \GL_r(L)$; the restriction of $\rho$ to $K_\Pi$ is both semisimple (because $K_\Pi$ is normal in $\pi_1(C)$) and unipotent (by
Chebotaryov and Brauer--Nesbitt), hence trivial. Hence $\rho \in \bigcap_{x \in W} H_x$ as desired.

We now have an $L$-local system $\rho$ on $X$ whose restriction to any curve $C$ in $X$ has the same Frobenius traces as $\rho_C$.
By Proposition~\ref{P:unit-root rep}(a), this gives rise to a unit-root object $\calF \in \FIsoc(X) \otimes L$, or equivalently in $\FIsoc^\dagger(X) \otimes L$ because $X = \overline{X}$,
which is a companion of $\calE$.
\end{proof}

\begin{cor} \label{C:companion p}
Theorem~\ref{T:companion} holds in all cases.
\end{cor}
\begin{proof}
As noted previously, for $\ell' \neq p$ this is included in \cite[Theorem~3.5.2]{kedlaya-companions}.
For $\ell \neq p, \ell'=p$, this is included in Theorem~\ref{T:algebraic companions}.
For $\ell = \ell' = p$, we may first apply Corollary~\ref{C:deligne prior} to change to the case $\ell \neq p$,
then proceed as before.
\end{proof}

We mention explicitly the following special case of Theorem~\ref{T:companion}.
\begin{cor} \label{C:Galois conjugates}
Let $\calE$ be an $E$-algebraic coefficient object on $X$. Then for any automorphism $\tau$ of $E$,
there exists a coefficient object $\calE_\tau$ such that for each $x \in X^\circ$, 
we have the equality $P((\calE_\tau)_x, T) = \tau(P(\calE_x,T))$ in $E[T]$ (where $\tau$ acts coefficientwise).
\end{cor}

\begin{cor} \label{C:deligne}
Theorem~\ref{T:deligne} holds in all cases.
\end{cor}
\begin{proof}
By Corollary~\ref{C:deligne prior}, parts (i)--(v) hold. To prove (vi), note that (ii) implies that $\calE$
is algebraic, so Corollary~\ref{C:companion p} implies the existence of a crystalline companion $\calF$.
By Lemma~\ref{L:companion irreducible}(a) and (c), $\calF$ is irreducible and $\det(\calF)$ is of finite order.
\end{proof}

\subsection{Newton polygons revisited}
\label{subsec:Newton}

With Theorem~\ref{T:companion} in hand, we can now assert much stronger properties of Newton polygons of Weil sheaves than were asserted in Lemma{-}Definition~\ref{L:generic Newton polygon}. These extend the Grothendieck--Katz semicontinuity theorem and the de Jong--Oort--Yang purity theorem to \'etale coefficients (see \cite[\S 3]{kedlaya-isocrystals}).
The first of these extensions had been informally conjectured by Drinfeld \cite[D.2.4]{drinfeld-pro-simple}.
\begin{theorem}[after Grothendieck--Katz] \label{T:stratification1}
Let $\calE$ be an $E$-algebraic $\ell$-adic coefficient object for some number field $E$, and fix an embedding $E \hookrightarrow \overline{\QQ}_p$. 
Then the function $x \mapsto N_x(\calE)$ from Lemma{-}Definition~\ref{L:generic Newton polygon} on $X$ is upper semicontinuous, with the endpoints being locally constant;
in particular, this function defines a locally closed stratification of $X$.
\end{theorem}
\begin{proof}
By Theorem~\ref{T:companion}, this follows from Lemma~\ref{L:convergent Newton polygon}(b) and (c).
\end{proof}

\begin{theorem}[after de Jong--Oort--Yang] \label{T:purity}
With notation as in Theorem~\ref{T:stratification1}, the Newton polygon stratification jumps purely in codimension $1$. More precisely, for $X$ irreducible,
each breakpoint of the generic Newton polygon disappears purely in codimension $1$.
\end{theorem}
\begin{proof}
By Theorem~\ref{T:companion}, this follows from Lemma~\ref{L:convergent Newton polygon}(c).
\end{proof}

\begin{remark}
Suppose that $X$ admits a good compactification $\overline{X}$ and that $\calE$ is a docile coefficient object on $X$.
Apply Theorem~\ref{T:companion} to construct a crystalline companion $\calF$ of $\calE$; by Corollary~\ref{C:companions are tame}, $\calF$ is again docile. We then define the Newton polygon function
$N_x(\calE) := N_x(\calF)$; by retracing through the arguments cited in \cite[\S 3]{kedlaya-isocrystals},
it can be shown that it satisfies the analogues of the theorems of Grothendieck--Katz and de Jong--Oort--Yang.
In particular, the conclusions of Theorem~\ref{T:stratification1} and Theorem~\ref{T:purity} can be seen to carry over to this definition; we leave a detailed development of this statement to a later occasion.
\end{remark}

\begin{example}
Take $X = \PP^1_k \setminus \{0,1,\infty\}$ with coordinate $\lambda$ and let $\calE$ be the middle cohomology of the Legendre elliptic curve $y^2 = x(x-1)(x-\lambda)$. (This is similar to \cite[Example~4.6]{kedlaya-isocrystals}
except with $N=2$.) Using the Tate uniformization of elliptic curves with split multiplicative reduction, one can show that $\calE$ is docile and for $\lambda \in \{0,1,\infty\}$,  $N_\lambda(\calE)$ equals the generic value (i.e., its slopes are 0 and 1).
\end{example}

We next extend some results of Koshikawa \cite{koshikawa} from crystalline to \'etale coefficients. (This is only meant to be a representative sample to illustrate the technique.)
\begin{defn} \label{D:absolutely unit-root}
Let $\calE$ be an $E$-algebraic coefficient object for some number field $E$, and fix an embedding $E \hookrightarrow \overline{\QQ}_p$. 
We say that $\calE$ is \emph{isoclinic} if for all $x \in X$, $N_x(\calE)$ consists of a single slope with some multiplicity;
if this slope is always 0, we say moreover that $\calE$ is \emph{unit-root}.
It suffices to check this condition at generic points; it implies local constancy of $N_x(\calE)$ on $X$. 

We say that $\calE$ is \emph{absolutely isoclinic/unit-root} if it is isoclinic/unit-root with respect to every choice of the embedding $E \hookrightarrow \overline{\QQ}_p$. If $\calE$ is absolutely isoclinic, then by Kronecker's theorem on roots of unity, we can make $\calE$ absolutely unit-root using a constant twist (at the expense of possibly enlarging the number field generated by Frobenius traces).
By way of contrast, see \cite[Example~2.2]{koshikawa} for an example of a coefficient object which is unit-root but not absolutely isoclinic.
\end{defn}

\begin{defn}
A coefficient object $\calE$ on $X$ is \emph{isotrivial} if there exists a finite \'etale cover $f\colon Y \to X$ such that $f^* \calE$ is trivial. This implies that $\calE$ is absolutely unit-root,
but not conversely \cite[Example~3.5]{koshikawa}. 
By Lemma~\ref{L:companion irreducible}(d), if two semisimple coefficient objects are companions, then one is isotrivial if and only if the other is.
\end{defn}

\begin{theorem}[after Koshikawa] \label{T:koshikawa}
Let $\calE$ be a coefficient object on $X$ which is absolutely unit-root. 
Assume in addition either that $\calE$ is semisimple, or that each irreducible constituent of $\calE$ has determinant of finite order.
Then $\calE$ is isotrivial.
\end{theorem}
\begin{proof}
In both cases, we may reduce to the case where $\calE$ is irreducible with determinant of finite order
(using Lemma~\ref{L:decompose to finite order} in the first case). 
Using Corollary~\ref{C:deligne}, we may reduce to the crystalline case.
We may then apply \cite[Theorem~1.4]{koshikawa} to conclude.
\end{proof}

\begin{theorem}[after Koshikawa] \label{T:koshikawa2}
Suppose that $X$ is proper.
Let $\calE$ be a semisimple, $\QQ$-algebraic coefficient object on $X$ 
with constant Newton polygon.
Then $\calE$ is isotrivial.
\end{theorem}
\begin{proof}
Using Theorem~\ref{T:algebraic companions}, we may reduce to the crystalline case.
We may then apply \cite[Theorem~1.6]{koshikawa} to conclude.
\end{proof}

\begin{remark}
In Theorem~\ref{T:koshikawa2}, \cite[Example~2.2]{koshikawa} shows that the hypothesis of $\QQ$-algebraicity cannot be relaxed to mere algebraicity.
\end{remark}

\subsection{Wan's theorem on fixed-slope \texorpdfstring{$L$}{L}-functions}

We give a statement on a certain factorization of the $L$-function of $\calE$.

\begin{theorem}[after Wan] \label{T:Wan}
Let $\calE$ be an algebraic coefficient object on $X$ and fix a $p$-adic valuation on the full coefficient field of $\calE$.
For $s \in \QQ$, let $P_{s}(\calE_x, T)$ be the factor of $P(\calE_x,T)$ with constant term $1$ corresponding to the slope-$s$ segment of $N_x(\calE)$.
Then the associated L-function 
\[
L_s(X, \calE, T) = \prod_{x \in X^\circ} P_{s}(\calE_x,T^{[\kappa(x):k]})^{-1}
\]
is $p$-adic meromorphic (i.e., a ratio of two $p$-adic entire series).
\end{theorem}
\begin{proof}
For any locally closed stratification of $X$, $L_s(X, \calE, T)$ equals the product of $L_s(Y, \calE, T)$
as $Y$ varies over the strata; we may thus assume that $X$ is affine.
Apply Theorem~\ref{T:companion} to construct a crystalline companion $\calF$ of $\calE$;
we then have $L_s(X,\calE,T) = L_s(X,\calF,T)$.
The $p$-adic meromorphicity of $L_s(X,\calF,T)$ is a theorem of Wan \cite[Theorem~1.1]{wan2}
(see also \cite{wan1, wan3}); this proves the claim.
\end{proof}

\begin{remark}
In the crystalline case, Theorem~\ref{T:Wan} had been conjectured by Dwork \cite{dwork};
this conjecture, originally resolved by the work of Wan cited above, 
was motivated by his original study of zeta functions of algebraic varieties via $p$-adic analytic methods. Indeed, Dwork's original proof of rationality of zeta functions of varieties over finite fields \cite{dwork-zeta1, dwork-zeta2} involved combining archimedean and $p$-adic analytic information about zeta functions as power series.
\end{remark}

\subsection{Skeleton sheaves}

We reformulate \cite[Theorem~2.5]{drinfeld-deligne}, restricted to smooth $k$-schemes, to include $p$-adic coefficients.

\begin{defn}
Let $\calE$ be a coefficient object on $X$. For $x \in X^\circ$, let $c(\calE, x)$ be the number of eigenvalues of Frobenius on $\calE_x$ which belong to $\overline{\QQ}$ (counted with multiplicity);
then $\calE$ is algebraic if and only if $c(\calE,x) = \rank(\calE_x)$ for all $x \in X^\circ$.
\end{defn}

\begin{lemma} \label{L:number of algebraic eigenvalues}
For any coefficient object $\calE$ on $X$, the function $c(\calE, \bullet)\colon X^\circ \to \ZZ$ is locally constant. In particular, if $X$ is irreducible then $c$ is constant.
\end{lemma}
\begin{proof}
Since we can pass a smooth curve through any two given closed points of $X$ in the same component, we may reduce to the case where $X$ is a curve.
We may then furher reduce to the case where $\calE$ is irreducible; in this case, we deduce the claim directly from Corollary~\ref{C:on curve constant twist}.
\end{proof}

\begin{lemma} \label{L:split algebraic eigenvalues1}
Suppose that $X$ is an affine curve.
Let
\[
0 \to \calE_1 \to \calE \to \calE_2 \to 0
\]
be a short exact sequence of coefficient objects on $X$ such that $c(\calE,x) = c(\calE_2,x) = \rank(\calE_{2,x})$ for all $x \in X^\circ$. Then this exact sequence splits.
\end{lemma}
\begin{proof}
We may use internal Homs to reduce to the case $\calE_2 = \calO$,
in which case we must show that $\Ext(\calO, \calE_1) = 0$. For this, we may reduce to the case where $\calE_1$ is irreducible. In this case,
by Corollary~\ref{C:on curve constant twist}, for some transcendental element $\lambda$ of the full coefficient field, 
the constant twist $\calF$ of $\calE_1$ by $\lambda$ has the property that the Frobenius eigenvalues of $\calF_x$ are algebraic for each $x \in X^\circ$.

The eigenvalues of Frobenius on $H^0(X, \calF)$ are all roots of unity and hence algebraic;
since $X$ is an affine curve, we can combine the algebraicity of the eigenvalues of $\calF$ with the Lefschetz trace formula \cite[(1.1.7.1)]{kedlaya-companions}
to deduce that the eigenvalues of Frobenius on $H^1(X, \calF)$ are algebraic. Since $\lambda$ is not algebraic, we deduce that none of the eigenvalues of Frobenius
on $H^0(X, \calE_1)$ or $H^1(X, \calE_1)$ is equal to 1; that is, in the Hochschild--Serre exact sequence
\[
H^0(X, \calE_1)_{\varphi} \to \Ext(\calO, \calE_1) \to H^1(X, \calE_1)^{\varphi}
\]
the two terms on the ends are both zero. We thus deduce that $\Ext(\calO, \calE_1) = 0$ as desired.
\end{proof}

\begin{cor} \label{C:split algebraic eigenvalues}
Suppose that $X$ is an affine curve. For any coefficient object $\calE$ on $X$, there exists a direct sum decomposition
$\calE = \calE_1 \oplus \calE_2$ such that $c(\calE,x) = c_1(\calE_2,x) = \rank(\calE_{2,x})$ for all $x \in X^\circ$.
\end{cor}
\begin{proof}
We proceed by induction on $\rank(\calE)$. If $\calE$ is nonzero, let $\calE_1$ be an irreducible subobject of $\calE$.
By first applying Corollary~\ref{C:on curve constant twist} to $\calE_1$, then applying Lemma~\ref{L:split algebraic eigenvalues1} and the induction hypothesis to $\calE/\calE_1$,
we deduce the claim.
\end{proof}

\begin{defn} \label{D:skeleton sheaf}
Fix a category $\calC$ of coefficient objects with full coefficient field $F$.
A \emph{skeleton sheaf} on $X$ valued in $\calC$ (of rank $n$)
is a function $\chi\colon X^\circ \to F[T]$ such that
for each morphism $f\colon C \to X$ of $k$-schemes with $C$ a curve over $k$,
there exists a coefficient object $\calE_C$ on $C$ in $\calC$ (of rank $n$)
such that $P(\calE_{C,x}, T) = \chi(f(x))$ for each $x \in C^\circ$.

We say that $\chi$ is \emph{tame} (resp. \emph{docile}) if the coefficient objects $\calE_C$ can all be taken to be tame (resp.\ docile).

We say that $\chi$ is \emph{representable} if there exists a coefficient object $\calE$ on $X$ in the specified category such that $P(\calE_x, T) = \chi(x)$ for each $x \in X^\circ$.
\end{defn}

\begin{theorem}
Let $\chi$ be a skeleton sheaf on $X$ such that for some alteration $g\colon X' \to X$,
the pullback of $\chi$ along $g$ is tame. Then $\chi$ is representable.
\end{theorem}
\begin{proof}
We may asume that $X$ is irreducible. The function $x \mapsto \deg(\chi(x))$ on $X^\circ$ is constant on smooth irreducible curves in $X$, and hence is constant;
we induct on this constant value.
For $x \in X^\circ$, let $c(x, \chi)$ be the number of algebraic eigenvalues of $\chi(x)$ (counted with multiplicity);
by Lemma~\ref{L:number of algebraic eigenvalues}, this function is also constant.
Since it is harmless to make a single constant twist on all of the $\calE_C$, we may assume that $c(x, \chi) \neq 0$.

For each $x \in X^\circ$, factor $\chi(x)$ as $\chi_1(x) \chi_2(x)$ where $\chi_1$ has all roots transcendental and $\chi_2$ has all roots algebraic (and $\deg \chi_2 > 0$).
By Lemma~\ref{L:split algebraic eigenvalues1}, we may split $\calE_C$ as a direct sum $\calE_{C,1} \oplus \calE_{C,2}$ such that
$P(\calE_{C,i,x}, T) = \chi_i(x)$ for $i \in \{1,2\}$ and $x \in C^\circ$.
By the induction hypothesis, there exists a coefficient object $\calE_1$ on $X$ in the specified category such that $P(\calE_{1,x}, T) = \chi_1(x)$ for each $x \in X^\circ$.
Meanwhile, after applying Corollary~\ref{C:algebraic companion} if needed to convert $\calE_{C,2}$ into an \'etale coefficient,
we may apply \cite[Th\'eor\`eme~3.1]{deligne-finite} (compare \cite[Th\'eor\`eme~5.1]{cadoret}) to deduce that the polynomials $\chi(x)$ for $x \in X^\circ$ are all defined over a single number field.
We may then apply \cite[Theorem~2.5]{drinfeld-deligne} to  construct an \'etale coefficient $\calE_2$
such that $P(\calE_{2,x}, T) = \chi_2(x)$ for each $x \in X^\circ$.  By applying Theorem~\ref{T:algebraic companions} if needed to convert $\calE_2$ into a crystalline coefficient,
we may take $\calE = \calE_1 \oplus \calE_2$ to conclude.
\end{proof}

\subsection{Reduction of the structure group}
\label{sec:reduction}

As mentioned in the introduction, one can interpret the problem of constructing companions, for coefficient objects of rank $r$, as a problem implicitly associated to the group $\GL_r$; it is then natural to consider the corresponding problem for other groups. A closely related question is the extent to which monodromy groups are preserved by the companion relation. In the \'etale-to-\'etale case, this question was studied in some detail by Chin
\cite{chin}. The adaptations of Chin's work to incorporate the crystalline case were originally made by 
P\'al \cite{pal}; we follow more closely the treatment by D'Addezio \cite{daddezio1}. 

\begin{hypothesis}
Throughout \S\ref{sec:reduction}, we impose \cite[Hypothesis~1.3.1]{kedlaya-companions}: assume that $X$ is connected, fix a closed point $x \in X^\circ$ and a geometric point
$\overline{x} = \Spec(\overline{k})$ of $X$ lying over $x$, and for each category of crystalline coefficient objects fix an embedding of $W(\overline{k})$ into the full coefficient field.
\end{hypothesis}

The following result illustrates a certain limitation of Corollary~\ref{C:deligne prior}; see Remark~\ref{R:base field}.

\begin{theorem} \label{T:daddezio1}
Let $\calE$ be an algebraic coefficient object on $X$. 
Then there exists a number field $E$ with the following property:
for every category of coefficient fields and every embedding of $E$ into the full coefficient field,
for $L$ the corresponding completion of $E$, 
$\calE$ admits an companion in $\Weil(X) \otimes L$ or $\FIsoc^\dagger(X) \otimes L$.
\end{theorem}
\begin{proof}
In the \'etale-to-\'etale case, this is due to Chin \cite[Main~Theorem]{chin};
we may add the crystalline case using Corollary~\ref{C:deligne}.
See also \cite[Theorem~3.7.2]{daddezio1}.
\end{proof}

\begin{remark} \label{R:base field}
The point of Theorem~\ref{T:daddezio1} is that given that $\calE$ is $E$-algebraic, any companion is also $E$-algebraic (this being a condition solely involving Frobenius traces) but may not be realizable over the completion of $E$; see \cite[Remark~3.7.3]{daddezio1} for an illustrative example.
The content of Chin's result is that this phenomenon can be eliminated everywhere at once by a single finite extension of $E$ itself; note that this does not follow from Lemma~\ref{L:bound coefficient extension} because that result does not give enough control on the local extensions.
\end{remark}

We next address the question of independence of $\ell$ in the formation of monodromy groups.
At the level of component groups, one has the following statement.
\begin{theorem} \label{T:same component group}
Let $\calE$ be a pure algebraic coefficient object on $X$ and let $\calF$ be a companion of $\calE$.
Then there exists an isomorphism $\pi_0(G(\calE)) \cong \pi_0(G(\calF))$
which is compatible with the surjections $\psi_\calE, \psi_\calF$ from $\pi_1^{\mathrm{et}}(X, \overline{x})$
(see \cite[Proposition~1.3.11]{kedlaya-companions}),
and which induces an isomorphism $\pi_0(\overline{G}(\calE)) \cong \pi_0(\overline{G}(\calF))$.
\end{theorem}
\begin{proof}
In the \'etale case, this is due to Serre \cite{serre} and Larsen--Pink \cite[Proposition~2.2]{larsen-pink};
the key step is to show that $G(\calE)$ (resp.\ $\overline{G}(\calE)$) is connected if and only if
$G(\calF)$ (resp.\ $\overline{G}(\calF)$) is.
For the adaptation of this result, see \cite[Proposition~8.22]{pal} for the case of $\overline{G}$,
or \cite[Theorem~4.1.1]{daddezio1} for both cases.
\end{proof}

At the level of connected components, one has the following result.

\begin{theorem} \label{T:same connected monodromy}
Let $\calE$ be a semisimple algebraic coefficient object on $X$. 
Then for some number field $E$ as in Theorem~\ref{T:daddezio1},
there exists a connected split reductive group $G_0$ (resp. $\overline{G}_0$) over $E$
such that for every place $\lambda$ of $E$ at which $\calE$ admits a companion $\calF$
(possibly including $\calE$ itself),
the group $G_0 \otimes_E E_\lambda$ is isomorphic to $G(\calF)^\circ$
(resp. the group $\overline{G}_0 \otimes_E E_\lambda$ is isomorphic to $\overline{G}(\calF)^\circ$.
Moreover, these isomorphisms can be chosen so that for some faithful $E$-linear representation $\rho_0$ of $G_0$ (resp. $\overline{\rho}_0$ of $\overline{G}_0$) chosen independently of $\lambda$, the representation $\rho_0 \otimes_E E_\lambda$ (resp. $\overline{\rho}_0 \otimes_E E_\lambda$)
corresponds to the restriction to $G(\calF)^\circ$ of the canonical representation of $G(\calF)$
(resp. the restriction to $\overline{G}(\calF)^\circ$ of the canonical representation of $\overline{G}(\calF)$).
\end{theorem}
\begin{proof}
As per \cite[Remark~4.3.9]{daddezio1}, it suffices to prove the statements about $G_0$ as they immediately imply the corresponding statements about $\overline{G}_0$.
Since we assume $\calE$ is semisimple, we may reduce to the case where it is irreducible;
then $\calE$ is pure by Lemma~\ref{L:irreducible is pure}, and by Corollary~\ref{C:deligne} it becomes $p$-plain after a constant twist.
We may then apply \cite[Theorem~1.4]{chin2} in the \'etale case and \cite[Theorem~8.23]{pal} or \cite[Theorem~4.3.2]{daddezio1} in the crystalline case.
\end{proof}

\begin{remark} \label{R:chin integrated}
It seems likely that one can integrate Theorem~\ref{T:same component group} and Theorem~\ref{T:same connected monodromy} into a single independence statement for arithmetic monodromy groups; see \cite[Conjecture~1.1]{chin}.
Such a result would give an alternate proof of \cite[Corollary~3.7.5]{kedlaya-companions}.

On a related note, it may be possible to treat some new cases of \cite[Conjecture~1.1]{chin} using the automorphic-to-Galois direction of the geometric Langlands correspondence for general reductive groups, which is currently available in the \'etale case
\cite{lafforgue}. If so, then incorporating the crystalline case would require an extension of \cite{lafforgue} to the crystalline case, which is likely to be available soon (see Remark~\ref{R:isocrystals with extra structure} for some background and \cite{kedlaya-xu} for a more recent development).
\end{remark}

\begin{remark} \label{R:same connected geometric monodromy}
As per \cite[Remark~4.3.9]{daddezio1}, Theorem~\ref{T:same connected monodromy} immediately implies the corresponding statement for geometric monodromy groups. 
\end{remark}

In connection with the discussion of Lefschetz slicing in \cite[\S 3.7]{kedlaya-companions}, we also mention the following result.
\begin{theorem} \label{T:preservation of monodromy group}
Let $f\colon Y \to X$ be a morphism of smooth connected $k$-schemes, and fix a consistent choice of base points.
Let $\calE, \calF$ be semisimple coefficient objects on $X$ which are companions. 
Then the inclusion $G(f^* \calE) \to G(\calE)$ (resp.\ $\overline{G}(f^* \calE) \to \overline{G}(\calE)$)
is an isomorphism if and only if $G(f^* \calF) \to G(\calF)$ 
(resp.\ $\overline{G}(f^* \calF) \to \overline{G}(\calF)$)
is an isomorphism.
\end{theorem}
\begin{proof}
By Lemma~\ref{L:monodromy}(a), $\overline{G}(\calE)^\circ$
and $\overline{G}(\calF)^\circ$ are semisimple. 
We may thus apply \cite[Theorem~4.4.2]{daddezio1} to conclude.
\end{proof}

\begin{remark} \label{R:isocrystals with extra structure}
On the crystalline side, there is a rich theory of \emph{isocrystals with additional structure} encoding the replacement of $\GL_r$ by a more general group $G$; 
in particular, this theory provides the analogue of a Newton polygon for a $G$-isocrystal at a point, which carries more information than just the Newton polygon of the underlying isocrystals.
The basic references for this are the papers of
Kottwitz \cite{kottwitz1,kottwitz2} and Rapoport--Richartz 
\cite{rapoport-richartz}. 
\end{remark}

\subsection{Extension along a fibration}

The following result is of interest because a direct proof in the crystalline case would provide an alternate construction of crystalline companions. 

\begin{theorem}
Let $f\colon \overline{X} \to S$ be a smooth curve fibration with pointed locus $Z$ and unpointed locus $X$. 
Choose $x \in S^\circ$ and let $\calE_x$ be an irreducible tame coefficient object on $X \times_S x$. Then there exist an \'etale morphism $S' \to S$, a point $y \in S^{\prime \circ}$ lying over $x$, and a tame coefficient object $\calE$ on $X \times_S S'$ such that the pullbacks of $\calE$ and $\calE_x$ to $X \times_S y$ are isomorphic.
\end{theorem}
\begin{proof}
By Lemma~\ref{L:decompose to finite order}, by making a constant twist we may ensure that $\calE_x$ has determinant of some finite order $n$, 
and hence is algebraic (Lemma~\ref{L:irreducible to algebraic}). 
In this case, we will ensure that $\calE$ is also algebraic.
By Theorem~\ref{T:companion}, we need only treat the \'etale case. 

Given the shape of the desired conclusion, we are free at any point to pull everything back along an \'etale morphism $S' \to S$ whose image contains $x$. To begin with, we may enforce that $X \to S$ admits a section $s$. Let $\overline{x}$ be a geometric point over $x$; we then have a sequence \cite[Expos\'e XIII, Proposition 4.1]{grothendieck-sga1-tame}
\begin{equation} \label{eq:tame sequence1}
1 \dashrightarrow \pi_1^{\tame}(X \times_S \overline{x}) \to \pi_1^{\tame}(X) \to \pi_1^{\tame}(S) \to 1
\end{equation}
which is split by $s$; we thus obtain a full short exact sequence (including the dashed arrow).

The group $\pi_1^{\tame}(S)$ itself fits into an exact sequence
\[
1 \to \pi_1^{\tame}(S \times_k \overline{k}) \to \pi_1^{\tame}(S) \to \pi_1(k) \to 1;
\]
we may assume that $x \cong \Spec k$, in which case this sequence is itself split by the inclusion $x \to S$.

The outer action of $\pi_1^{\tame}(S)$ on
$\pi_1^{\tame}(X \times_S \overline{x})$ via \eqref{eq:tame sequence1} induces an action on the set $T$ of isomorphism classes of continuous irreducible $\overline{\QQ}_\ell$-representations of $\pi_1^{\tame}(X \times_S \overline{x})$ of rank $r$ with determinant of finite order $n$.
Then the \'etale coefficient $\calE_x$ corresponds to an isomorphism class $\overline{\rho}$ in $T$ which is fixed by the 
action of $\pi_1(k) \subset \pi_1^{\tame}(S)$. 
Any isomorphism class in the orbit of $\overline{\rho}$ under $\pi_1^{\tame}(S)$ is also fixed by the action of some copy of $\pi_1(k)$ in $ \pi_1^{\tame}(S)$, and thus extends to a representation of $\pi_1^{\tame}(X \times_S \overline{x})$ with determinant of order $n$.
By Deligne's finiteness theorem (e.g., see \cite[Theorem~2.1]{esnault-kerz} or \cite[Corollary~3.7.6]{kedlaya-companions}), we deduce that the orbit of $\overline{\rho}$ in $T$ is in fact a finite set.

From the previous paragraph, by making an \'etale base change on $S$, we may ensure that the isomorphism class of $\overline{\rho}$ is fixed by the outer action of $\pi_1^{\tame}(S)$ via \eqref{eq:tame sequence1}. In the case $r=1$, using the splitting via $s$ we immediately obtain an extension of $\overline{\rho}$ to $\pi_1^{\tame}(X)$. 

In the general case, we know from the previous paragraph that $\det(\overline{\rho})$ extends to $\pi_1^{\tame}(X)$.
Hence the remaining obstruction to extending $\overline{\rho}$ is the Clifford obstruction valued in the group $H^2(\pi_1^{\tame}(S), \mu_r)$. This can also be trivialized by making an \'etale base  change on $S$.
\end{proof}


\begin{thebibliography}{999}
\bibitem{abe-companion}
T. Abe, Langlands correspondence for isocrystals and existence of
crystalline companion for curves, \textit{J. Amer. Math. Soc.} \textbf{31} (2018), 921--1057.

\bibitem{abe-esnault}
T. Abe and H. Esnault, A Lefschetz theorem for overconvergent isocrystals with Frobenius
structure, \textit{Ann. Sci. Éc. Norm. Supér.} \textbf{52} (2019), 1243--1264. 

\bibitem{agrawal}
S. Agrawal, Deformations of overconvergent isocrystals on the projective line,
\textit{J. Number Theory} \textbf{237} (2022), 167--241.

\bibitem{arabia}
A. Arabia, Rel\`evements des alg\`ebres lisses et de leurs morphismes,
\textit{Comment. Math. Helv.} \textbf{76} (2001), 607--639.

\bibitem{artin-sga4}
M. Artin, Comparaison avec la cohomologie classique: cas d'un schema lisse,
Expos\'e XI in \textit{Th\'eorie des Topos et Cohomologie \'Etale des Schemas (SGA 4 III)}, Lecture Notes in Math. 305, Springer-Verlag, Berlin, 1973.

\bibitem{atiyah}
M. Atiyah, Complex analytic connections in fibre bundles, \textit{Trans. Amer. Math. Soc.} \textbf{85} (1957), 181--207.

\bibitem{berthelot-mem}
P. Berthelot, G\'eom\'etrie rigide et cohomologie des vari\'et\'es alg\'ebriques de caract\'eristique $p$, in
Introductions aux cohomologies $p$-adiques (Luminy, 1982), 
\textit{M\'em. Soc. Math. France (N.S.)} \textbf{23} (1986), 7--32.

\bibitem{bhatt-lurie}
B. Bhatt, Crystals and Chern classes, arXiv:2310.00676v1 (2023).

\bibitem{biswas-heu}
I. Biswas and V. Heu, Meromorphic connections on vector bundles over curves,
\textit{Proc. Indian Math. Soc.} \textbf{124} (2014), 487--496.

\bibitem{bosch}
S. Bosch, \textit{Lectures on Formal and Rigid Geometry}, 
Lecture Notes in Math. 2105, Springer, Cham, 2014.

\bibitem{bgr}
S. Bosch, U. G\"untzer, and R. Remmert,
\textit{Non-Archimedean Analysis},
Grundlehren der Math. Wiss. 261, Springer-Verlag, Berlin, 1984.

\bibitem{cadoret}
A. Cadoret, La conjecture des compagnons d'apr\`es Deligne, Drinfeld, L. Lafforgue, T. Abe., ...,
Expos\'e Bourbaki 1155, \textit{Ast\'erisque} \textbf{422} (2020), 173--223.

\bibitem{chin}
C. Chin, Independence of $\ell$ in Lafforgue's theorem, \textit{Adv. Math.} \textbf{180} (2003), 64--86.

\bibitem{chin2}
C. Chin, Independence of $\ell$ of monodromy groups, \textit{J. Amer. Math. Soc.} \textbf{17} (2004), 723--747.

\bibitem{conrad}
B. Conrad, Relative ampleness in rigid geometry, \textit{Ann. Inst. Fourier} \textbf{56} (2006), 1049--1126.

\bibitem{crew-mono}
R. Crew, $F$-isocrystals and their monodromy groups,
\textit{Ann. Scient. \'Ec. Norm. Sup.} \textbf{25} (1992), 429--464.

\bibitem{daddezio1}
M. D'Addezio, The monodromy groups of lisse sheaves and overconvergent F-isocrystals,
\textit{Sel. Math.} \textbf{26} (2020).

\bibitem{daddezio}
M. D'Addezio, Parabolicity conjecture of $F$-isocrystals, \textit{Annals of Math.} \textbf{198} (2023), 619--656.

\bibitem{dejong-alterations}
A.J. de Jong, Smoothness, semi-stability, and alterations,
\textit{Publ. Math. IH\'ES} \textbf{83} (1996), 51--93.

\bibitem{deligne-weil2}
P. Deligne, La conjecture de Weil, II, \textit{Publ. Math. IH\'ES} \textbf{52} (1981), 313--428.

\bibitem{deligne-finite}
P. Deligne, Finitude de l'extension de $\mathbb{Q}$ engendr\'ee par des traces 
de Frobenius, en caract\'eristique finie, \textit{Moscow Math. J.} \textbf{12} (2012),
497--514.

\bibitem{drinfeld}
V.G. Drinfeld,
Langlands’ Conjecture for $\GL(2)$ over functional Fields, 
Proceedings of the International Congress of Mathematicians, Helsinki, 565--574 (1978).

\bibitem{drinfeld-deligne}
V. Drinfeld, On a conjecture of Deligne, \textit{Moscow Math. J.} \textbf{12} (2012),
515--542, 668.

\bibitem{drinfeld-pro-simple}
V. Drinfeld, On the pro-semisimple completion of the fundamental group of a smooth variety
over a finite field, \textit{Adv. Math.} \textbf{327} (2018), 708--788.

\bibitem{drinfeld-kedlaya}
V. Drinfeld and K.S. Kedlaya, Slopes of indecomposable $F$-isocrystals,
\textit{Pure Appl. Math. Quart.} \textbf{13} (2017), 131--192.

\bibitem{dwork-zeta1}
B. Dwork, On the rationality of the zeta function of an algebraic variety, \textit{Amer. J. Math.} \textbf{82} (1960), 631--648.

\bibitem{dwork-zeta2}
B. Dwork, On the zeta function of a hypersurface, \textit{Publ. Math. IH\'ES} \textbf{12} (1962), 5--68.

\bibitem{dwork}
B. Dwork, Normalized period matrices II, \textit{Ann. of Math.} \textbf{98} (1973), 1--57.

\bibitem{elkik}
R. Elkik, Solutions d'\'equations \`a coefficients dans un anneau hens\'elien, \textit{Ann. Scient.
\'Ec. Norm. Sup.} \textbf{6} (1973), 553--604.

\bibitem{esnault}
H. Esnault, \textit{Local Systems in Algebraic-Arithmetic Geometry}, 
Lecture Notes in Math. 2337, Springer Nature, 2023.

\bibitem{esnault-groechenig}
H. Esnault and M. Groechenig, Crystallinity of rigid flat connections revisited, arXiv:2309.15949v2 (2023).

\bibitem{esnault-kerz}
H. Esnault and M. Kerz, A finiteness theorem for Galois representations of function fields over finite fields
(after Deligne), \textit{Acta Math. Vietnam.} \textbf{37} (2012), 531--562.

\bibitem{esnault-shiho2}
H. Esnault and A. Shiho, Chern classes of crystals, \textit{Trans. Amer. Math. Soc.}
\textbf{371} (2019), 1333--1358.

\bibitem{gabber-ramero}
O. Gabber and L. Ramero,
\textit{Almost Ring Theory}, Lecture Notes in Math. 1800,
Springer-Verlag, Berlin, 2003.

\bibitem{gieseker}
D. Gieseker, Stable vector bundles and the Frobenius morphism, \textit{Ann. Scient. \'Ec. Norm. Sup.}
\textbf{6} (1973), 95--101.

\bibitem{grosse-klonne}
E. Grosse-Kl\"onne, Rigid analytic spaces with overconvergent structure sheaf,
\textit{J. reine angew. Math.} \textbf{519} (2000), 73--95.

\bibitem{grothendieck-sga1-tame}
A. Grothendieck and M. Raynaud, \textit{Rev\^etements \'Etales et Groupe Fondamental (SGA 1)}, Lecture Notes in Math.\ 224,
Springer-Verlag, Berlin, 1971.

\bibitem{grothendieck-ogg-shafarevich}
A. Grothendieck, Formule d'Euler--Poincar\'e en cohomologie \'etale, Expos\'e X in \textit{S\'eminaire de G\'eom\'etrie Alg\'ebriques du Bois--Marie (SGA 5)}, Lecture Notes in Math.\ 589,
Springer, Berlin, 1977.

\bibitem{gku}
T. Grubb, K.S. Kedlaya, and J. Upton, 
A cut-by-curves criterion for overconvergence of $F$-isocrystals,
arXiv:2202.03604v1 (2022). 

\bibitem{hartshorne}
R. Hartshorne, Ample vector bundles on curves, \textit{Nagoya Math. J.} \textbf{43} (1971), 73--89.

\bibitem{kato-log}
K. Kato, Logarithmic structures of Fontaine--Illusie, in \textit{Algebraic Analysis, Geometry, and Number Theory},
Johns Hopkins University, Baltimore, 1988, 191--224.

\bibitem{katz}
N.M. Katz,
Nilpotent connections and the monodromy theorem: applications of a result of Turrittin,
\textit{Publ. Math. IH\'ES} \textbf{39} (1970), 175--232.

\bibitem{katz-slope}
N. Katz, Slope filtrations of $F$-crystals, Journ\'ees de G\'eom\'etrie Alg\'ebriques
(Rennes, 1978), \textit{Ast\'erisque} \textbf{63} (1979), 113--164.

\bibitem{kedlaya-semi1}
K.S. Kedlaya, Semistable reduction for overconvergent $F$-isocrystals, I: Unipotence and logarithmic extensions, \textit{Compos. Math.} 143 (2007), 1164--1212;
errata, \cite{shiho-log}.

\bibitem{kedlaya-semi2}
K.S. Kedlaya, Semistable reduction for overconvergent $F$-isocrystals, II: A valuation-theoretic approach, \textit{Compos. Math.} \textbf{144} (2008), 657--672.

\bibitem{kedlaya-semi3}
K.S. Kedlaya, Semistable reduction for overconvergent $F$-isocrystals, III: local semistable reduction at monomial valuations, \textit{Compos. Math.}
\textbf{145} (2009), 143--172.

\bibitem{kedlaya-ainf}
K.S. Kedlaya, Some ring-theoretic properties of $\mathbf{A}_{\mathrm{inf}}$, 
\textit{$p$-adic Hodge Theory}, Simons Symposia, Springer, 2020, 129--141.

\bibitem{kedlaya-semi4}
K.S. Kedlaya, Semistable reduction for overconvergent F-isocrystals, 
IV: Local semistable reduction at nonmonomial valuations, 
\textit{Compos. Math.} \textbf{147} (2011), 467--523;
erratum, \cite[Remark~4.6.7]{kedlaya-connections}.

\bibitem{kedlaya-connections}
K.S. Kedlaya,
Local and global structure of connections on nonarchimedean curves, 
\textit{Compos. Math.} \textbf{151} (2015), 1096--1156;
errata (with Atsushi Shiho), \textit{ibid.} 153 (2017), 2658--2665. 

\bibitem{kedlaya-aws}
K.S. Kedlaya, Sheaves, stacks, and shtukas, lecture notes from the 2017 Arizona Winter School: Perfectoid Spaces,
Math. Surveys and Monographs 242, Amer. Math. Soc., 2019. 

\bibitem{kedlaya-course}
K.S. Kedlaya, \textit{$p$-adic Differential Equations}, second edition, Cambridge Univ. Press, Cambridge, 2022.

\bibitem{kedlaya-isocrystals}
K.S. Kedlaya, Notes on isocrystals, \textit{J. Num. Theory} \textbf{237} (2022), 353--394.

\bibitem{kedlaya-companions}
K.S. Kedlaya, \'Etale and crystalline companions, I, 
\textit{\'Epijournal G\'eom. Alg.} \textbf{6} (2022), article no. 20. 

\bibitem{kedlaya-drinfeld-isoc}
K.S. Kedlaya, Drinfeld's lemma for $F$-isocrystals, I,
\textit{Intl. Math. Res. Notices} (2024), article ID rnae039. 
	
\bibitem{kedlaya-liu}
K.S. Kedlaya and R. Liu, Relative $p$-adic Hodge theory: Foundations, \textit{Ast\'erisque} \textbf{371} (2015).

\bibitem{kedlaya-xu}
K.S. Kedlaya and D. Xu, Drinfeld's lemma for F-isocrystals, II: Tannakian approach, 
\textit{Compos. Math.} \textbf{160} (2024), 90--119.

\bibitem{koshikawa}
T. Koshikawa, Overconvergent unit-root $F$-isocrystals and isotriviality, 
\textit{Math. Res. Lett.} \textbf{24} (2017), 1707--1727.

\bibitem{kottwitz1}
R.E. Kottwitz, Isocrystals with additional structure, \textit{Compos. Math.} \textbf{56} (1985), 201--220.

\bibitem{kottwitz2}
R.E. Kottwitz, Isocrystals with additional structure, II, \textit{Compos. Math.} \textbf{109} (1997), 255--339.

\bibitem{krishnamoorthy-pal}
R. Krishnamoorthy and A. P\'al, Rank 2 local systems and abelian varieties, \textit{Sel. Math.} \textbf{27} (2021), paper no. 51.

\bibitem{krishnamoorthy-pal2}
R. Krishnamoorthy and A. P\'al, Rank 2 local systems and abelian varieties II,
\textit{Compos. Math.} \textbf{158} (2022), 868--892.

\bibitem{lafforgue}
L. Lafforgue, Chtoucas de Drinfeld et correspondance de Langlands,
\textit{Invent. Math.} \textbf{147} (2002), 1--241.

\bibitem{lafforgue-hecke}
V. Lafforgue,
Estim\'ees pour les valuations $p$-adiques des valeurs propres des op\'erateurs de Hecke,
\textit{Bull. Soc. Math. France} \textbf{139} (2011), 455--477. 

\bibitem{lafforgue-v}
V. Lafforgue, Chtoucas pour les groupes r\'eductifs et param\'etrisation de Langlands globale,
\textit{J. Amer. Math. Soc.} \textbf{31} (2018), 719--891.

\bibitem{lang-weil}
S. Lang and A. Weil, Number of points of varieties in finite fields, \textit{Amer. J. Math.} \textbf{76} (1954), 819–-827.

\bibitem{langton}
S. Langton, Valuative criteria for families of vector bundles on algebraic varieties,
\textit{Ann. Math.} \textbf{101} (1975), 88--110.

\bibitem{larsen-pink}
M. Larsen and R. Pink, Abelian varieties, $\ell$-adic representations, and $\ell$-independence,
\textit{Math. Ann.} \textbf{302} (1995), 561--580.

\bibitem{laumon-moret-bailly}
G. Laumon and L. Moret-Bailly, \textit{Champs Alg\'ebriques}, Springer, Berlin, 2000.

\bibitem{lieblich}
M. Lieblich, Remarks on the stack of coherent algebras, \textit{Int. Math. Res. Notices} (2006), article ID
75273.

\bibitem{matsumura}
H. Matsumura, \textit{Commutative Algebra}, second edition, Benjamin/Cummings, Reading, Mass., 1980.

\bibitem{mehta-pauly}
V.B. Mehta and C. Pauly, Semistability of Frobenius direct images over curves,
\textit{Bull. Soc. Math. France} \textbf{135} (2007), 105--117.

\bibitem{meredith}
D. Meredith, Weak formal schemes, \textit{Nagoya Math. J.} \textbf{45} (1971), 1--38.

\bibitem{milne}
J.S. Milne, \textit{\'Etale Cohomology}, Princeton Univ. Press, Princeton, 1980.

\bibitem{narasimhan-seshadri}
M.S. Narasimhan and C.S. Seshadri, Stable and unitary vector bundles on a compact Riemann surface,
\textit{Ann. of Math.} \textbf{82} (1965), 540--567.

\bibitem{oh-shimizu}
G. Oh and K. Shimizu, Moduli stacks of crystals and isocrystals, arXiv:2504.14801v1 (2025).

\bibitem{ohtsuki}
M. Ohtsuki, A residue formula for Chern classes associated with logarithmic connections,
\textit{Tokyo J. Math.} \textbf{5} (1982), 13--21.

\bibitem{pal}
A. P\'al, The $p$-adic monodromy group of abelian varieties over global function fields of characteristic $p$, \textit{Doc. Math.} \textbf{27} (2022), 1509--1579.

\bibitem{rapoport-richartz}
M. Rapoport and M. Richartz, On the classification and specialization of $F$-isocrystals with additional structure,
\textit{Compos. Math.} \textbf{103} (1996), 153--181.

\bibitem{schiffmann}
O. Schiffmann, Indecomposable vector bundles and stable Higgs bundles over
smooth projective curves, \textit{Ann. of Math.} \textbf{183} (2016), 297--362.

\bibitem{serre-gaga}
J.-P. Serre, G\'eom\'etrie alg\'ebrique et g\'eom\'etrie analytique,
\textit{Ann. Inst. Fourier} \textbf{6} (1956), 1--42.

\bibitem{serre-galois}
J.-P. Serre, \textit{Galois Cohomology}, corrected second printing, Springer-Verlag, Berlin, 1997.

\bibitem{serre}
J.-P. Serre, Lettres \`a Ken Ribet du 1/1/1981 et du 29/1/1981, in \textit{\OE uvres -- Collected Papers IV},
Springer--Verlag, Heidelberg, 2000, 1--12.

\bibitem{shatz}
S. Shatz, The decomposition and specialization of algebraic families of vector bundles, \textit{Compos. Math.} \textbf{35} (1977), 163--187.

\bibitem{shepherd-barron}
N. Shepherd-Barron, Semi-stability and reduction mod $p$, \textit{Topology}
\textbf{37} (1998), 659--664.

\bibitem{shiho-log}
A. Shiho, On logarithmic extension of overconvergent isocrystals,
\textit{Math. Ann.} \textbf{348} (2010), 467--512.

\bibitem{shiho-cut-over}
A. Shiho, Cut-by-curves criterion for the overconvergence
of $p$-adic differential equations, \textit{Manuscripta Math.}
\textbf{132} (2010), 517--537.

\bibitem{silverman}
J.H. Silverman, \textit{The Arithmetic of Elliptic Curves}, second edition,
Graduate Texts in Math. 106, Springer, Dordrecht, 2009.

\bibitem{simpson1}
C.T. Simpson, Moduli of representations of the fundamental group of a smooth
projective variety, I, \textit{Publ. Math. IH\'ES} \textbf{79} (1994), 47--129.

\bibitem{simpson2}
C.T. Simpson, Moduli of representations of the fundamental group of a smooth projective variety,
II, \textit{Publ. Math. IH\'ES} \textbf{80} (1994), 5--79.

\bibitem{stacks-project}
The Stacks Project Authors, \textit{Stacks Project},
\url{http://stacks.math.columbia.edu}, retrieved March 2026.

\bibitem{sun}
X. Sun, Direct images of bundles under Frobenius morphism, \textit{Invent. Math.} \textbf{173} (2008), 427--447.

\bibitem{tsuzuki}
N. Tsuzuki, Constancy of Newton polygons of $F$-isocrystals on Abelian varieties and isotriviality of families of curves, 
\textit{J. Inst. Math. Jussieu} \textbf{20} (2021), 587--625.

\bibitem{tsuzuki-minimal}
N. Tsuzuki,
Minimal slope conjecture of $F$-isocrystals,
\textit{Invent. Math.} \textbf{205} (2022).

\bibitem{wan1}
D. Wan, Dwork's conjecture on unit root zeta functions, \textit{Ann. Math.} \textbf{150} (1999), 867--927.

\bibitem{wan2}
D. Wan, Higher rank case of Dwork's conjecture,  \textit{J. Amer. Math. Soc.} \textbf{13} (2000), 
807--852.

\bibitem{wan3}
D. Wan, Rank one case of Dwork's conjecture, \textit{J. Amer. Math. Soc.} \textbf{13} (2000), 853--908.

\bibitem{wiesend1}
G. Wiesend, A construction of covers of arithmetic schemes,
\textit{J. Number Theory} \textbf{121} (2006), 118--131.


\end{thebibliography}
\end{document}